\documentclass[onefignum,onetabnum]{siamart251216}

\usepackage{amssymb}
\usepackage{bm}
\usepackage{booktabs}
\usepackage{makecell}
\usepackage{multirow}
\usepackage[export]{adjustbox}
\usepackage{tikz}
\usepackage{pgfplots}
\usepgfplotslibrary{fillbetween}
\usepackage{subfig}
\pgfplotsset{compat=1.18}

\newcommand{\TheTitle}{Computing statistical solutions of a Mach 2000 astrophysical jet}

\title{\TheTitle
\thanks{\funding{The work of the first author was supported by the PRIME programme of the German Academic Exchange Service (DAAD), with funding from the Federal Ministry of Research, Technology and Space. The initiation of this work at UZH was supported by a KHYS ConYS grant at KIT in 2024.}}}

\author{Stephan Simonis\thanks{Seminar for Applied Mathematics, ETH Zurich, 8092 Zurich, Switzerland, and Institute for Applied and Numerical Mathematics, Karlsruhe Institute of Technology, 76131 Karlsruhe, Germany (\email{ssimonis@ethz.ch}, \email{stephan.simonis@kit.edu}).}
\and Gauthier Wissocq\thanks{CEA CESTA, BP 2, 33114 Le Barp, France.}}

\headers{Computing statistical solutions of a Mach 2000 astrophysical jet}{S. Simonis and G. Wissocq}

\ifpdf
\hypersetup{
  pdftitle={\TheTitle},
  pdfauthor={S. Simonis and G. Wissocq}
}
\fi

\begin{document}

\maketitle

\begin{abstract}
The multi-dimensional compressible Euler equations admit non-unique entropy solutions in turbulent regimes, and extreme-Mach astrophysical flows are a natural setting in which this breakdown of deterministic well-posedness becomes computationally visible. We compute statistical solutions of a Mach~2000 astrophysical jet, defined as the pushforward of an initial probability measure through a vectorial lattice Boltzmann method, by Monte Carlo sampling with $M=1000$ realizations on grids of up to $12.8$ million cells. Under mesh refinement the individual realizations diverge pathwise, while the statistical solution converges: Wasserstein distances of the one- and two-point marginals, the ensemble mean, and the ensemble standard deviation all exhibit stable positive convergence rates. A spatially resolved analysis along the jet axis traces this dichotomy to the structure of the one-point laws, which are numerically Dirac in the undisturbed core, skewed in the sheared turbulent regions, and intermittent two-state mixtures at the random leading front. We conclude that the computed statistical solution is non-Dirac and remains stable in the extreme compressible regime, in which no strong solution is expected to exist.
\end{abstract}




\section{Introduction}
Astrophysical jets at Mach numbers exceeding $10^3$ represent a significant challenge in computational fluid dynamics.
The Mach~2000 parameter regime serves two purposes.
First, it is directly physically representative of protostellar jets (Herbig--Haro objects) propagating into ultra-cold ($T \sim 10\,\text{K}$) molecular clouds, where ambient sound speeds drop below $0.3\,\text{km/s}$, yielding classical Mach numbers in excess of $10^3$ \cite{bally2016}.
Second, it serves as a canonical mathematical benchmark for the high-energy limits of the classical Euler equations.
By pushing the deterministic breakdown to the extreme compressible limit, we provide a mathematical proxy for the shock-dominated, hyper-turbulent topologies found in more complex relativistic flows, such as active galactic nuclei and gamma-ray bursts \cite{piran2004, blandford2019}.

While traditional numerical research utilizes the Mach~2000 jet primarily as a short-duration stress test to validate positivity-preserving limiters \cite{zhang2010, wu2023}, the underlying mathematical question of uniqueness in the long-time, fully turbulent regime remains challenging.
State-of-the-art high-order frameworks, such as subcell-limited discontinuous Galerkin methods, typically restrict this benchmark to early transient phases (e.g., $t \le 0.0015$) \cite{rueda2022, rueda2024}.
Kinetic-theory-derived maximum-principle bounds have recently been used to stabilize deterministic simulations of astrophysical jets up to Mach 800 \cite{dzanic2024}.
In this chaotic regime, foundational convex integration results \cite{delellis} demonstrate that the multi-dimensional Euler equations may admit infinitely many entropy solutions for the exact same initial data.
To recover a well-posed theoretical framework, DiPerna introduced measure-valued solutions \cite{diperna1985}, which were subsequently refined into statistical solutions \cite{fjordholm}.
These shift the focus from individual deterministic trajectories to a time-evolving probability measure.
Although conceptually promising, the computation of such statistical solutions for two- and three-dimensional fluid flows remains challenging due to the large number of samples required for approximating the underlying Young measures.
The computation of measure-valued solutions for the incompressible Euler equations was initiated by Lanthaler and Mishra~\cite{lanthaler2015}, and the first results of computing Wasserstein-convergent statistical solutions to three-dimensional compressible flows have been obtained by Lye~\cite{lye2020}.
The use of efficient lattice Boltzmann method (LBM) implementations for approximating statistical solutions of fluid flows was first established by Simonis and Mishra~\cite{simonis2024}.
Nonetheless, to date, the LBM has been applied to compute statistical solutions of incompressible fluid flows only.

In this paper, we transition the Mach~2000 astrophysical jet problem from a deterministic stability test into a study of statistical hyper-compressible turbulence.
To investigate the statistical convergence, we propose a probabilistic LBM that combines Monte Carlo sampling with the vectorial LBM (VLBM) for compressible Euler equations proposed by Wissocq et al.~\cite{wissocq2024}.
We extend the temporal horizon to $t=0.0035$ and compute an ensemble of \(M=1000\) samples on highly resolved spatial grids of up to $4000 \times 3200$ cells, using optimized CUDA kernels for the solver and memory-mapped CPU streaming for the post-processing of the empirical measures.
We contrast the sample-wise strong $L^1$ error with the convergence of the probability measure in Wasserstein metrics, measured both between successive resolutions (\textit{Cauchy} rates) and against the finest grid (\textit{reference} rates).
The strong error ceases to converge under refinement, its \textit{Cauchy} rate dropping to $\approx 0.03$ at the final time, whereas the statistical observables converge: against the finest grid we measure rates of $\approx 0.67$ in the one-point Wasserstein distance, $\approx 0.66$ in the mean, $\approx 0.61$ in the standard deviation, and $\approx 1.03$ in the two-point Wasserstein distance.
A spatially resolved analysis of the one-point law along the jet axis locates the origin of this dichotomy: the marginals are numerically Dirac in the undisturbed core, skewed in the Kelvin--Helmholtz-sheared turbulent regions, and intermittent two-state mixtures at the leading fringe, governed by the random front position, yet all of these laws remain stable under refinement even as the individual realizations sampling them decorrelate.
In this sense we numerically probe the weak--strong dichotomy for a flow in which no strong solution is expected to exist: the scheme selects, through its numerical dissipation, a particular sequence of approximations whose pathwise limit does not exist but whose statistical limit is stable and non-degenerate, while pointwise deterministic convergence is lost.

The remainder of this paper is structured as follows.
Section~\ref{sec:methodology} introduces the numerical methodology: the vectorial lattice Boltzmann scheme for the compressible Euler equations (Section~\ref{sec:methods}), the setup of the Mach~2000 jet including the stochastically perturbed inlet (Section~\ref{sec:setup}), and the notion of statistical solutions together with the construction of the initial probability measure (Section~\ref{sec:statInit}).
Section~\ref{sec:weak-strong} develops the convergence framework: the sampling schedule (Section~\ref{subsec:sampling}), the restriction of all fields to a common grid (Section~\ref{subsec:restriction}), the \textit{Cauchy} and \textit{reference} comparison operators (Section~\ref{subsec:comparison}), and the statistical observables (Section~\ref{subsec:observables}).
Section~\ref{sec:results} presents the numerical results: the implementation, the visual evolution of the density fields (Section~\ref{subsec:visualization}), the spatially resolved analysis along the jet axis (Sections~\ref{subsec:centerline} and~\ref{subsec:pdfheatmap}), the convergence-rate study contrasting pathwise divergence with statistical convergence (Section~\ref{subsec:discussion}), and the temporal trends of the fitted rates (Section~\ref{subsec:ratetrends}).
Section~\ref{sec:conclusions} concludes the paper and suggests further research directions. 

\section{Numerical methodology}\label{sec:methodology}

\subsection{Formulation of the vectorial lattice Boltzmann method}\label{sec:methods}

The VLBM scheme from \cite{wissocq2024} approximates two-dimensional systems of hyperbolic conservation laws using a five-wave ($D2Q5$) kinetic advection--relaxation model.
The method relies on a set of distribution functions $\bm{u}_k$, whose evolution naturally splits into a collision phase and a streaming phase, upholding the general structure of an LBM.
We briefly summarize the deterministic scheme below and refer the reader to \cite{wissocq2024} for further details.

Our main target is the set of Euler equations for compressible flows, i.e.,
\begin{equation}\label{eq:euler}
    \partial_t \bm{u}(\bm{x},t) + \partial_x \bm{f}(\bm{u}(\bm{x},t)) + \partial_y \bm{g}(\bm{u}(\bm{x},t)) = \bm{0} , \qquad \bm{x} \in \mathcal{D} \subseteq \mathbb{R}^{d}, t \in (0, T],
\end{equation}
where the vector of conserved variables is $\bm{u} = (\rho, \rho v_1, \rho v_2, E)^\top$.
The physical fluxes in \eqref{eq:euler} are respectively,
\begin{align}
    \bm{f}(\bm{u}) &= (\rho v_1, \rho v_1^2 + p, \rho v_1 v_2, (E+p)v_1)^\mathrm{T}, \label{eq:fluxF} \\
    \bm{g}(\bm{u}) &= (\rho v_2, \rho v_1 v_2, \rho v_2^2 + p, (E+p)v_2)^\mathrm{T} . \label{eq:fluxG}
\end{align}
The pressure $p$ is related to the conserved variables by the ideal gas equation of state
\begin{equation}\label{eq:eos}
p = (\gamma-1)\left(E - \frac{1}{2}\rho(v_1^2 + v_2^2)\right),
\end{equation}
where $\gamma > 1$ is the constant heat capacity ratio.

The approximation of the macroscopic conserved variables in \eqref{eq:euler} is defined as the zeroth-order moment of the distribution functions, i.e.,
\begin{equation}
\bm{u}(\bm{x}, t) \approx \bm{u}^\epsilon(\bm{x},t) = \sum_{k=1}^{5} \bm{u}_k^\epsilon(\bm{x},t) .
\end{equation}
The collision operator drives the distributions toward a local equilibrium state defined by a set of Maxwellian functions, $\bm{M}_k(\bm{u}^\epsilon)$, designed to satisfy consistency conditions with the target hyperbolic fluxes $\bm{f}(\bm{u}^\epsilon)$ and $\bm{g}(\bm{u}^\epsilon)$, \eqref{eq:fluxF} and \eqref{eq:fluxG}, respectively.
For the $D2Q5$ lattice, the Maxwellians are defined as
\begin{align}
\bm{M}_{1,2} (\bm{u}^\epsilon) &= \left(\frac{1-\alpha}{4}\right)\bm{u}^\epsilon \pm \frac{\bm{f}(\bm{u}^\epsilon)}{2a} ,  \label{eq:maxwellians12} \\
\bm{M}_{3,4}(\bm{u}^\epsilon) &= \left(\frac{1-\alpha}{4}\right)\bm{u}^\epsilon \pm \frac{\bm{g}(\bm{u}^\epsilon)}{2a} , \label{eq:maxwellians34} \\
\bm{M}_5(\bm{u}^\epsilon) &= \alpha \bm{u}^\epsilon , \label{eq:maxwellians5}
\end{align}
where $a$ is the kinetic speed and $\alpha \in [0,1]$ is a free parameter.
To satisfy the Bouchut criterion \cite{bouchut1999}, the kinetic speed must bound the spectral radius of the flux Jacobian matrices.
The continuous kinetic equations rely on a relaxation parameter $\epsilon$.
The numerical schemes are derived by integrating along characteristics and taking the stiff relaxation limit $\epsilon \to 0$.
To capture strong discontinuities while avoiding the severe numerical oscillations typical of second-order schemes, a local flux-limiting blending parameter $\theta_{i,j}^n \in [0,1]$ is introduced.
This parameter effectively bridges the highly dissipative first-order update ($\theta=0$) and the oscillatory second-order update ($\theta=1$).
The resulting blended stream--collide algorithm reads
\begin{align}
\bm{u}_{1,i,j}^{n+1} &= \bm{M}_1(\bm{u}_{i-1,j}^{n}) + \theta_{i-1/2,j}^n[\bm{M}_1(\bm{u}_{i-1,j}^{n}) - \bm{u}_{1,i-1,j}^{n}] , \label{eq:vlbm1} \\
\bm{u}_{2,i,j}^{n+1} &= \bm{M}_2(\bm{u}_{i+1,j}^{n}) + \theta_{i+1/2,j}^n[\bm{M}_2(\bm{u}_{i+1,j}^{n}) - \bm{u}_{2,i+1,j}^{n}] , \label{eq:vlbm2} \\
\bm{u}_{3,i,j}^{n+1} &= \bm{M}_3(\bm{u}_{i,j-1}^{n}) + \theta_{i,j-1/2}^n[\bm{M}_3(\bm{u}_{i,j-1}^{n}) - \bm{u}_{3,i,j-1}^{n}] , \label{eq:vlbm3} \\
\bm{u}_{4,i,j}^{n+1} &= \bm{M}_4(\bm{u}_{i,j+1}^{n}) + \theta_{i,j+1/2}^n[\bm{M}_4(\bm{u}_{i,j+1}^{n}) - \bm{u}_{4,i,j+1}^{n}] , \label{eq:vlbm4} \\
\bm{u}_{i,j}^{n+1} &= \bm{u}_{1,i,j}^{n+1} + \bm{u}_{2,i,j}^{n+1} + \bm{u}_{3,i,j}^{n+1} + \bm{u}_{4,i,j}^{n+1} + \alpha \bm{u}_{i,j}^{n} \nonumber \\
& ~~~~~ -\theta^{n}_{i-1/2,j} [\bm{M}_{2}(\bm{u}_{i,j}^{n}) - \bm{u}_{2,i,j}^{n}] - \theta^{n}_{i+1/2,j}[\bm{M}_{1}(\bm{u}_{i,j}^{n}) - \bm{u}_{1,i,j}^{n}] \nonumber \\
& ~~~~~ - \theta^{n}_{i,j-1/2} [\bm{M}_{4}(\bm{u}_{i,j}^{n}) - \bm{u}_{4,i,j}^{n}] - \theta^{n}_{i,j+1/2}[\bm{M}_{3}(\bm{u}_{i,j}^{n}) - \bm{u}_{3,i,j}^{n}] .\label{eq:vlbm5}
\end{align}
Note that we entirely skipped the decoupled evolution rule for \(\bm{u}_{5}\) for efficiency, as proposed in \cite{wissocq2024}.
Further algorithmic and numerical details, such as the definition of the blending parameters and the RLMP limiter used here, are given in \cite{wissocq2024}.

\subsection{Case setup of the Mach 2000 astrophysical jet}\label{sec:setup}

We model a highly supersonic, shock-dominated jet propagating into a stationary ambient medium.
The jet injection velocity is set to $v_{\mathrm{jet},1} = 800$, with an ambient pressure $p_{\mathrm{amb}} = 0.4127$ and a heat capacity ratio $\gamma = 5/3$.
The density ratio between the jet core and the ambient medium is $10$, with base values $\bar{\rho}_{\mathrm{jet}} = 5.0$ and $\rho_{\mathrm{amb}} = 0.5$.
Tracking the primary bow shock yields an empirical jet head velocity $v_{\mathrm{head}} \approx 666.7$.
The simulation is conducted on a two-dimensional Cartesian grid $\mathcal{D} = [0, x_{\mathrm{max}}] \times [y_{\mathrm{min}}, y_{\mathrm{max}}] = [0, 2.5] \times [-1, 1]$.

At the physical inlet $x=0$, we impose a Dirichlet boundary condition by setting the incoming distribution functions $\bm{u}_k$ to the local equilibrium $\bm{M}_k$ corresponding to a stochastically perturbed inlet state $\bm{u}_{\mathrm{inlet}}$, i.e.,
\begin{equation}\label{eq:inletBC}
    \bm{u}_k(0, y, t, \omega) = \bm{M}_k\Big(\bm{u}_{\mathrm{inlet}}(y, \omega)\Big),
\end{equation}
where $\bm{u}_{\mathrm{inlet}} = (\rho(y, \omega), \rho(y, \omega)v_1(y), 0, E(y, \omega))^\top$.
The velocity profile is defined as a step function $v_1(y) = v_{\mathrm{jet},1}$ for $|y| \le r_{\mathrm{jet}}$ (with $r_{\mathrm{jet}} = 0.05$) and $v_1(y)=0$ otherwise.
The density $\rho(y, \omega)$ is defined by the stochastic perturbation \eqref{eq:karhunenExp} (see Section~\ref{sec:statInit}) inside the core and the ambient density $\rho_{\mathrm{amb}}$ elsewhere, i.e.,
\begin{equation}\label{eq:rho_inlet_def}
    \rho(y, \omega) = \mathbb{I}_{\{|y| \le r_{\mathrm{jet}}\}} \rho'(y, \omega) + \mathbb{I}_{\{|y| > r_{\mathrm{jet}}\}} \rho_{\mathrm{amb}}.
\end{equation}
The total energy at the inlet $E(y, \omega)$ is derived from the ideal gas equation of state \eqref{eq:eos} using the prescribed ambient pressure $p_{\mathrm{amb}}$ and the inlet density \(\rho(y,\omega)\).
At $t=0$, the interior domain is initialized to the quiescent ambient state $\bm{u}_{\mathrm{amb}} = (\rho_{\mathrm{amb}}, 0, 0, p_{\mathrm{amb}}/(\gamma-1))^\top$, i.e.,
\begin{equation}\label{eq:u0_def}
    \bm{u}_0(x, y, \omega) = \mathbb{I}_{\{x=0\}} \bm{u}_{\mathrm{inlet}}(y, \omega) + \mathbb{I}_{\{x>0\}} \bm{u}_{\mathrm{amb}}.
\end{equation}

The outflow boundary at $x = x_{\mathrm{max}}$ carries a zero-gradient (homogeneous Neumann) condition on the full conserved state, $\partial_x \bm{u} = \bm{0}$, enforced by first-order extrapolation,
\begin{equation}\label{eq:outletBC}
    \bm{u}(x_{\mathrm{max}}, y, t,\omega) = \bm{u}(x_{\mathrm{max}} - \triangle x, y, t,\omega).
\end{equation}
The transverse boundaries at $y = \pm 1$ are treated identically, with a zero-gradient (homogeneous Neumann) condition on the full conserved state, $\partial_y \bm{u} = \bm{0}$, so that outgoing structures leave the transverse edges without reflection,
\begin{equation}\label{eq:transverseBC}
    \bm{u}(x, \pm 1, t,\omega) = \bm{u}(x, \pm 1 \mp \triangle x, t,\omega).
\end{equation}
Mesoscopically, both Neumann conditions are realized by copying the nearest interior nodal state onto the boundary layer of distribution functions.

Collecting the above, the deterministic problem solved for each fixed realization $\omega\in\Omega$ is the following random initial--boundary--value problem (IBVP) for the two-dimensional compressible Euler system on the space--time cylinder $\mathcal{D}\times(0,T]$:
\begin{subequations}\label{eq:ibvp}
\begin{alignat}{2}
    \partial_t \bm{u} + \partial_x \bm{f}(\bm{u}) + \partial_y \bm{g}(\bm{u}) &= \bm{0}, &\qquad& (\bm{x},t)\in\mathcal{D}\times(0,T], \label{eq:ibvp-pde}\\
    \bm{u}(x,y,0,\omega) &= \bm{u}_0(x,y,\omega), && (x,y)\in\mathcal{D}, \label{eq:ibvp-ic}\\
    \bm{u}(0,y,t,\omega) &= \bm{u}_{\mathrm{inlet}}(y,\omega), && y\in[-1,1],\ t\in(0,T], \label{eq:ibvp-inlet}\\
    \partial_x \bm{u}(x_{\mathrm{max}},y,t,\omega) &= \bm{0}, && y\in[-1,1],\ t\in(0,T], \label{eq:ibvp-outlet}\\
    \partial_y \bm{u}(x,\pm1,t,\omega) &= \bm{0}, && x\in[0,x_{\mathrm{max}}],\ t\in(0,T], \label{eq:ibvp-wall}
\end{alignat}
\end{subequations}
with fluxes $\bm{f},\bm{g}$ from \eqref{eq:fluxF}--\eqref{eq:fluxG}, closed by the equation of state \eqref{eq:eos}. The only source of randomness in \eqref{eq:ibvp} is the perturbed inlet state $\bm{u}_{\mathrm{inlet}}(\cdot,\omega)$ entering through the Dirichlet datum \eqref{eq:ibvp-inlet} and, via \eqref{eq:u0_def}, the initial datum \eqref{eq:ibvp-ic}; the domain geometry, the equation of state, and the outflow and transverse conditions \eqref{eq:ibvp-outlet}--\eqref{eq:ibvp-wall} are deterministic and identical across realizations. We denote by $\mathcal{S}_t\colon \bm{u}_0(\cdot,\omega)\mapsto \bm{u}(\cdot,t,\omega)$ the (formal) exact solution operator of the IBVP~\eqref{eq:ibvp}, advancing an admissible initial state to the entropy solution at time~$t$.

\subsection{Statistical solutions and stochastic initial data}\label{sec:statInit}
Let $(\Omega, \mathcal{F}, \mathbb{P})$ be a complete probability space, where $\Omega$ is the sample space of all realizations, $\mathcal{F}$ is the $\sigma$-algebra of measurable events, and $\mathbb{P}$ is the probability measure. The randomness of the IBVP~\eqref{eq:ibvp} enters solely through the perturbed inlet density, and hence through the inlet state $\bm{u}_{\mathrm{inlet}}(\cdot,\omega)$ and the initial datum $\bm{u}_0(\cdot,\omega)$ of \eqref{eq:u0_def}. We regard the initial datum as a measurable map $\bm{u}_0\colon \Omega \to L^1(\mathcal{D};\mathbb{R}^q)$, where $L^1(\mathcal{D};\mathbb{R}^q)$ is the Lebesgue space of integrable functions valued in the $q$-dimensional state space of conserved variables, with $q=4$ for the two-dimensional Euler system~\eqref{eq:euler}. This map pushes the abstract measure $\mathbb{P}$ forward to an initial probability measure on state space,
\begin{align}\label{eq:mu0}
    \mu_0 = (\bm{u}_0)_{\#}\mathbb{P} \in \mathcal{P}(L^1),
\end{align}
where the subscript $\cdot_{\#}$ denotes the pushforward (image measure) operator and $\mathcal{P}(L^1)$ is the space of Borel probability measures on $L^1$. Concretely, for any Borel set $B\subseteq L^1$ one has $\mu_0(B)=\mathbb{P}(\{\omega\in\Omega : \bm{u}_0(\cdot,\omega)\in B\})$, so that $\mu_0$ encodes the distribution of the random inlet uncertainty over admissible initial states.

The \emph{statistical solution} of the IBVP~\eqref{eq:ibvp} is then defined as the time-evolved image of this initial uncertainty under the exact solution operator $\mathcal{S}_t$ of \eqref{eq:ibvp}. That is, transporting each realization forward by $\mathcal{S}_t$ and collecting the resulting distribution, the statistical solution at time $t$ is the pushforward measure
\begin{align}\label{eq:mut_exact}
    \mu_t = (\mathcal{S}_t)_{\#}\mu_0 \in \mathcal{P}(L^1),
    \qquad t\in[0,T],
\end{align}
which assigns to each Borel set $B\subseteq L^1$ the probability $\mu_t(B)=\mathbb{P}(\{\omega : \mathcal{S}_t(\bm{u}_0(\cdot,\omega))\in B\})$ that the flow occupies a state in $B$ at time $t$. The measure $\mu_t$ is the central object whose convergence we study: weak--strong uniqueness asserts that if a strong solution $\bar{\bm{u}}$ of \eqref{eq:ibvp} exists then $\mu_t=\delta_{\bar{\bm{u}}(t)}$, whereas in the turbulent regime $\mu_t$ is expected to remain non-degenerate.

Because the exact operator $\mathcal{S}_t$ is unavailable, we approximate $\mu_t$ by replacing $\mathcal{S}_t$ with the discrete VLBM solution operator. Let $\mathcal{S}_t^{\triangle x}\colon L^1\to L^1$ denote the nonlinear forward operator induced at mesh width $\triangle x$ by the sequence of stream--collide updates \eqref{eq:vlbm1}--\eqref{eq:vlbm5}, together with the boundary treatment~\eqref{eq:ibvp-inlet}--\eqref{eq:ibvp-wall}, advancing the discrete state from $0$ to $t$. The time step is not an independent parameter: it is slaved to the mesh width through the CFL condition $\triangle t = \triangle x / a_n$, where $a_n$ denotes the kinetic speed of \eqref{eq:maxwellians12} and \eqref{eq:maxwellians34} re-evaluated at every time step $n$ so as to bound the current spectral radius of the flux Jacobians in the sense of Bouchut~\cite{bouchut1999}. Consequently $\triangle t = \mathcal{O}(\triangle x)$ and $\triangle t \to 0$ automatically as $\triangle x \to 0$, so that the single parameter $\triangle x$ indexes the entire space--time discretization. The \emph{approximate statistical solution} at mesh width $\triangle x$ is the pushforward of the same initial measure $\mu_0$ under this discrete operator,
\begin{align}\label{eq:pushforward}
    \mu_t^{\triangle x} = (\mathcal{S}_t^{\triangle x})_{\#}\mu_0 \in \mathcal{P}(L^1),
\end{align}
which is the discrete counterpart of $\mu_t$ in \eqref{eq:mut_exact}, differing from it only through the replacement of $\mathcal{S}_t$ by $\mathcal{S}_t^{\triangle x}$. Assuming consistency and stability of the scheme, $\mu_t^{\triangle x}\to\mu_t$ in the Wasserstein metric as $\triangle x\to0$; this is the discretization limit, controlled solely by $\triangle x$ through the CFL coupling above.

Computationally, neither $\mu_t$ nor $\mu_t^{\triangle x}$ is directly accessible, since both are pushforwards of the full measure $\mu_0$; instead we sample the inlet uncertainty and form an empirical measure. Drawing $M$ independent realizations $\omega_1,\dots,\omega_M\in\Omega$ (fixed by a pseudo-random seed) and propagating each by the discrete operator, $\bm{u}^{\triangle x}(\omega_m,t)=\mathcal{S}_t^{\triangle x}(\bm{u}_0(\omega_m))$, we obtain the empirical (Monte Carlo) measure
\begin{align}\label{eq:empirical}
    \mu_t^{\triangle x, M} = \frac{1}{M}\sum_{m=1}^{M} \delta_{\bm{u}^{\triangle x}(\omega_m,t)} \in \mathcal{P}(L^1),
\end{align}
where $M$ is the ensemble size, $m$ the sample index, and $\delta$ the Dirac measure. This empirical measure is the only object actually computed, and it is separated from the exact statistical solution $\mu_t$ by two independent errors. First, at fixed mesh width, the law of large numbers gives the \emph{sampling limit} $\mu_t^{\triangle x,M}\to\mu_t^{\triangle x}$ as $M\to\infty$, the empirical measure converging to the true discrete pushforward \eqref{eq:pushforward} in the Wasserstein metric at the Monte Carlo rate. Second, the \emph{discretization limit} $\mu_t^{\triangle x}\to\mu_t$ as $\triangle x\to0$ established above removes the scheme error. Convergence of the computed measure to the exact statistical solution therefore requires the \emph{joint limit}
\begin{align}\label{eq:jointlimit}
    \mu_t^{\triangle x, M} \longrightarrow \mu_t
    \qquad\text{as}\qquad
    \triangle x \to 0 \ \ \text{and}\ \ M \to \infty,
\end{align}
in which both the mesh is refined and the ensemble is enlarged. In practice we cannot take either limit to completion; we fix $M=1000$ and vary $\triangle x$ over the sequence of grids introduced in Section~\ref{sec:weak-strong}, and we verify a posteriori (Section~\ref{subsec:sampling}) that at these resolutions $M=1000$ lies above the sampling-error floor, so that the residual sampling error in \eqref{eq:jointlimit} is negligible relative to the discretization error and the measured rates reflect the discretization limit $\triangle x\to0$. All convergence metrics in Section~\ref{sec:weak-strong} are accordingly computed between empirical measures of the form \eqref{eq:empirical}.

\begin{figure}[ht!]
    \centering
    \begin{tikzpicture}
    \begin{axis}[
        width=8cm, 
        height=4cm,
        xlabel={$\overline{y}$},
        ylabel={$\rho'(\overline{y}, \omega)$},
        ylabel style={font=\footnotesize},
        yticklabel style={font=\footnotesize},
        xticklabel style={font=\footnotesize},
        xlabel style={font=\footnotesize},
        ymin=4.6, 
        ymax=5.4,
        grid=major,
        grid style={solid, black!60},
        legend pos=outer north east,
        legend style={font=\footnotesize, legend cell align=left},
        enlarge x limits=false
    ]
    
    \addplot[name path=upper, draw=none, forget plot] coordinates {
        (-1.00,5.000) (-0.96,5.001) (-0.92,5.005) (-0.88,5.011) (-0.84,5.019) (-0.80,5.030) (-0.76,5.042) (-0.72,5.056) (-0.68,5.072) (-0.64,5.089) (-0.60,5.106) (-0.56,5.124) (-0.52,5.143) (-0.47,5.161) (-0.43,5.179) (-0.39,5.196) (-0.35,5.212) (-0.31,5.228) (-0.27,5.241) (-0.23,5.254) (-0.19,5.264) (-0.15,5.273) (-0.11,5.281) (-0.07,5.287) (-0.03,5.290) (0.01,5.290) (0.05,5.288) (0.09,5.284) (0.13,5.277) (0.17,5.268) (0.21,5.257) (0.25,5.245) (0.29,5.231) (0.33,5.216) (0.37,5.199) (0.41,5.183) (0.45,5.166) (0.49,5.149) (0.54,5.132) (0.58,5.114) (0.62,5.097) (0.66,5.080) (0.70,5.065) (0.74,5.050) (0.78,5.037) (0.82,5.025) (0.86,5.015) (0.90,5.008) (0.94,5.003) (1.00,5.000)
    };
    \addplot[name path=lower, draw=none, forget plot] coordinates {
        (-1.00,5.000) (-0.96,4.999) (-0.92,4.995) (-0.88,4.990) (-0.84,4.982) (-0.80,4.972) (-0.76,4.960) (-0.72,4.947) (-0.68,4.932) (-0.64,4.915) (-0.60,4.898) (-0.56,4.879) (-0.52,4.860) (-0.47,4.839) (-0.43,4.818) (-0.39,4.797) (-0.35,4.778) (-0.31,4.760) (-0.27,4.744) (-0.23,4.730) (-0.19,4.719) (-0.15,4.709) (-0.11,4.701) (-0.07,4.695) (-0.03,4.690) (0.01,4.688) (0.05,4.689) (0.09,4.694) (0.13,4.702) (0.17,4.712) (0.21,4.723) (0.25,4.737) (0.29,4.753) (0.33,4.771) (0.37,4.790) (0.41,4.809) (0.45,4.830) (0.49,4.850) (0.54,4.870) (0.58,4.889) (0.62,4.907) (0.66,4.924) (0.70,4.940) (0.74,4.954) (0.78,4.966) (0.82,4.977) (0.86,4.986) (0.90,4.993) (0.94,4.997) (1.00,5.000)
    };
    \addplot[fill=gray, opacity=0.3] fill between[of=upper and lower];
    \addlegendentry{$\pm 1 \sigma$}
    
    \addplot[black, line width=1.5pt, densely dashed] coordinates {
        (-1.00,5.000) (-0.96,5.000) (-0.92,5.000) (-0.88,5.000) (-0.84,5.000) (-0.80,5.001) (-0.76,5.001) (-0.72,5.001) (-0.68,5.002) (-0.64,5.002) (-0.60,5.002) (-0.56,5.002) (-0.52,5.001) (-0.47,5.000) (-0.43,4.999) (-0.39,4.997) (-0.35,4.995) (-0.31,4.994) (-0.27,4.993) (-0.23,4.992) (-0.19,4.992) (-0.15,4.991) (-0.11,4.991) (-0.07,4.991) (-0.03,4.990) (0.01,4.989) (0.05,4.989) (0.09,4.989) (0.13,4.989) (0.17,4.990) (0.21,4.990) (0.25,4.991) (0.29,4.992) (0.33,4.993) (0.37,4.995) (0.41,4.996) (0.45,4.998) (0.49,4.999) (0.54,5.001) (0.58,5.002) (0.62,5.002) (0.66,5.002) (0.70,5.002) (0.74,5.002) (0.78,5.001) (0.82,5.001) (0.86,5.001) (0.90,5.000) (0.94,5.000) (1.00,5.000)
    };
    \addlegendentry{Mean ($M=1000$)}
    
    \addplot[blue, opacity=0.8, line width=1.5pt] coordinates {
        (-1.00,5.000) (-0.96,5.000) (-0.92,4.999) (-0.88,4.997) (-0.84,4.994) (-0.80,4.989) (-0.76,4.977) (-0.72,4.959) (-0.68,4.938) (-0.64,4.921) (-0.60,4.910) (-0.56,4.899) (-0.52,4.886) (-0.47,4.878) (-0.43,4.884) (-0.39,4.899) (-0.35,4.906) (-0.31,4.895) (-0.27,4.874) (-0.23,4.865) (-0.19,4.879) (-0.15,4.910) (-0.11,4.946) (-0.07,4.987) (-0.03,5.038) (0.01,5.094) (0.05,5.143) (0.09,5.169) (0.13,5.174) (0.17,5.170) (0.21,5.168) (0.25,5.170) (0.29,5.173) (0.33,5.169) (0.37,5.154) (0.41,5.128) (0.45,5.097) (0.49,5.073) (0.54,5.061) (0.58,5.057) (0.62,5.053) (0.66,5.045) (0.70,5.033) (0.74,5.021) (0.78,5.012) (0.82,5.006) (0.86,5.002) (0.90,5.001) (0.94,5.000) (1.00,5.000)
    };
    \addlegendentry{Sample $\omega_1$}
    
    \addplot[red, opacity=0.8, line width=1.5pt] coordinates {
        (-1.00,5.000) (-0.96,5.001) (-0.92,5.002) (-0.88,5.004) (-0.84,5.003) (-0.80,4.994) (-0.76,4.976) (-0.72,4.948) (-0.68,4.913) (-0.64,4.873) (-0.60,4.828) (-0.56,4.778) (-0.52,4.727) (-0.47,4.685) (-0.43,4.661) (-0.39,4.660) (-0.35,4.675) (-0.31,4.695) (-0.27,4.708) (-0.23,4.712) (-0.19,4.714) (-0.15,4.724) (-0.11,4.745) (-0.07,4.774) (-0.03,4.815) (0.01,4.873) (0.05,4.951) (0.09,5.041) (0.13,5.120) (0.17,5.175) (0.21,5.205) (0.25,5.224) (0.29,5.239) (0.33,5.251) (0.37,5.256) (0.41,5.252) (0.45,5.239) (0.49,5.217) (0.54,5.188) (0.58,5.155) (0.62,5.124) (0.66,5.098) (0.70,5.077) (0.74,5.059) (0.78,5.042) (0.82,5.027) (0.86,5.015) (0.90,5.006) (0.94,5.002) (1.00,5.000)
    };
    \addlegendentry{Sample $\omega_2$}
    
    \addplot[green!70!black, opacity=0.8, line width=1.5pt] coordinates {
        (-1.00,5.000) (-0.96,4.999) (-0.92,4.998) (-0.88,4.996) (-0.84,4.994) (-0.80,4.995) (-0.76,5.000) (-0.72,5.010) (-0.68,5.028) (-0.64,5.052) (-0.60,5.081) (-0.56,5.113) (-0.52,5.143) (-0.47,5.168) (-0.43,5.187) (-0.39,5.201) (-0.35,5.208) (-0.31,5.207) (-0.27,5.200) (-0.23,5.196) (-0.19,5.199) (-0.15,5.209) (-0.11,5.216) (-0.07,5.210) (-0.03,5.191) (0.01,5.164) (0.05,5.137) (0.09,5.108) (0.13,5.072) (0.17,5.028) (0.21,4.984) (0.25,4.950) (0.29,4.933) (0.33,4.930) (0.37,4.934) (0.41,4.940) (0.45,4.945) (0.49,4.948) (0.54,4.948) (0.58,4.943) (0.62,4.937) (0.66,4.936) (0.70,4.941) (0.74,4.951) (0.78,4.963) (0.82,4.973) (0.86,4.983) (0.90,4.991) (0.94,4.997) (1.00,5.000)
    };
    \addlegendentry{Sample $\omega_3$}
    \end{axis}
    \end{tikzpicture}
    \caption{Realizations of the truncated Karhunen--Lo\`eve density perturbation $\rho'(y, \omega)$ \eqref{eq:karhunenExp} at the inlet interface across the normalized transverse coordinate $\overline{y} = y/r_{\mathrm{jet}}$. Shown are the baseline mean, the $\pm 1\sigma$ variance envelope, computed from \(M=1000\) samples, respectively, and three independent sample trajectories.}
    \label{fig:density_perturbation}
\end{figure}
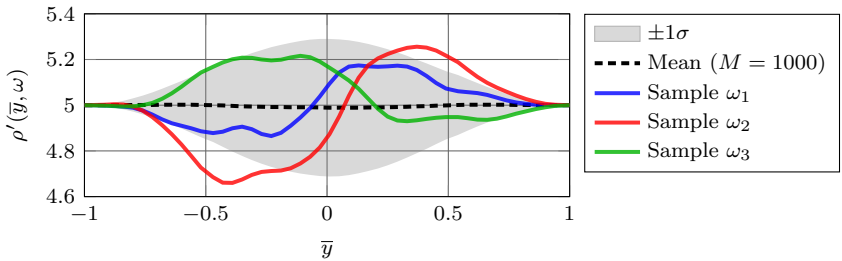

To appropriately evaluate the convergence of the statistical solution, the initial probability space must be sampled via correlated perturbations.
Randomized forcing of this kind mirrors established practice in computational astrophysics, where jet simulations are routinely seeded with stochastic perturbations: white noise injected continuously at the jet inlet to excite Kelvin--Helmholtz modes \cite{viallet2007}, random or sinusoidal velocity fluctuations of the launching disk \cite{kigure2005}, monochromatic and white-noise seeding of thin-shell instabilities in supersonic colliding flows \cite{mcleod2013}, and stochastically clumped ambient winds \cite{lopezmiralles2022}.
In the present setting the perturbation serves the additional purpose of sampling the initial probability measure.
We apply a smooth, spatially correlated random field to the macroscopic fluid density $\rho$ across the transverse spatial coordinate $y$ (orthogonal to the jet axis).
For a given realization $\omega \in \Omega$, the density perturbation in \eqref{eq:rho_inlet_def} takes the form of a truncated Karhunen--Lo\`eve-type spectral expansion
\begin{equation}\label{eq:karhunenExp}
    \rho'(y, \omega) = \bar{\rho}_{\mathrm{jet}} \left[ 1 + A \sum_{k=1}^{K} \frac{1}{k^p} \left( Y_k(\omega) \cos(k \pi y) + Z_k(\omega) \sin(k \pi y) \right) W(y) \right],
\end{equation}
where $\rho'$ is the stochastic density fluctuation and $A$ acts as a global amplitude scaling factor, applied directly to the random coefficients without normalization by the sample maximum; this preserves the statistical independence and natural variance of the modes, which is required for a well-defined Wasserstein convergence study.
To break symmetry while maintaining stability and positivity, we set the non-normalized amplitude to $A=0.1$.
Because \eqref{eq:karhunenExp} modulates only the inlet density while the inlet velocity and pressure remain deterministic, the internal jet Mach number $v_{\mathrm{jet},1}/\sqrt{\gamma p_{\mathrm{amb}}/\rho}$ scales as $\sqrt{\rho}$ and inherits half the relative density fluctuation: on the axis it varies by $\pm3\%$ about its nominal value at one standard deviation, and by at most $+7.5\%$ to $-8.1\%$ in the extreme case that all coefficients attain their bounds, while the Mach number referenced to the ambient sound speed is unaffected.
The sum in \eqref{eq:karhunenExp} iterates over the discrete wave number index $k$ up to the truncation mode limit $K=10$, with the exponent $p=2$ dictating the spatial decay rate of the higher frequencies.
The coefficients $Y_k(\omega)$ and $Z_k(\omega)$ are independent random variables drawn from a continuous uniform probability distribution $\mathcal{U}(-1,1)$ and $W(y)$ serves as a spatially localized windowing function to smoothly taper the perturbation at the jet boundaries; each realization of \eqref{eq:karhunenExp} is thus a smooth, symmetry-breaking modulation of the deterministic core density $\bar{\rho}_{\mathrm{jet}}$, with mode-wise variance decaying quadratically in the wave number.
Figure~\ref{fig:density_perturbation} illustrates the spatial structure of the stochastic extension of the inlet boundary condition \eqref{eq:karhunenExp}, displaying the deterministic mean density, the $\pm 1\sigma$ variance envelope, and three independent sample realizations of the truncated Karhunen--Lo\`eve perturbation across the normalized transverse coordinate.

This specific perturbation structure is mandatory for evaluating Wasserstein convergence.
The quadratic decay ($1/k^2$ for $p=2$) guarantees that the expected value of the Sobolev norm is bounded, $\mathbb{E}[\|\rho'\|_{H^s}^2] < \infty$, for an appropriate fractional Sobolev space $H^s$ of order $s$.
This ensures the initial probability measure is mean-square differentiable, preventing immediate unphysical Riemann singularities while successfully breaking the artificial symmetry of the Cartesian grid to trigger natural Kelvin--Helmholtz instabilities.
By evaluating the convergence of the approximated ensemble \eqref{eq:pushforward} via the exact 1-point discrete optimal transport formulation (see Section~\ref{sec:weak-strong}), we formally shift the evaluation paradigm from point-by-point flow tracking to macroscopic measure transport.

\section{Weak--strong uniqueness and convergence metrics}\label{sec:weak-strong}

Having constructed in Section~\ref{sec:statInit} the exact statistical solution $\mu_t=(\mathcal{S}_t)_{\#}\mu_0$, its discrete pushforward $\mu_t^{\triangle x}=(\mathcal{S}_t^{\triangle x})_{\#}\mu_0$, and the empirical measure $\mu_t^{\triangle x,M}$ of \eqref{eq:empirical} actually computed at mesh width $\triangle x$ and ensemble size $M$, we now describe how convergence toward $\mu_t$ under the joint limit \eqref{eq:jointlimit} is assessed. The guiding theoretical notion is the weak--strong uniqueness principle: whenever a classical strong solution $\bar{\bm{u}}$ of the Euler system~\eqref{eq:euler} exists, the exact statistical solution must reduce to the Dirac measure concentrated on it, $\mu_t=\delta_{\bar{\bm{u}}(t)}$, in which $\bar{\bm{u}}(t)$ is the exact deterministic state at time~$t$ and $\delta_{\bar{\bm{u}}(t)}$ is the Dirac distribution centered at $\bar{\bm{u}}(t)$, a measure of zero variance. In the turbulent, hyper-compressible regime studied here no global strong solution is expected to exist; individual realizations $\bm{u}^{\triangle x}(\omega_m,\cdot,t)=\mathcal{S}_t^{\triangle x}(\bm{u}_0(\omega_m))$ do not converge under refinement, and the limiting measure is generically non-degenerate. We stress that the principle is conditional and yields no analytical convergence guarantee here: the discontinuous inlet data rule out the regular strong solutions required by its hypothesis, and convex-integration constructions \cite{delellis} show that entropy admissibility alone does not restore uniqueness. Any convergence observed below is therefore a statement about the sequence of approximations selected by the VLBM discretization, a numerical counterpart of the weak--strong dichotomy, and not a theorem-backed statement about the underlying continuous problem. Convergence must therefore be assessed not through a single pathwise error but through a family of statistical observables, each probing a distinct functional of the ensemble $\{\bm{u}^{\triangle x}(\omega_m,\cdot,t)\}_{m=1}^{M}$ underlying $\mu_t^{\triangle x,M}$. Throughout we hold the ensemble size fixed at $M=1000$ and refine $\triangle x$, having established in Section~\ref{sec:statInit} that at these resolutions the sampling error is negligible relative to the discretization error, so that the measured rates probe the discretization limit $\triangle x\to0$ of the joint limit~\eqref{eq:jointlimit}; the role of $M$ is quantified separately in Section~\ref{subsec:sampling}. We describe the observables below, together with the two comparison operators through which their inter-resolution errors are formed, and we make explicit that the measured convergence rate of any observable depends on the operator chosen to relate successive resolutions.

Throughout this section the probability space $(\Omega,\mathcal{F},\mathbb{P})$, the ensemble size~$M$, the realizations~$\omega_m$, and the grid-indexed samples $\bm{u}^{\triangle x}(\omega_m,\bm{x},t)$ are exactly those defined in Section~\ref{sec:statInit}; we reuse them without redefinition. Unless stated otherwise, all results below use the full ensemble of $M=1000$ realizations at every resolution.

\subsection{Sampling}\label{subsec:sampling}
Two schedules relating the sample size to the mesh are conceivable. In the \emph{fixed-sample} schedule the same number~$M$ of realizations is used at every resolution, whereas in the \emph{resolution-tied} schedule, adopted in parts of the literature to accelerate the apparent convergence of statistical quantities, the sample size at a resolution of $N$ cells per transverse direction is set to $M(N)=N$ and thus grows with the mesh. We have examined both and report all results below for the fixed-sample schedule with $M=1000$. The resolution-tied schedule was found to have no appreciable effect on the measured convergence rates, for two reasons. First, since the sample sets are nested and share a common seed, the reduced coarse-grid ensembles used under the resolution-tied schedule are prefixes of the full thousand-sample ensemble rather than independent draws, so the two schedules are strongly correlated. Second, and more importantly, at the resolutions considered here the estimator is already sample-converged: $M=1000$ lies above the sampling-error floor, so reducing the coarse-grid sample count neither exposes nor obscures the underlying discretization rate. In the language of the joint limit~\eqref{eq:jointlimit}, fixing $M=1000$ freezes the sampling error $\mu_t^{\triangle x,M}-\mu_t^{\triangle x}$ at a level below the discretization error $\mu_t^{\triangle x}-\mu_t$ across the entire range of $\triangle x$ examined, so that the $\triangle x$-sweep isolates the discretization limit while the sampling limit $M\to\infty$ remains effectively resolved. Across all observables the two schedules agreed to within $0.07$ in fitted rate, a discrepancy comparable to the scatter of the fits themselves. Thus, in the present setting the fixed-sample estimate is assumed to be already converged, and we do not distinguish the two schedules further.

\subsection{Restriction to a common grid}\label{subsec:restriction}
Comparing solutions computed on different lattices requires a common representation. We introduce a conservative restriction operator~$\mathcal{R}$ that maps a macroscopic lattice field onto a fixed coarse grid by local block-averaging, assigning to each coarse cell the volume-weighted average of the fine-lattice nodal values it contains. Because the VLBM stores a single macroscopic state $\bm{u}^{\triangle x}=\sum_{k=1}^{5}\bm{u}_k^{\triangle x}$ at each lattice node, this averaging is exact and conserves the integral of every component, introducing no interpolation error. In all comparisons reported below, the macroscopic fields at every resolution are restricted by~$\mathcal{R}$ onto the coarsest grid, of $N_{\bm{x}}=500\times400$ cells, employed in the study, and the resulting coarse-grained fields are compared cell by cell. We denote the cell centers of this common grid by $\{\bm{x}_k\}_{k=1}^{N_{\bm{x}}}$. Restricting to the coarsest common grid confines the comparison to the large-scale, resolved content of the flow: structure below the scale of the common grid is averaged out prior to comparison, which is the appropriate choice when assessing convergence of the statistical, coarse-grained state rather than of the under-resolved fine-scale features that, by the pathwise argument above, are not expected to converge.

\subsection{Comparison operators}\label{subsec:comparison}
Let the available grids be ordered by increasing resolution and indexed by their streamwise cell count $N_x\in\{500,1000,2000,4000\}$ on the domain $[0,2.5]\times[-1,1]$, so that the four grids are $500\times400$, $1000\times800$, $2000\times1600$, and $4000\times3200$, the last carrying $12.8$ million cells. Writing $N_x^{(\ell)}$ for the $\ell$-th grid, $\ell=0,\dots,L$ with $L=3$, the corresponding mesh width is $\triangle x_\ell = 2.5/N_x^{(\ell)}$, giving $\triangle x_\ell\in\{5.0\times10^{-3},\,2.5\times10^{-3},\,1.25\times10^{-3},\,6.25\times10^{-4}\}$. All error functionals are evaluated after restriction to the coarsest grid $500\times400$, on which the fields are compared cell by cell; in the tables below each compared pair is identified by the streamwise cell counts $N_x$ of its two members. For an observable~$\mathfrak{O}$, which takes as its two arguments a coarse and a restricted fine field and returns a non-negative error, we form its inter-resolution error in one of two ways. In the \emph{Cauchy} comparison, successive resolutions are compared against one another,
\begin{equation}\label{eq:cauchy}
    \mathcal{E}^{\mathrm{C}}_{\mathfrak{O}}(\ell)
      = \mathfrak{O}\!\left(\mu_t^{\triangle x_{\ell},M},\, \mu_t^{\triangle x_{\ell-1},M}\right),
      \qquad \ell=1,\dots,L,
\end{equation}
and the error at level~$\ell$ is associated with the coarser mesh width~$\triangle x_{\ell-1}$ of the pair. This construction requires no external reference and measures whether the sequence of approximations is internally convergent in the sense of Cauchy. In the \emph{reference} comparison, each resolution is instead compared against the finest available empirical solution~$\mu_t^{\triangle x_L,M}$, which serves as a surrogate for the unknown exact statistical solution,
\begin{equation}\label{eq:ref}
    \mathcal{E}^{\mathrm{R}}_{\mathfrak{O}}(\ell)
      = \mathfrak{O}\!\left(\mu_t^{\triangle x_{\ell},M},\, \mu_t^{\triangle x_{L},M}\right),
      \qquad \ell=0,\dots,L-1,
\end{equation}
the error at level~$\ell$ being associated with the mesh width~$\triangle x_{\ell}$ of the coarser member. For each observable and each comparison operator we estimate a convergence rate~$\mathfrak{r}_{\mathfrak{O}}$ by least-squares regression in doubly logarithmic coordinates,
\begin{equation}\label{eq:ratefit}
    (\mathfrak{r}_{\mathfrak{O}},\, c)
      \;=\; \operatorname*{arg\,min}_{(\mathfrak{r},\,c)\,\in\,\mathbb{R}^{2}}
        \sum_{\ell}
        \bigl(\log \mathcal{E}_{\mathfrak{O}}(\ell) - \mathfrak{r}\,\log \triangle x_{(\ell)} - c\bigr)^{2},
\end{equation}
where $\triangle x_{(\ell)}$ denotes the mesh width associated with the error at level~$\ell$ under the chosen comparison operator, so that $\mathcal{E}_{\mathfrak{O}}\sim(\triangle x)^{\mathfrak{r}_{\mathfrak{O}}}$ to leading order; the rate $\mathfrak{r}_{\mathfrak{O}}$ is thus the slope of the error in doubly logarithmic coordinates. 

\subsection{Observables}\label{subsec:observables}
We evaluate five observables. In each case the arguments of~$\mathfrak{O}$ are the coarse field and the restricted fine field selected by the chosen comparison operator, both on the common grid; to lighten the notation we write $\bm{u}^{\triangle x}$ and $\bm{u}^{\triangle x/2}$ for the coarser and finer member of a compared pair, with the understanding that under the \textit{reference} comparison the finer member is always the finest solution~$\bm{u}^{\triangle x_L}$.

The first observable is purely pathwise and serves as a control. For each realization~$\omega_m$ we form the relative $L^1(\mathcal{D})$ difference between the coarse and restricted fine density fields and average over the ensemble,
\begin{equation}\label{eq:strongError}
    \mathcal{E}_{\mathrm{strong}}(t)
      = \frac{1}{M}\sum_{m=1}^{M}
        \frac{\bigl\| \rho^{\triangle x}(\omega_m,\cdot,t)
              - \mathcal{R}\,\rho^{\triangle x/2}(\omega_m,\cdot,t)\bigr\|_{L^{1}(\mathcal{D})}}
             {\bigl\| \mathcal{R}\,\rho^{\triangle x/2}(\omega_m,\cdot,t)\bigr\|_{L^{1}(\mathcal{D})}},
\end{equation}
where $\|\cdot\|_{L^1(\mathcal{D})}$ is the spatial $L^1$ norm over the common grid and the same realization index~$m$ appears in both fields, so that the difference isolates the effect of refinement on an individual sample. Because the two fields share the pseudo-random seed~$\omega_m$, a non-vanishing limit of $\mathcal{E}_{\mathrm{strong}}$ as $\triangle x\to0$ is the numerical signature of pathwise non-uniqueness: individual realizations do not converge, even though, as shown below, their statistics do.

The remaining four observables are statistical, and are evaluated on the density component $\rho$ of the state (the perturbed variable; see the scoping remark at the end of this subsection). Writing the ensemble mean of the coarse density field at cell~$\bm{x}_k$ as the empirical expectation $\mathbb{E}_M[\rho^{\triangle x}](\bm{x}_k)=\tfrac{1}{M}\sum_{m=1}^{M}\rho^{\triangle x}(\omega_m,\bm{x}_k,t)$, consistent with the expectation operator $\mathbb{E}$ used in Section~\ref{sec:statInit} but taken over the finite ensemble, the mean error is the cell-averaged absolute difference of the two mean fields,
\begin{equation}\label{eq:mean}
    \mathcal{E}_{\mathrm{mean}}(t)
      = \frac{1}{N_{\bm{x}}}\sum_{k=1}^{N_{\bm{x}}}
        \bigl| \mathbb{E}_M[\rho^{\triangle x}](\bm{x}_k)
             - \mathcal{R}\,\mathbb{E}_M[\rho^{\triangle x/2}](\bm{x}_k)\bigr|.
\end{equation}
Likewise, denoting by $\sigma^{\triangle x}(\bm{x}_k)$ the cell-wise ensemble standard deviation of the coarse density field, that is the square root of the population variance $\mathbb{E}_M\bigl[(\rho^{\triangle x}-\mathbb{E}_M[\rho^{\triangle x}])^{2}\bigr]$ evaluated at $\bm{x}_k$, the standard-deviation error is
\begin{equation}\label{eq:std}
    \mathcal{E}_{\mathrm{std}}(t)
      = \frac{1}{N_{\bm{x}}}\sum_{k=1}^{N_{\bm{x}}}
        \bigl| \sigma^{\triangle x}(\bm{x}_k)
             - \mathcal{R}\,\sigma^{\triangle x/2}(\bm{x}_k)\bigr|.
\end{equation}
The mean~\eqref{eq:mean} and standard deviation~\eqref{eq:std} quantify convergence of the first- and second-moment one-point statistics of the density.

To probe the full one-point law rather than only its first two moments, we employ the Wasserstein distance between single-point marginals. At each cell~$\bm{x}_k$ let $\nu^{\triangle x}_{\bm{x}_k,t}$ denote the empirical law of the density over the ensemble, that is the one-dimensional distribution formed by the $M$ sampled values $\{\rho^{\triangle x}(\omega_m,\bm{x}_k,t)\}_{m=1}^{M}$ at that cell. The one-point Wasserstein error is the spatial average of the first Wasserstein distance between the coarse and fine marginals,
\begin{equation}\label{eq:w1cont}
    \mathcal{W}_{1,1}(t)
      = \frac{1}{N_{\bm{x}}}\sum_{k=1}^{N_{\bm{x}}}
        \inf_{\pi\in\Pi(\nu^{\triangle x}_{\bm{x}_k,t},\,\nu^{\triangle x/2}_{\bm{x}_k,t})}
        \int_{\mathbb{R}\times\mathbb{R}} |\xi-\zeta|\,\mathrm{d}\pi(\xi,\zeta),
\end{equation}
in which $\xi$ and $\zeta$ are admissible scalar values of the density, $|\xi-\zeta|$ is the one-dimensional cost of transporting a unit of probability mass from~$\xi$ to~$\zeta$, and $\Pi(\nu_1,\nu_2)$ is the set of couplings, that is joint probability measures on $\mathbb{R}\times\mathbb{R}$ with marginals $\nu_1$ and~$\nu_2$; the infimum over couplings is the optimal-transport (earth-mover's) cost. For empirical marginals of equal size~$M$ this infimum admits the closed form
\begin{equation}\label{eq:w1disc}
    \mathcal{W}_{1,1}(t)
      \approx \frac{1}{N_{\bm{x}}}\sum_{k=1}^{N_{\bm{x}}}
        \frac{1}{M}\sum_{m=1}^{M}
        \bigl| \rho^{\triangle x}_{(m)}(\bm{x}_k)
             - \mathcal{R}\,\rho^{\triangle x/2}_{(m)}(\bm{x}_k)\bigr|,
\end{equation}
where $\rho_{(1)}\le\rho_{(2)}\le\dots\le\rho_{(M)}$ denote the ascending order statistics of the ensemble of density values at the cell in question, so that the distance reduces to the mean absolute difference between the sorted samples of the two marginals. 
Note that the order-statistic index~$m$ here plays the role of the discrete inverse cumulative distribution function. 
Equivalently, by the standard one-dimensional optimal-transport identity, for any two probability measures $\nu_1,\nu_2$ on $\mathbb{R}$ with cumulative distribution functions $F_{\nu_1},F_{\nu_2}$ and quantile functions $F^{-1}_{\nu_1},F^{-1}_{\nu_2}$,
\begin{equation}\label{eq:w1cdf}
    W_1(\nu_1,\nu_2)
      = \int_{\mathbb{R}} \bigl|F_{\nu_1}(\xi)-F_{\nu_2}(\xi)\bigr|\,\mathrm{d}\xi
      = \int_{0}^{1} \bigl|F^{-1}_{\nu_1}(\alpha)-F^{-1}_{\nu_2}(\alpha)\bigr|\,\mathrm{d}\alpha,
\end{equation}
of which \eqref{eq:w1disc} is the empirical instance of the right-hand side. The identity \eqref{eq:w1cdf} connects the metric assessed here to the quantile diagnostics of Section~\ref{subsec:centerline} and to the cumulative-distribution reading of the one-point law developed in Section~\ref{subsec:pdfheatmap}.

Finally, to detect discrepancies in the joint spatial structure that the one-point marginals cannot resolve, we evaluate a two-point Wasserstein distance. Following the construction of Lye~\cite{lye2020}, we subsample the domain to $n_{\mathrm{pd}}$ points per spatial direction, and for each ordered pair of retained points~$(\bm{x}_k,\bm{x}_{k'})$ we form the empirical joint law $\rho^{\triangle x}_{k,k'}$ of the density pair $\bigl(\rho^{\triangle x}(\cdot,\bm{x}_k,t),\rho^{\triangle x}(\cdot,\bm{x}_{k'},t)\bigr)$ over the ensemble, a distribution on~$\mathbb{R}^{2}$. The two-point error is the average, over all $P=\bigl(n_{\mathrm{pd}}^{\,d}\bigr)^{2}$ such point-pairs in $d$ spatial dimensions, of the Wasserstein distance between the coarse and fine joint laws,
\begin{equation}\label{eq:w2}
    \mathcal{W}_{1,2}(t)
      = \frac{1}{P}\sum_{(k,k')}
        \mathcal{W}\!\Bigl(\rho^{\triangle x}_{k,k'},\,
          \mathcal{R}\,\rho^{\triangle x/2}_{k,k'}\Bigr),
\end{equation}
where $\mathcal{W}(\cdot,\cdot)$ denotes the first Wasserstein distance on~$\mathbb{R}^{2}$ with Euclidean ground cost, computed exactly by linear programming. We use $n_{\mathrm{pd}}=10$ points per direction throughout.

All five observables are evaluated on the density field~$\rho$, the component through which the random inlet perturbation~\eqref{eq:karhunenExp} enters the problem and the primary diagnostic for the Kelvin--Helmholtz roll-up and shock structure of the jet. Restricting the statistical analysis to the density marginal of the statistical solution is a deliberate reduction: the marginals $\mu_t^{\triangle x,M}\!\restriction_\rho$ are themselves well-defined probability measures whose convergence is a genuine and self-contained statement about $\mu_t$, and the full-state optimal-transport problem, posed on up to $12.8$ million cells with a four-dimensional ground space, is computationally prohibitive at the ensemble sizes required for a converged Wasserstein estimate. The convergence of the density marginal established here does not by itself imply convergence of the momentum or energy marginals, which are coupled to~$\rho$ through the flux \eqref{eq:fluxF}--\eqref{eq:fluxG} but need not share its rate. 
Extending the same analysis to the remaining conserved components is straightforward within the present pipeline and is left to future work.

\section{Numerical results}\label{sec:results}

\subsection{Implementation}
The methodology from \cite{wissocq2024} (see Section~\ref{sec:methods}) has been ported to GPU for the present work.
The Monte Carlo ensemble was implemented using Python and CuPy for CUDA-accelerated array processing.
The solver and the complete post-processing are available as the open-source library acropora~\cite{acropora}, which regenerates the ensembles of this paper seed for seed on GPUs. 
The ensembles restricted to the coarsest grid, together with the convergence and centerline results behind all tables and figures, are archived in~\cite{acroporadata}. 
To isolate the spatial convergence of the probability measure and separate grid discretization errors from statistical sampling noise, the Monte Carlo ensemble is computed across the sequence of progressively refined grids introduced in Section~\ref{subsec:comparison}.
The spatial resolution scales incrementally, doubling from the coarse baseline of $500 \times 400$ cells to the highly resolved grid of $4000 \times 3200$ discrete lattice cells ($12.8$ million cells).
The time-stepping is implemented as described in \cite{wissocq2024}.
The entire pipeline was executed in an array of 20 batch jobs with 50 samples each per GPU node on HoreKa Ruby (4x NVIDIA H200 GPUs each).
A single sample of the highest resolution ($12.8$ million grid points) took less than six minutes on a single GPU including I/O.
One batch took below 6.6 hours on average to compute all grids for 50 sample trajectories on one GPU node; the subsequent parallel CPU post-processing of the data into error metrics and visualizations required roughly 21 hours.

\subsection{Visualization of the density evolution}\label{subsec:visualization}

\begin{figure}[ht!]
    \centering
    \includegraphics[width=0.99\textwidth]{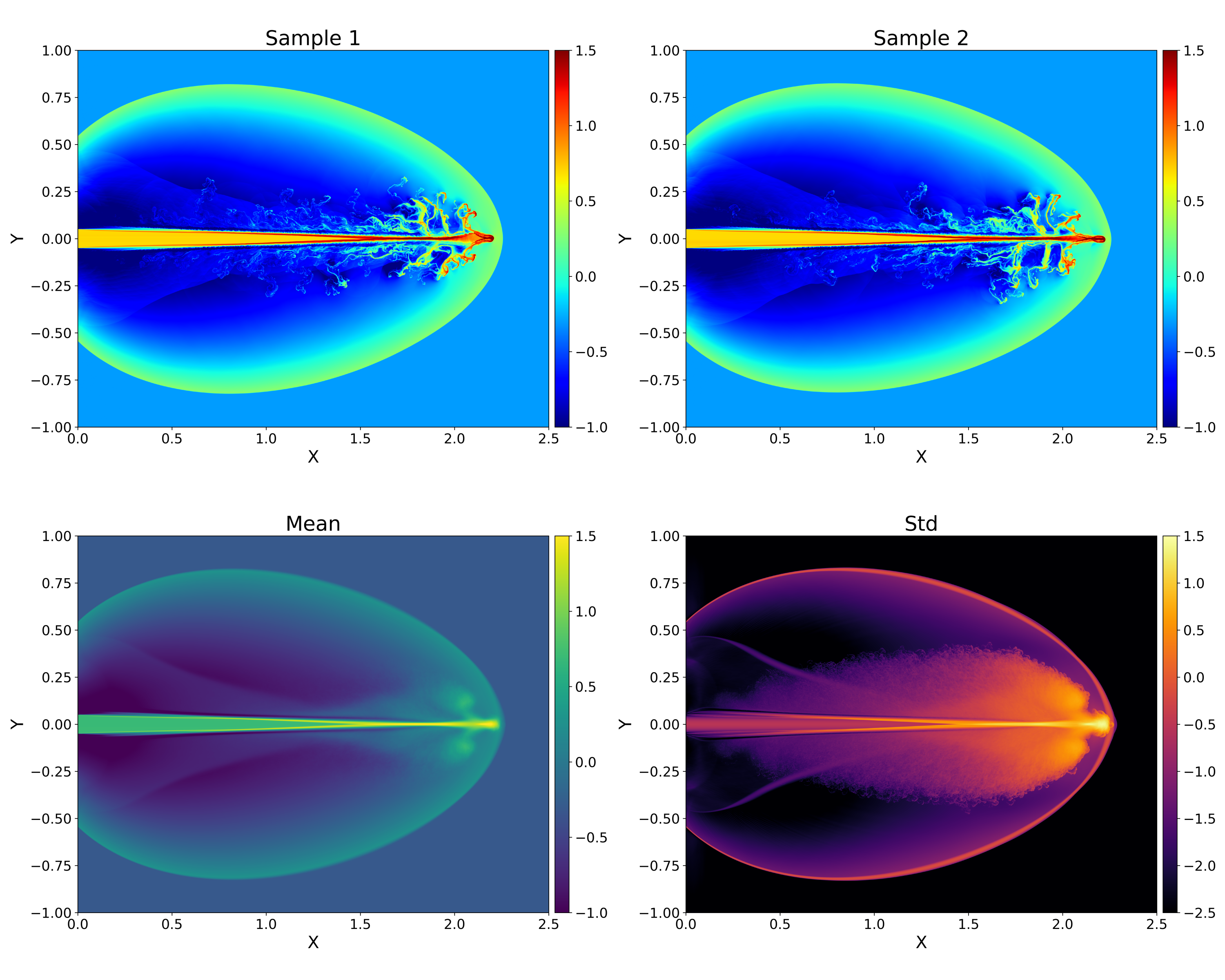}\\
    \caption{Two different samples (top row) and statistical moments of \(M=1000\) samples (bottom row: mean and standard deviation) of the density at $t=0.0035\,\text{s}$ on the $4000\times 3200$ grid.
    Colorbars are in logarithmic scale.
    }
    \label{fig:overview_final}
\end{figure}

\begin{figure}[ht!]
    \centering
    \includegraphics[width=0.32\textwidth]{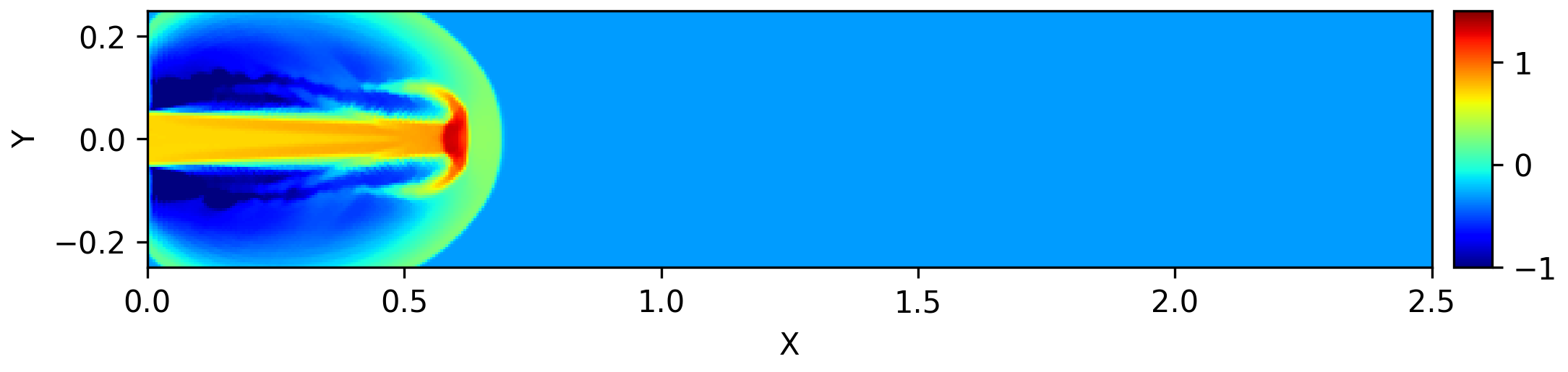}
    \includegraphics[width=0.32\textwidth]{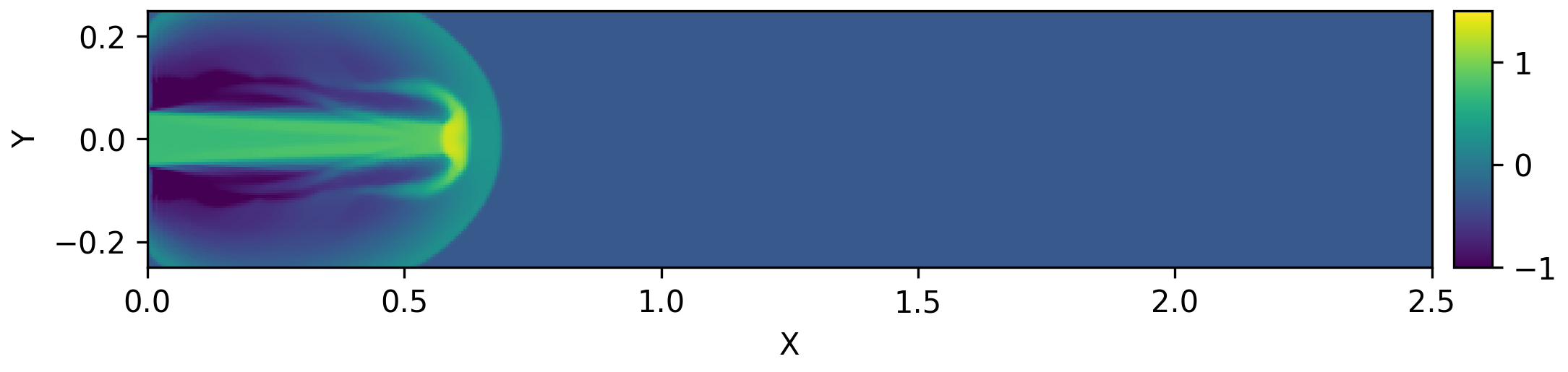}
    \includegraphics[width=0.32\textwidth]{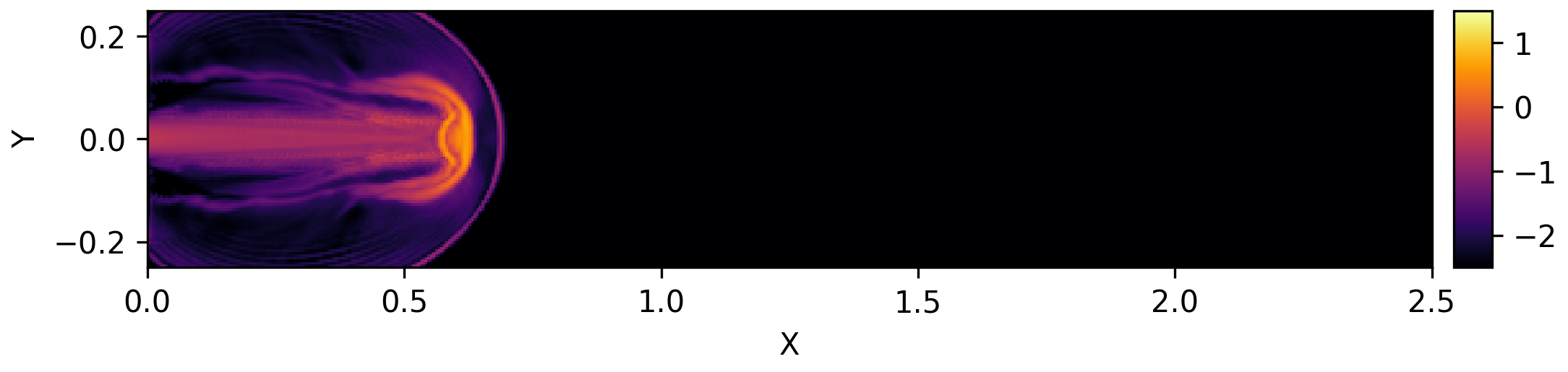} \\
    \includegraphics[width=0.32\textwidth]{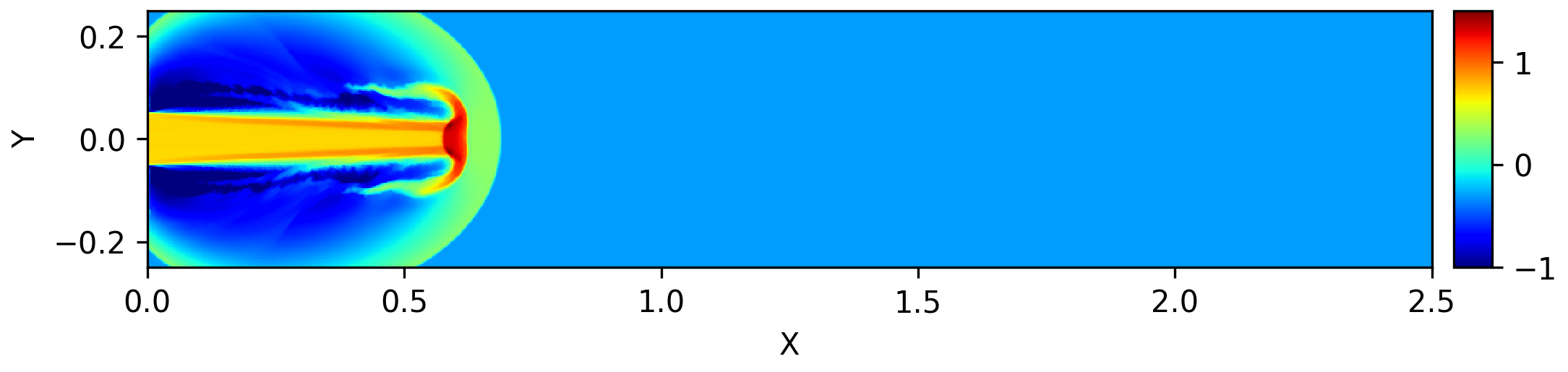}
    \includegraphics[width=0.32\textwidth]{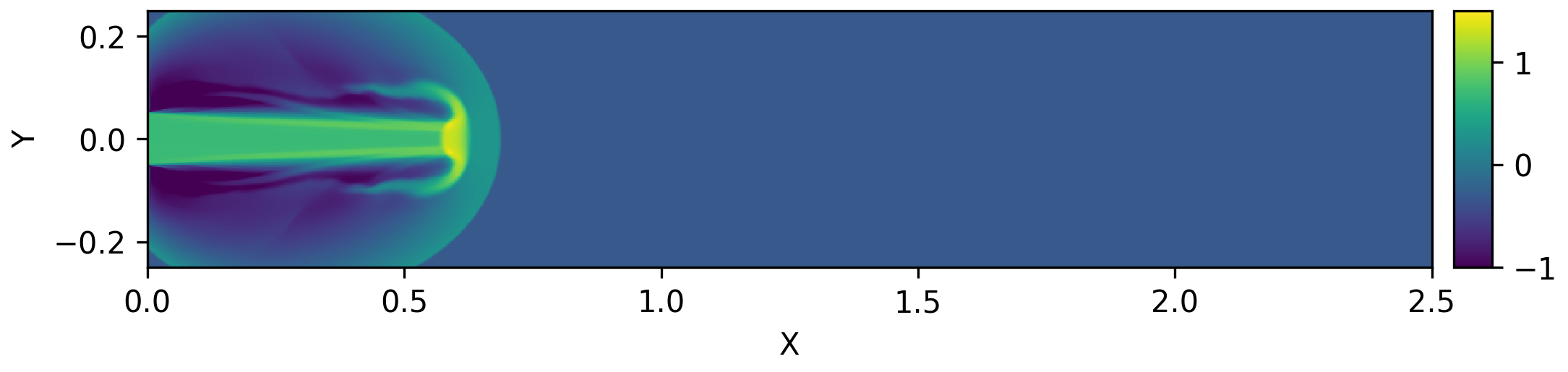}
    \includegraphics[width=0.32\textwidth]{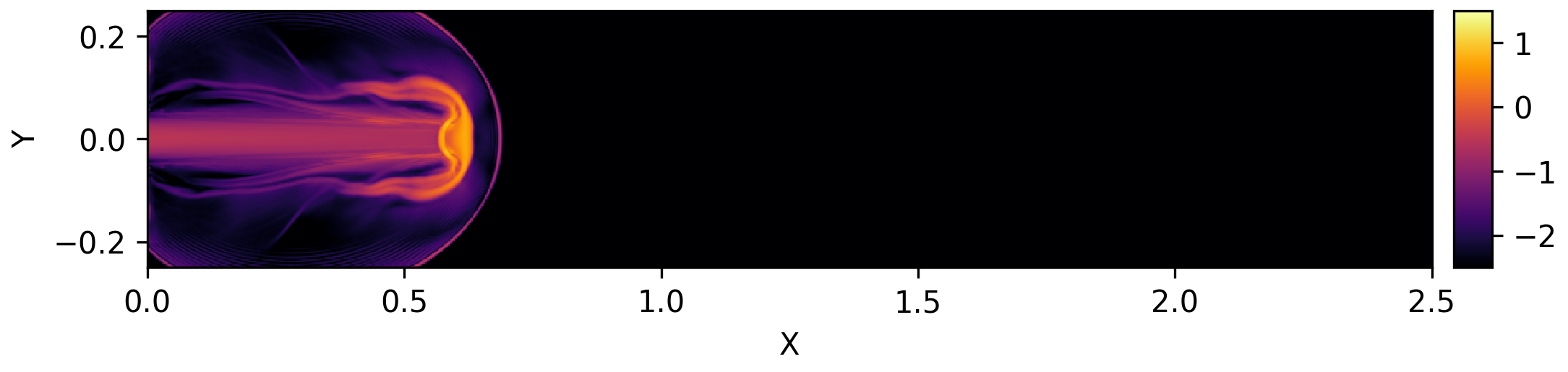}\\
    \includegraphics[width=0.32\textwidth]{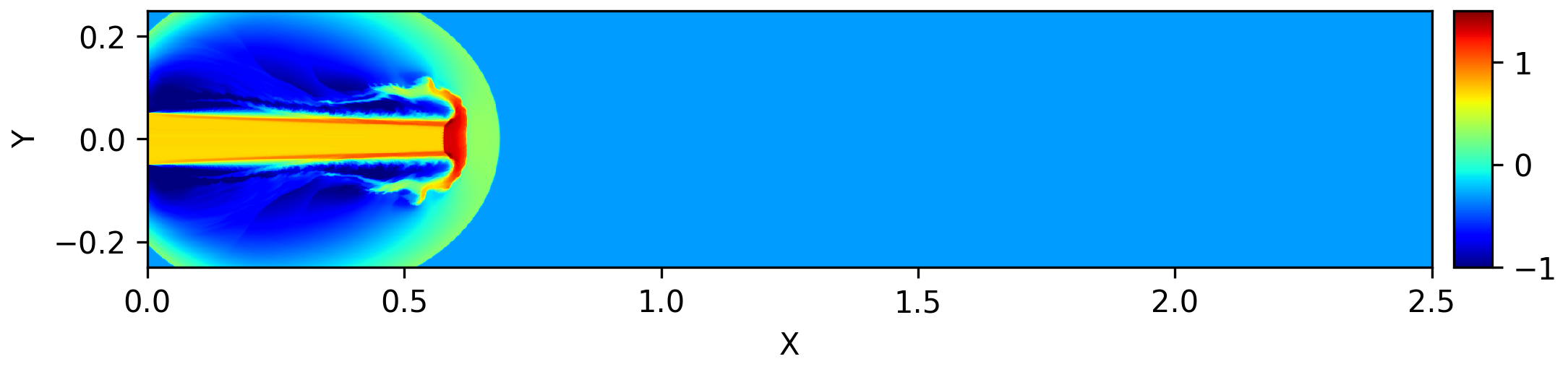}
    \includegraphics[width=0.32\textwidth]{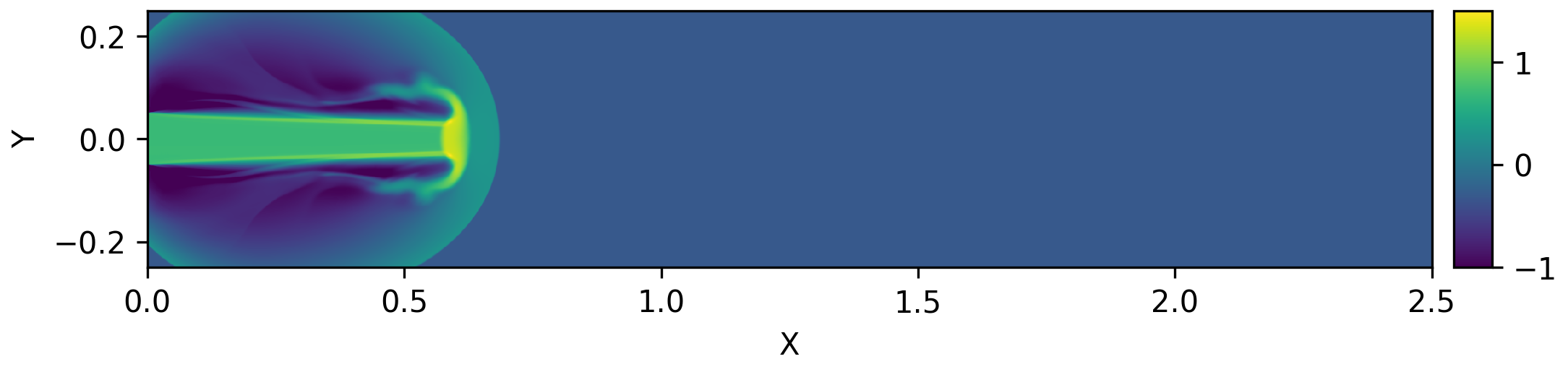}
    \includegraphics[width=0.32\textwidth]{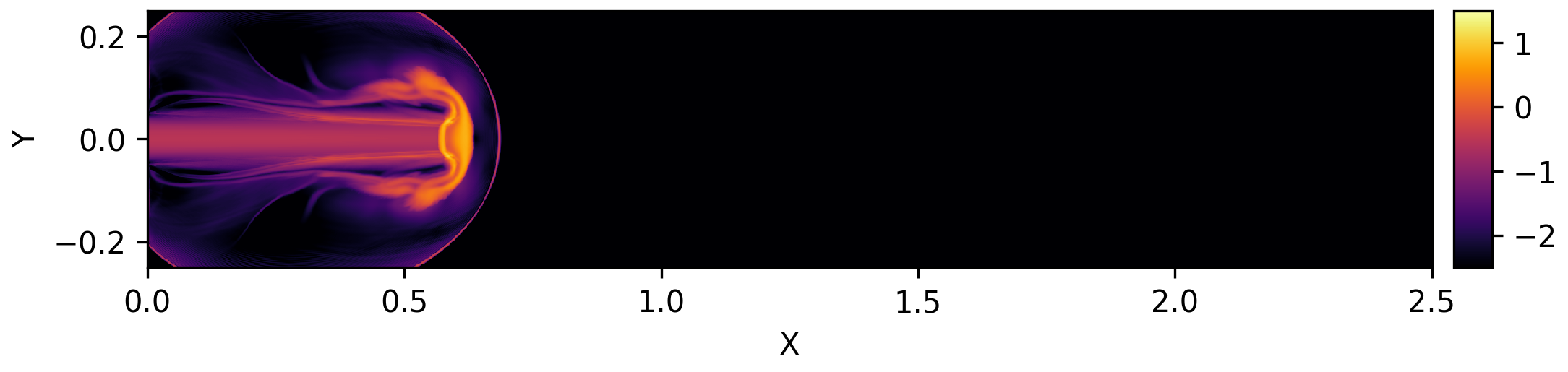}\\
    \includegraphics[width=0.32\textwidth]{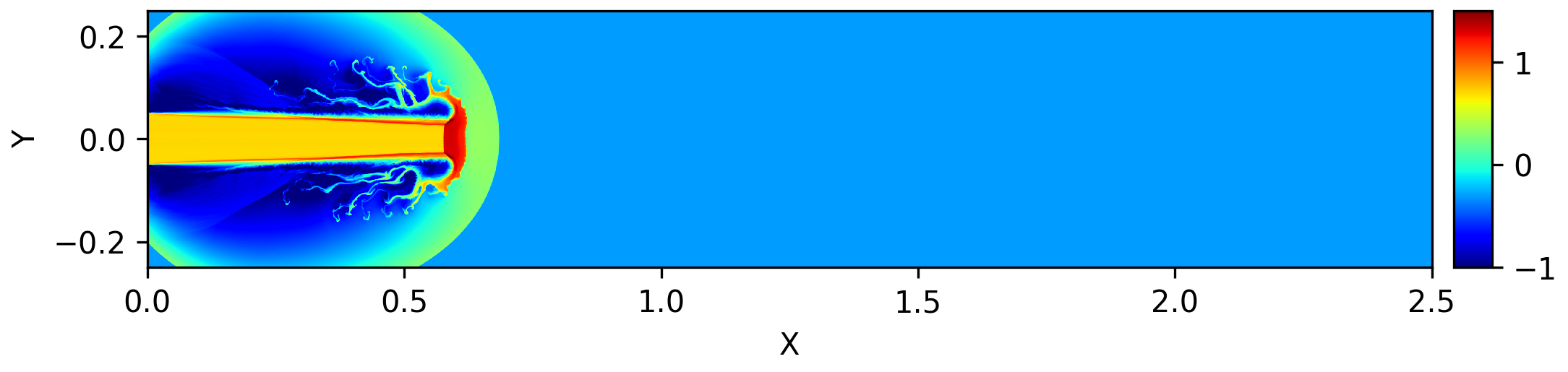}
    \includegraphics[width=0.32\textwidth]{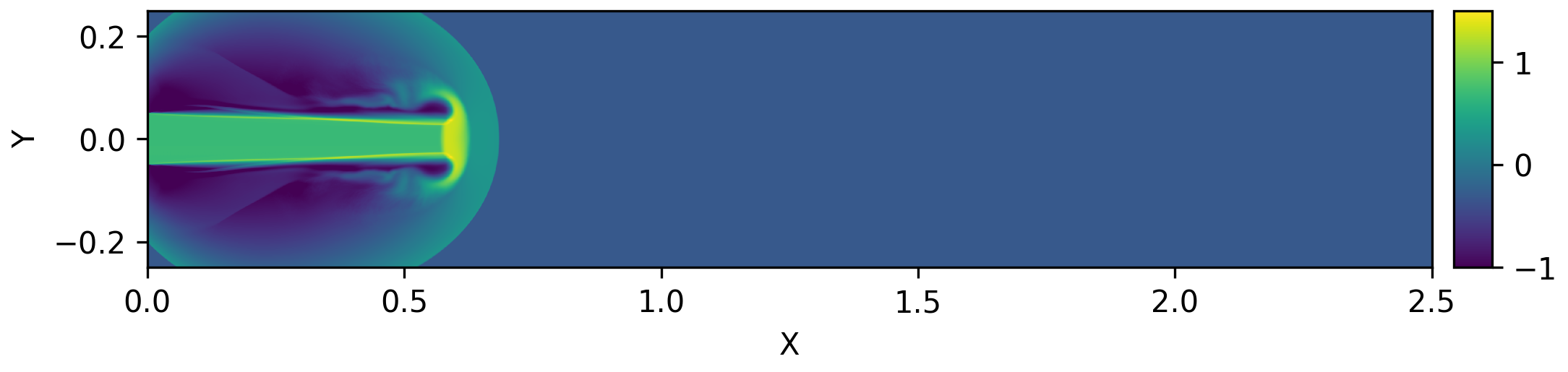}
    \includegraphics[width=0.32\textwidth]{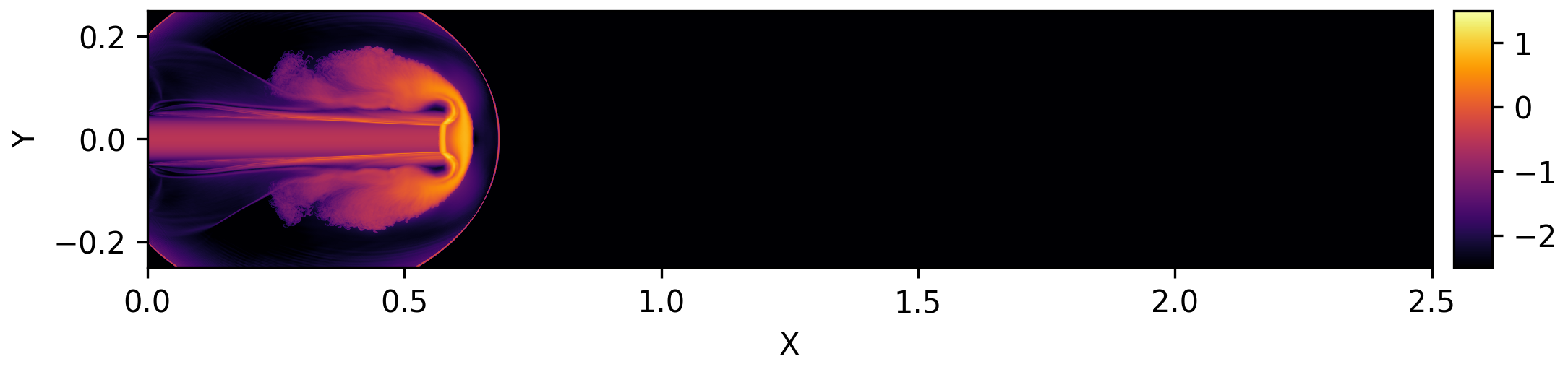}\\
    \caption{
        Single sample and statistical moments (mean and standard deviation for \(M=1000\)) of the density at $t=0.001\,\text{s}$ (left to right).
        The layout displays \(x\)-resolutions \(N=500, 1000, 2000, 4000\) (top to bottom) on the clipped domain \(\overline{\mathcal{D}} = [0,2.5] \times [-0.25, 0.25]\).
        Colorbars are in logarithmic scale.
    }
    \label{fig:jet_convergence_t1}
\end{figure}

\begin{figure}[ht!]
    \centering
    \includegraphics[width=0.32\textwidth]{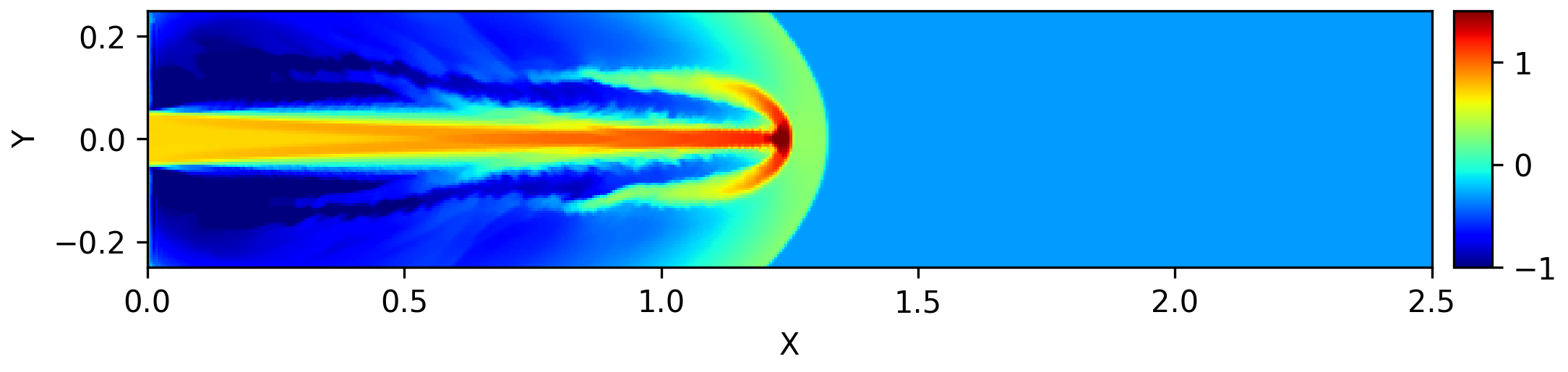}
    \includegraphics[width=0.32\textwidth]{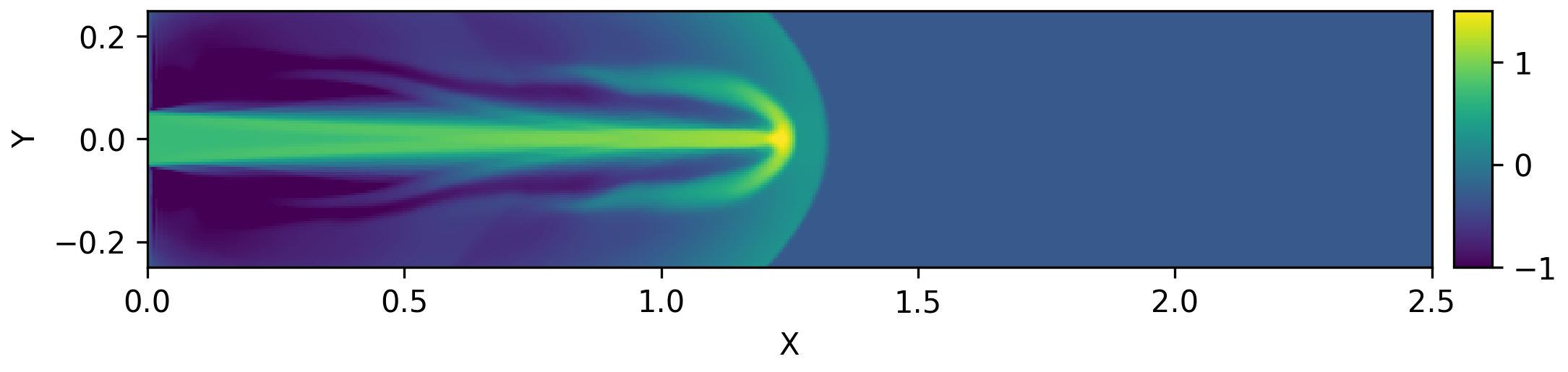}
    \includegraphics[width=0.32\textwidth]{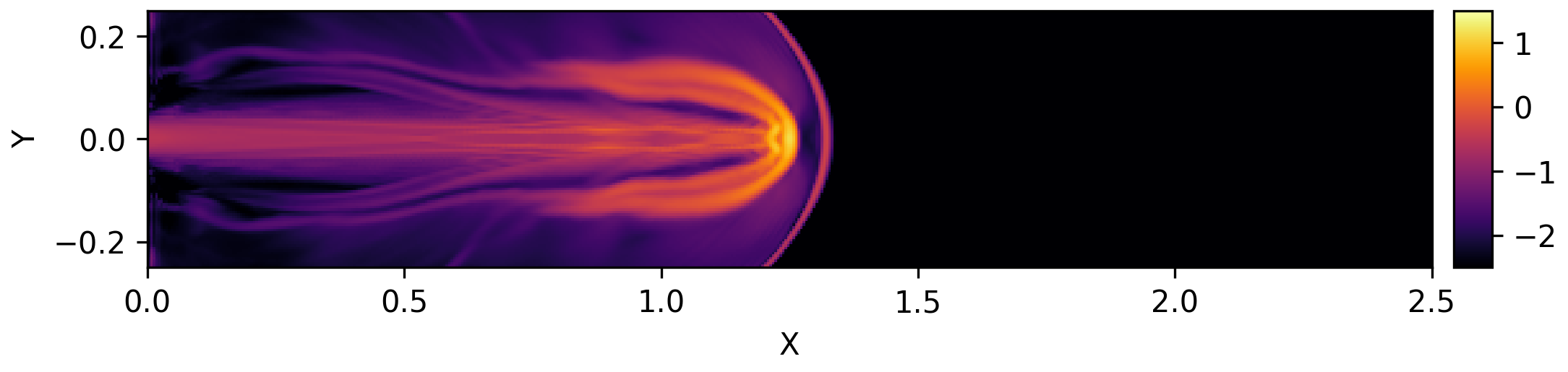} \\
    \includegraphics[width=0.32\textwidth]{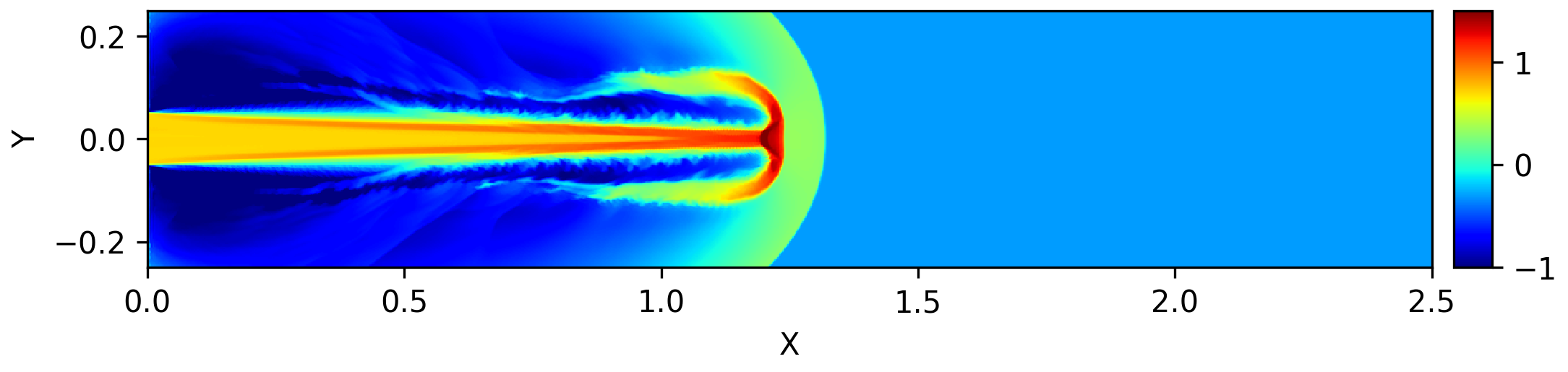}
    \includegraphics[width=0.32\textwidth]{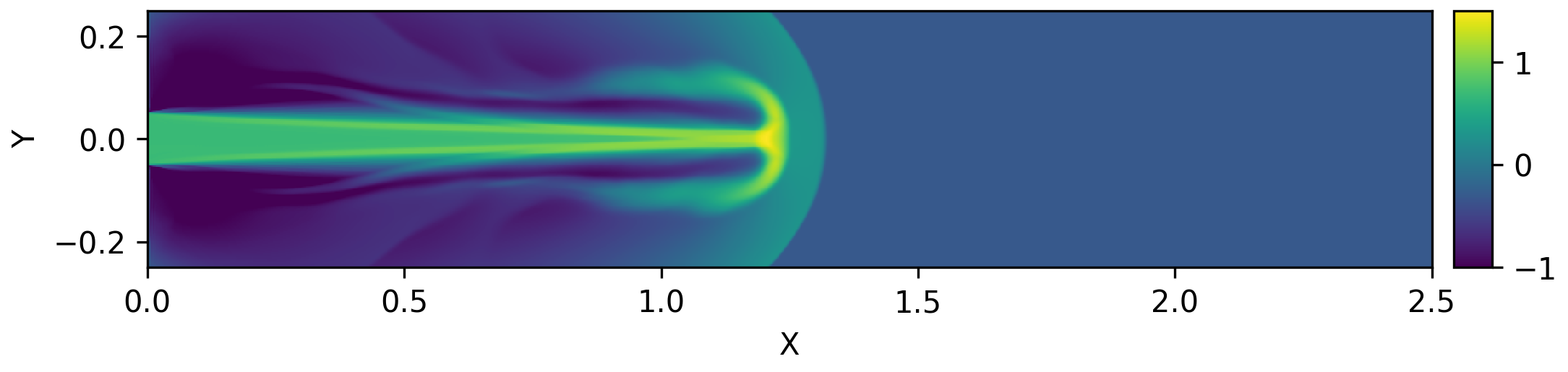}
    \includegraphics[width=0.32\textwidth]{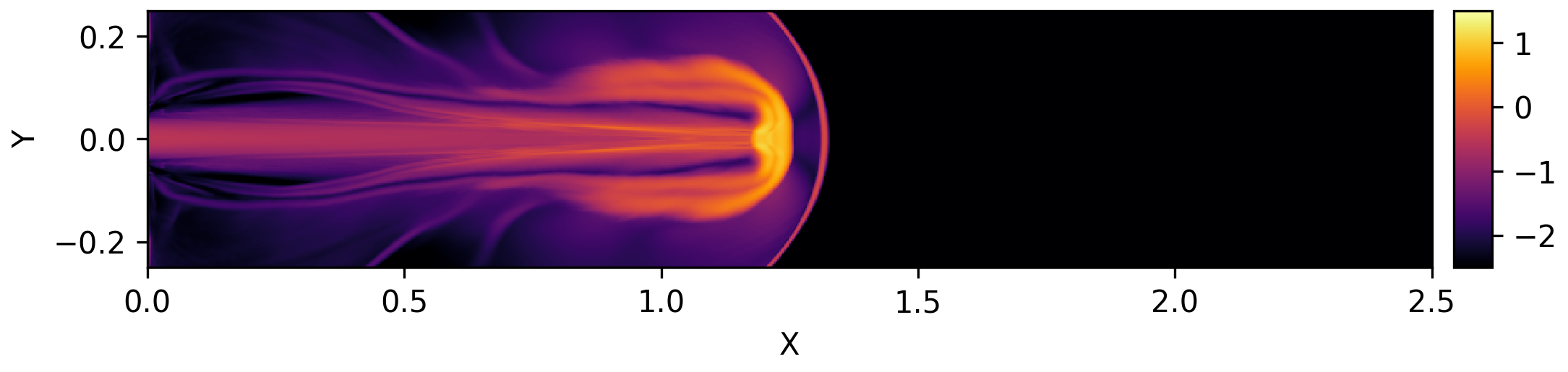}\\
    \includegraphics[width=0.32\textwidth]{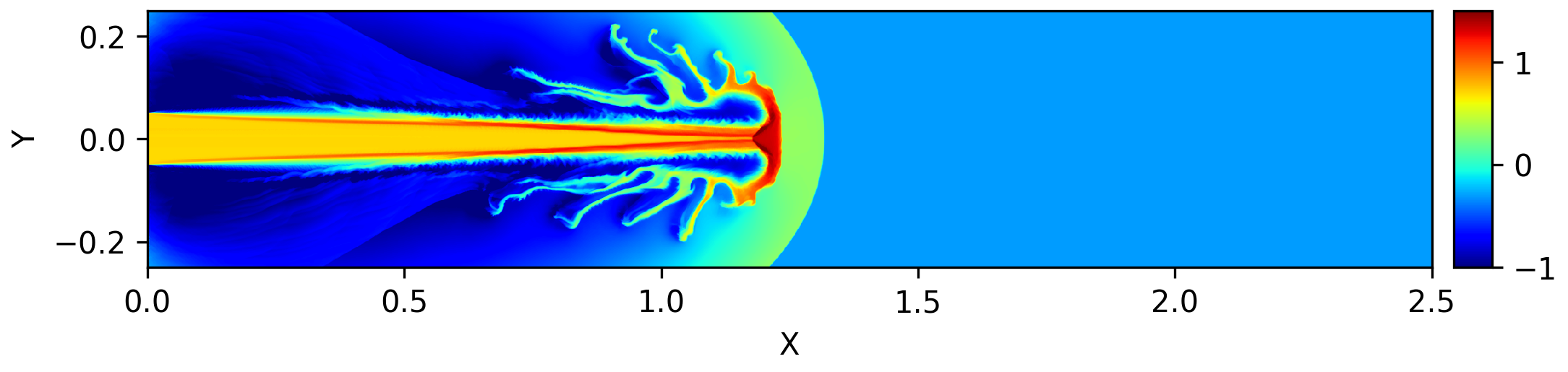}
    \includegraphics[width=0.32\textwidth]{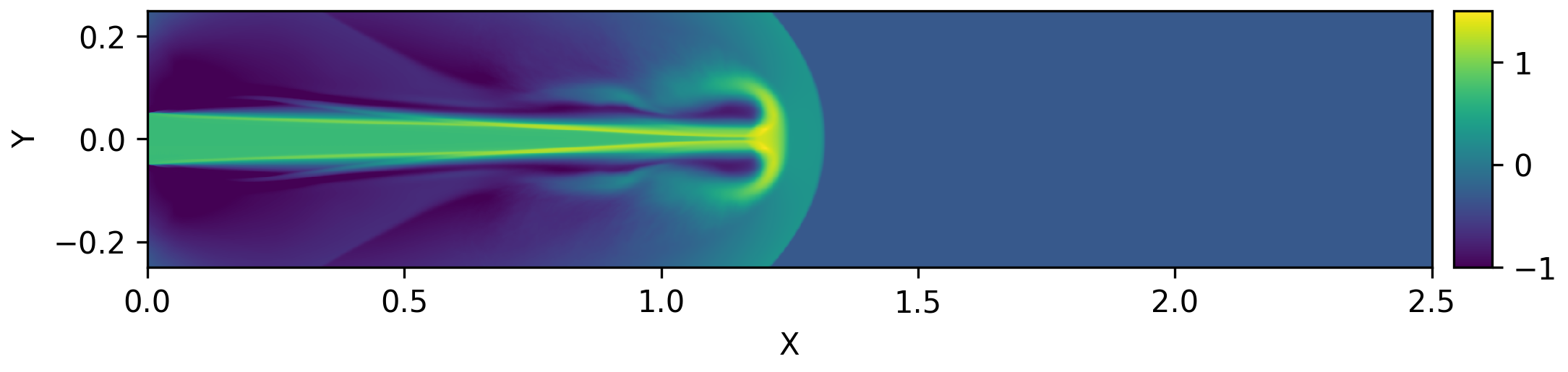}
    \includegraphics[width=0.32\textwidth]{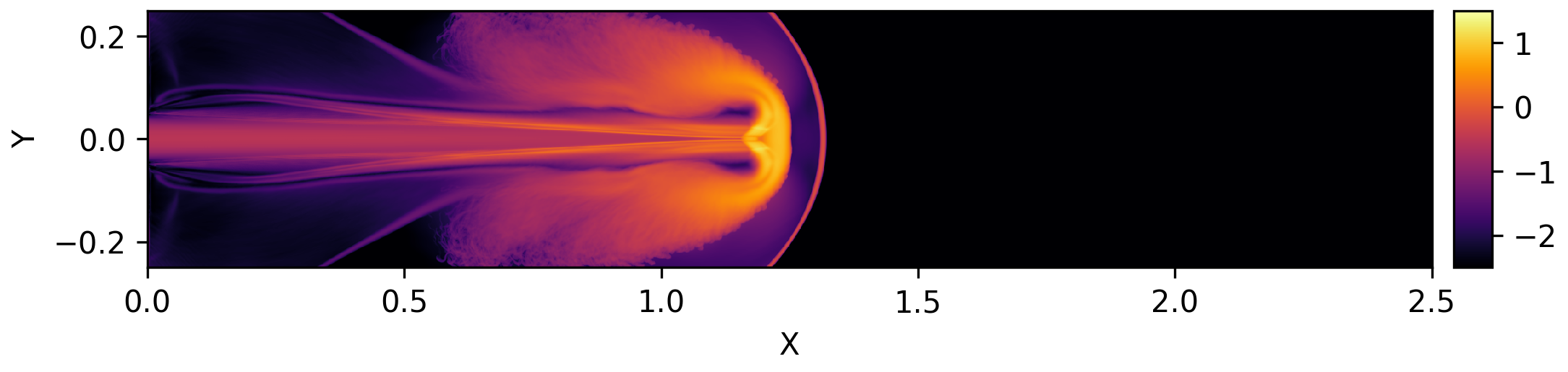}\\
    \includegraphics[width=0.32\textwidth]{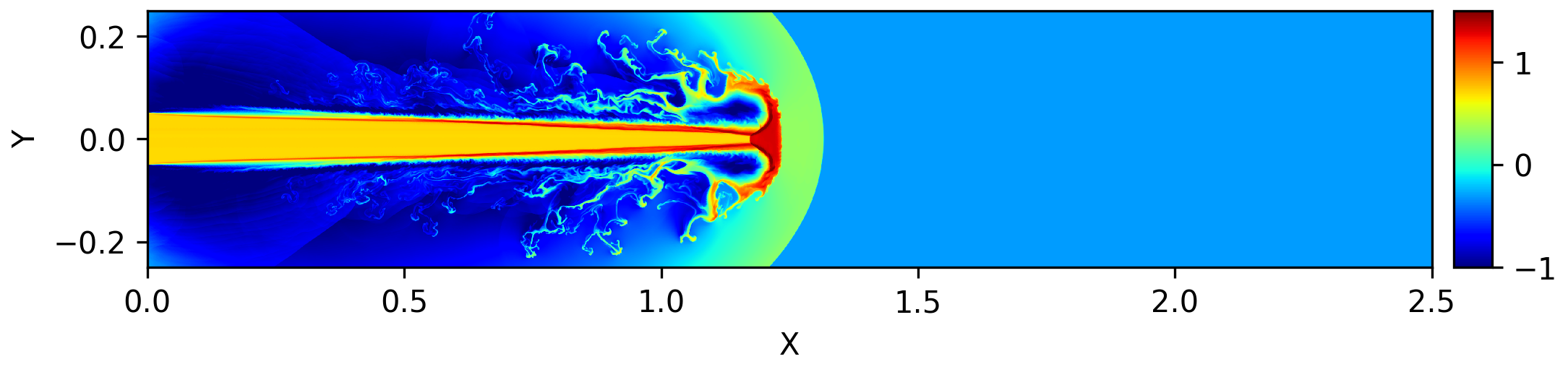}
    \includegraphics[width=0.32\textwidth]{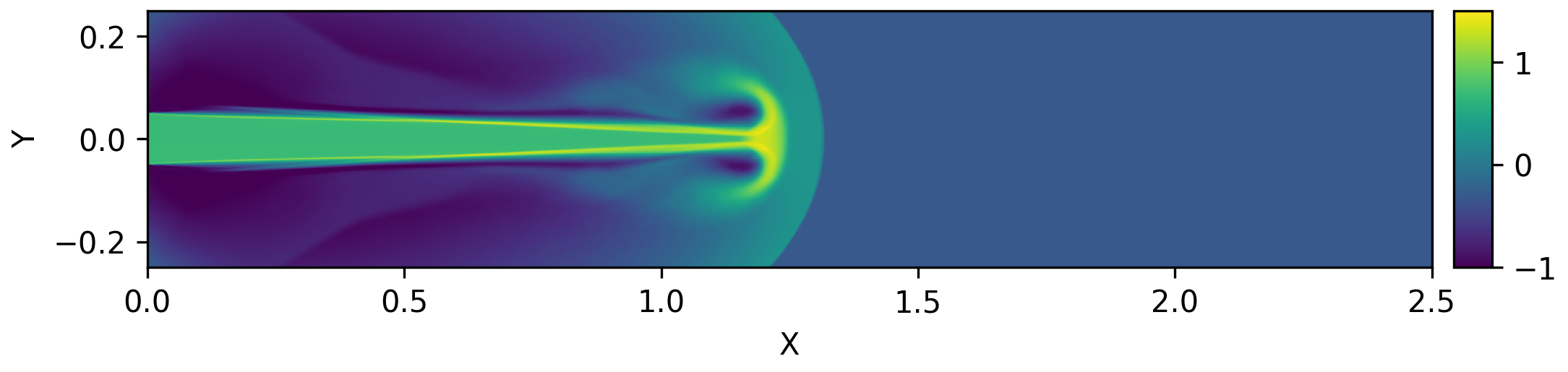}
    \includegraphics[width=0.32\textwidth]{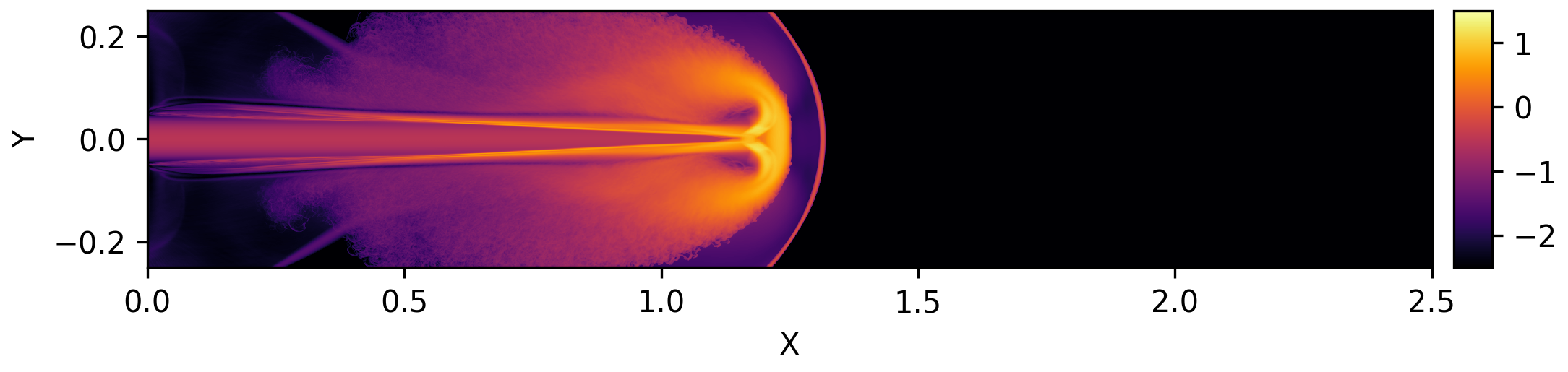}\\
    \caption{
        Single sample and statistical moments (mean and standard deviation for \(M=1000\)) of the density at $t=0.002\,\text{s}$ (left to right).
        The layout displays \(x\)-resolutions \(N=500, 1000, 2000, 4000\) (top to bottom) on the clipped domain \(\overline{\mathcal{D}} = [0,2.5] \times [-0.25, 0.25]\).
        Colorbars are in logarithmic scale.
    }
    \label{fig:jet_convergence_t2}
\end{figure}

\begin{figure}[ht!]
    \centering
    \includegraphics[width=0.32\textwidth]{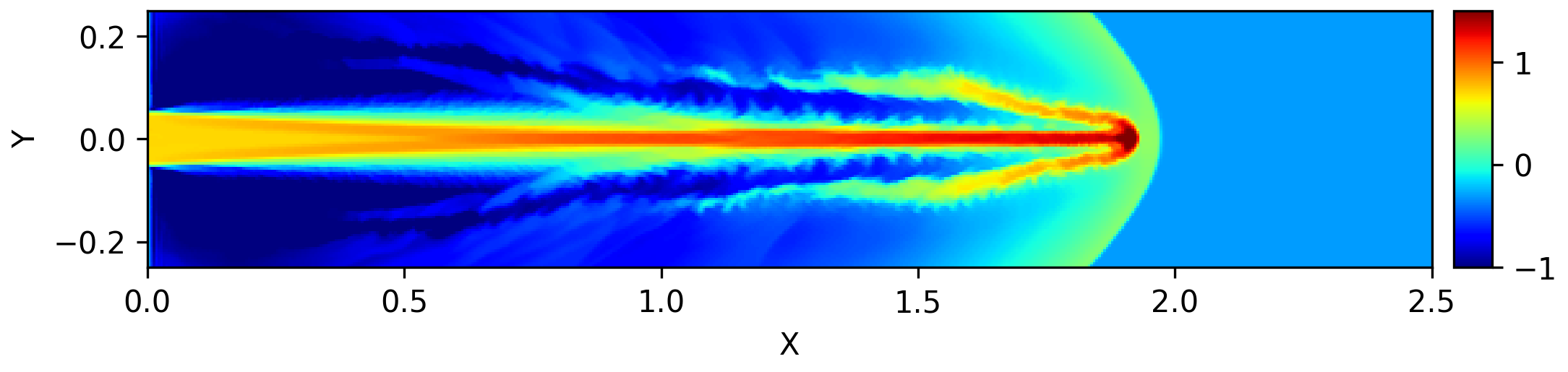}
    \includegraphics[width=0.32\textwidth]{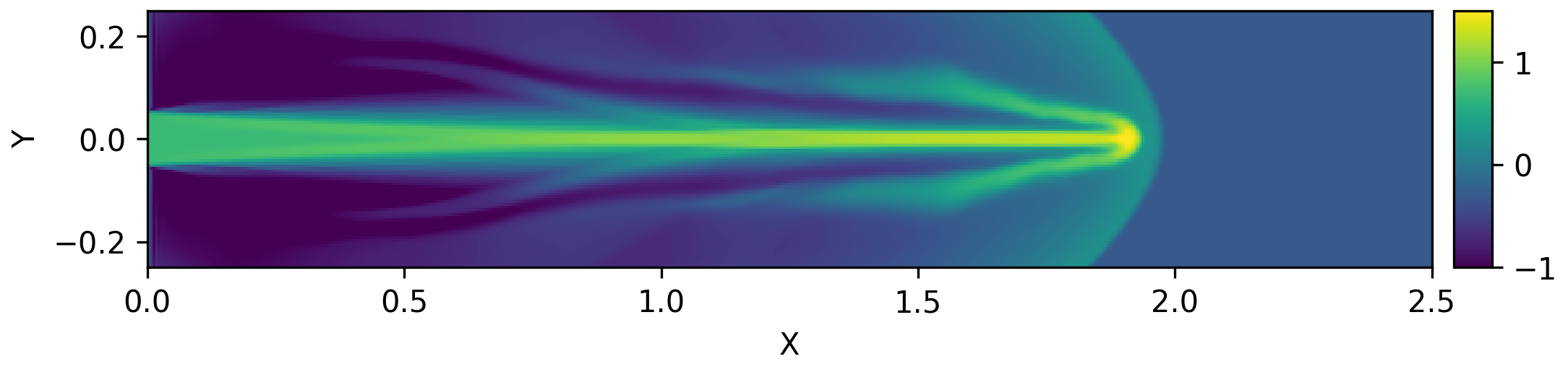}
    \includegraphics[width=0.32\textwidth]{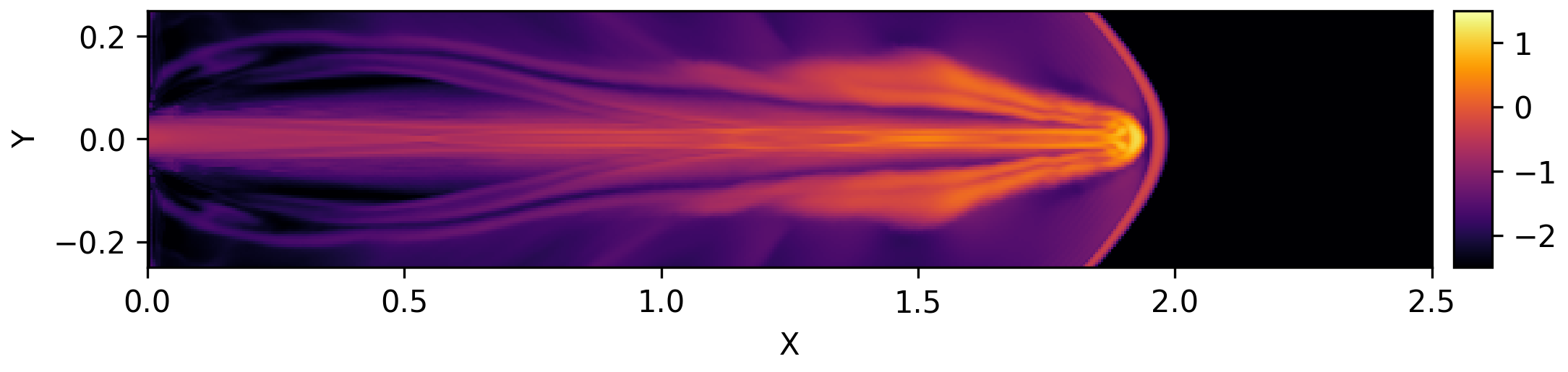} \\
    \includegraphics[width=0.32\textwidth]{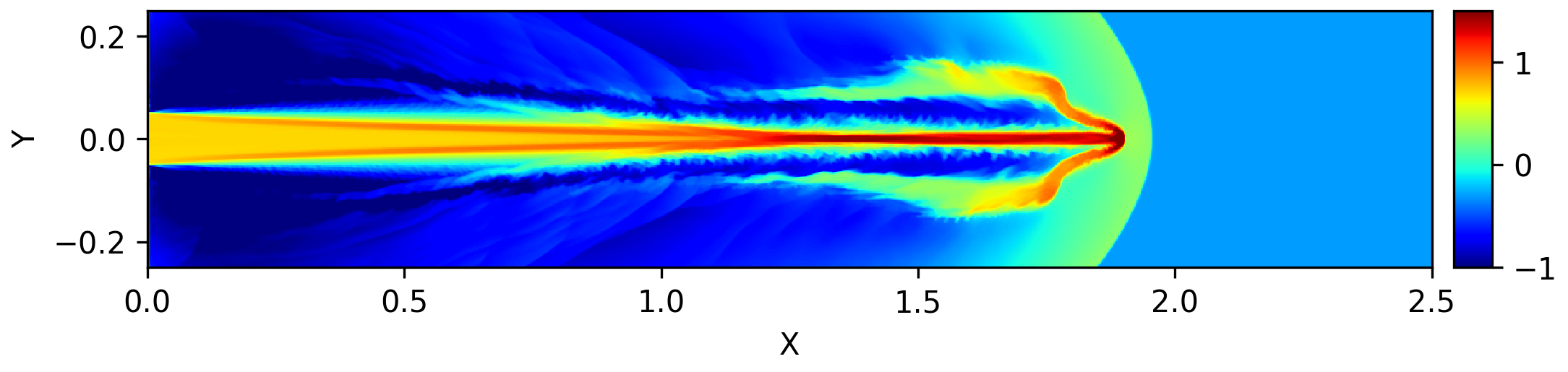}
    \includegraphics[width=0.32\textwidth]{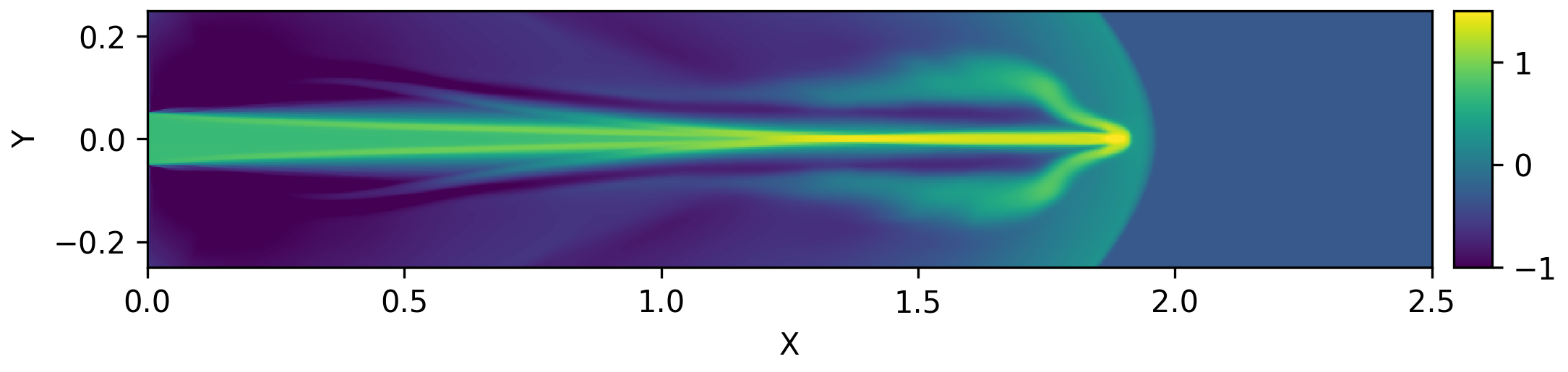}
    \includegraphics[width=0.32\textwidth]{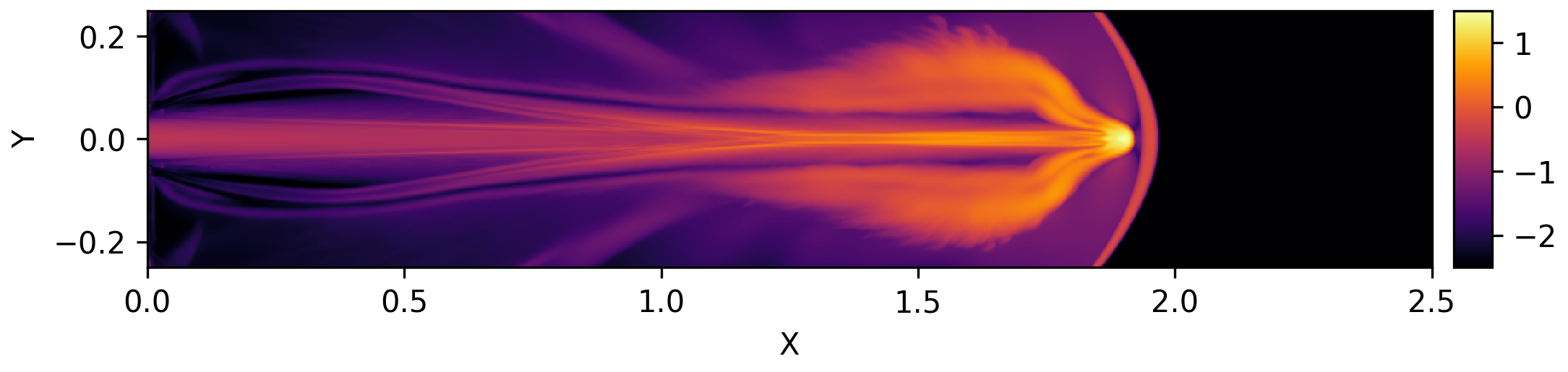}\\
    \includegraphics[width=0.32\textwidth]{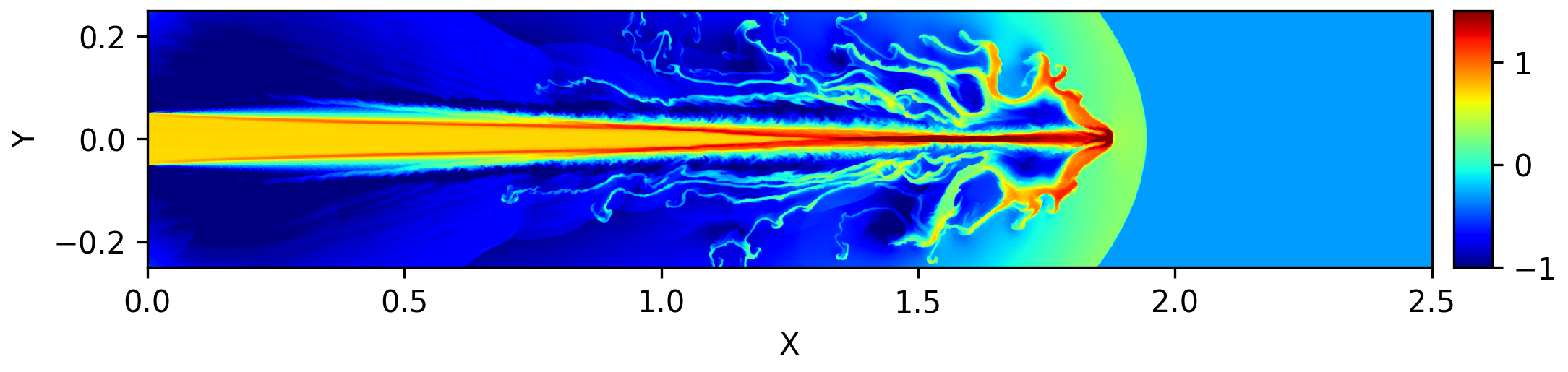}
    \includegraphics[width=0.32\textwidth]{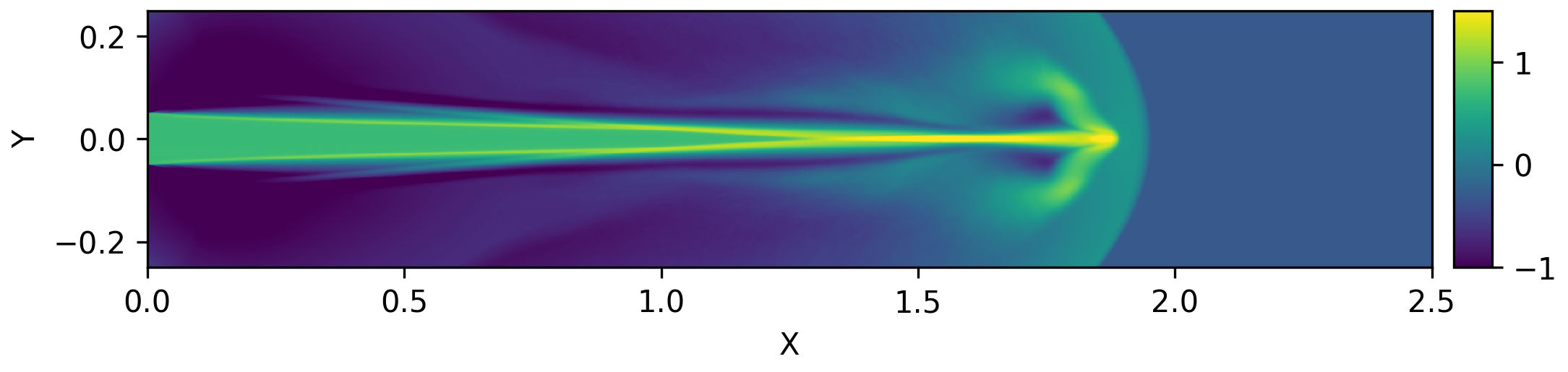}
    \includegraphics[width=0.32\textwidth]{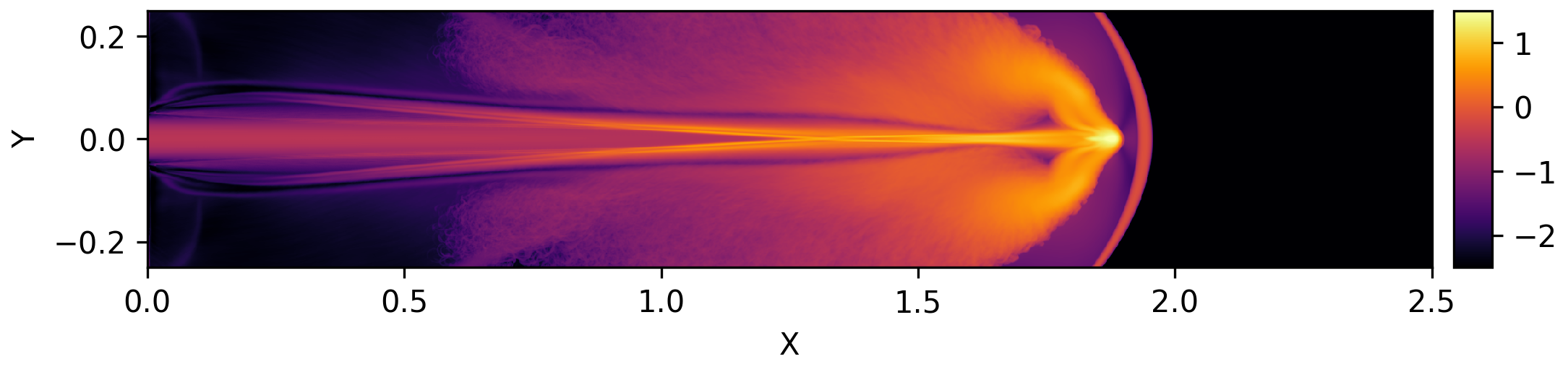}\\
    \includegraphics[width=0.32\textwidth]{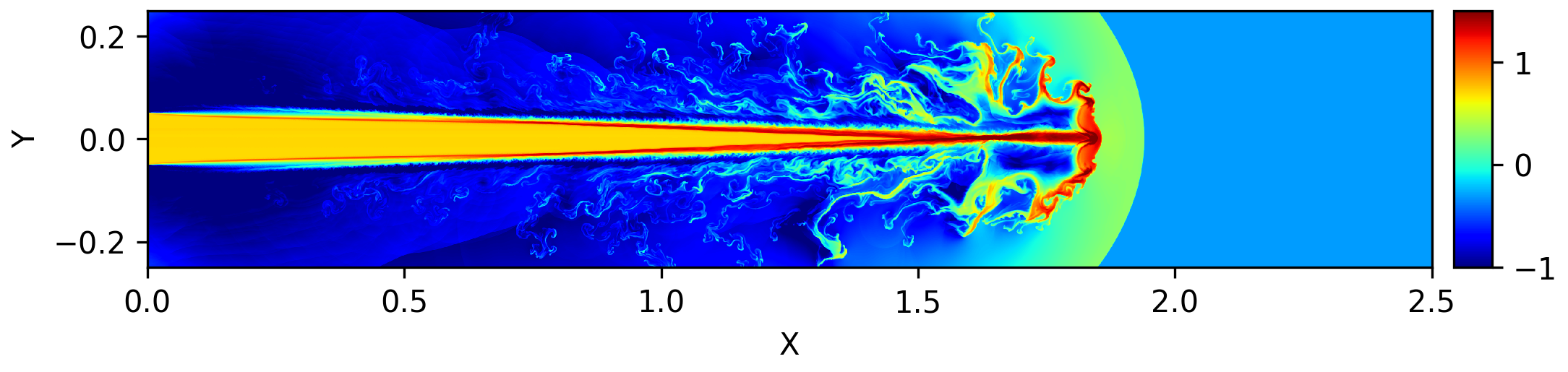}
    \includegraphics[width=0.32\textwidth]{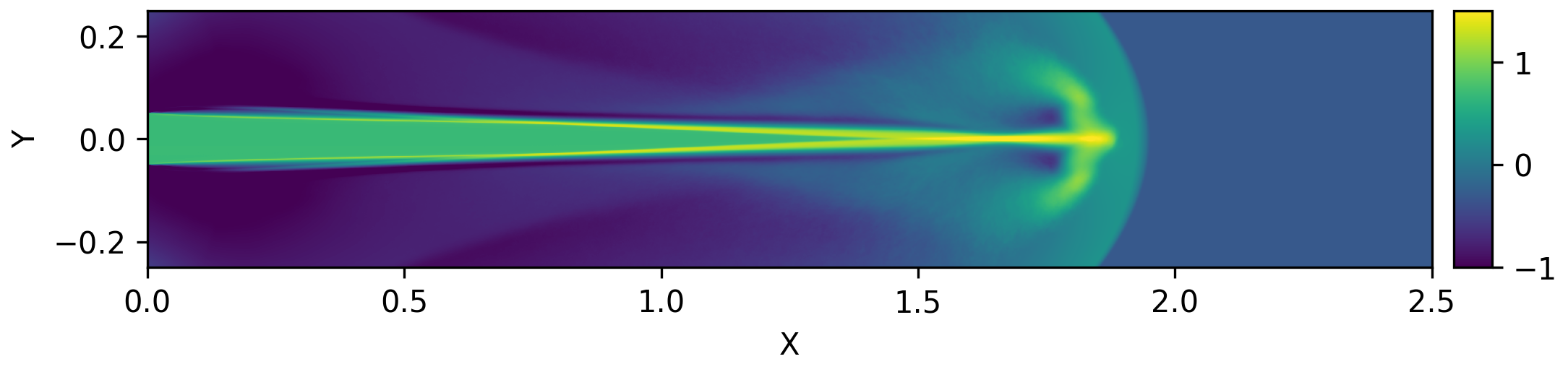}
    \includegraphics[width=0.32\textwidth]{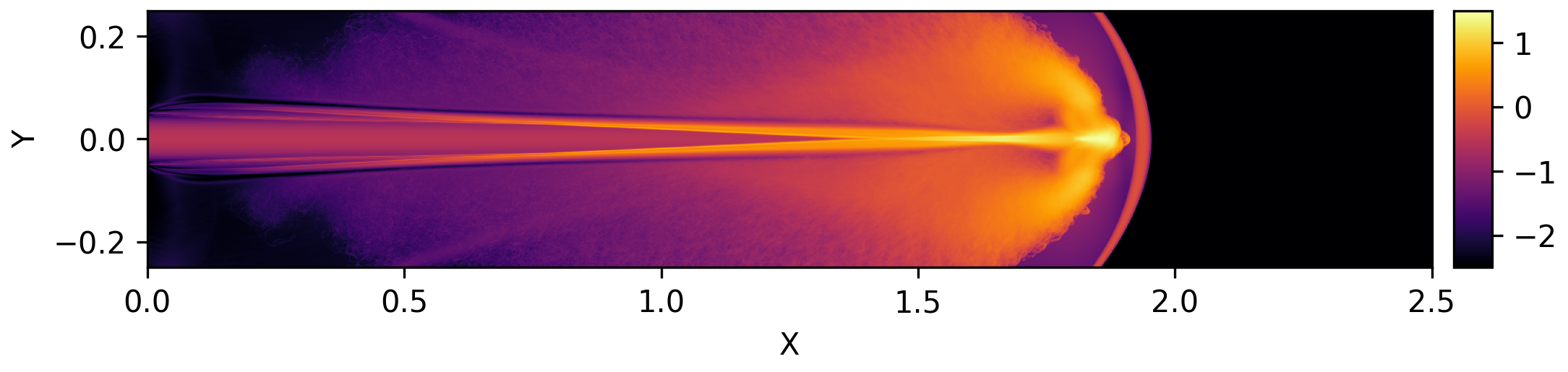}\\
    \caption{
        Single sample and statistical moments (mean and standard deviation for \(M=1000\)) of the density at $t=0.003\,\text{s}$ (left to right).
        The layout displays \(x\)-resolutions \(N=500, 1000, 2000, 4000\) (top to bottom) on the clipped domain \(\overline{\mathcal{D}} = [0,2.5] \times [-0.25, 0.25]\).
        Colorbars are in logarithmic scale.
    }
    \label{fig:jet_convergence_t3}
\end{figure}

\begin{figure}[ht!]
    \centering
    \includegraphics[width=0.32\textwidth]{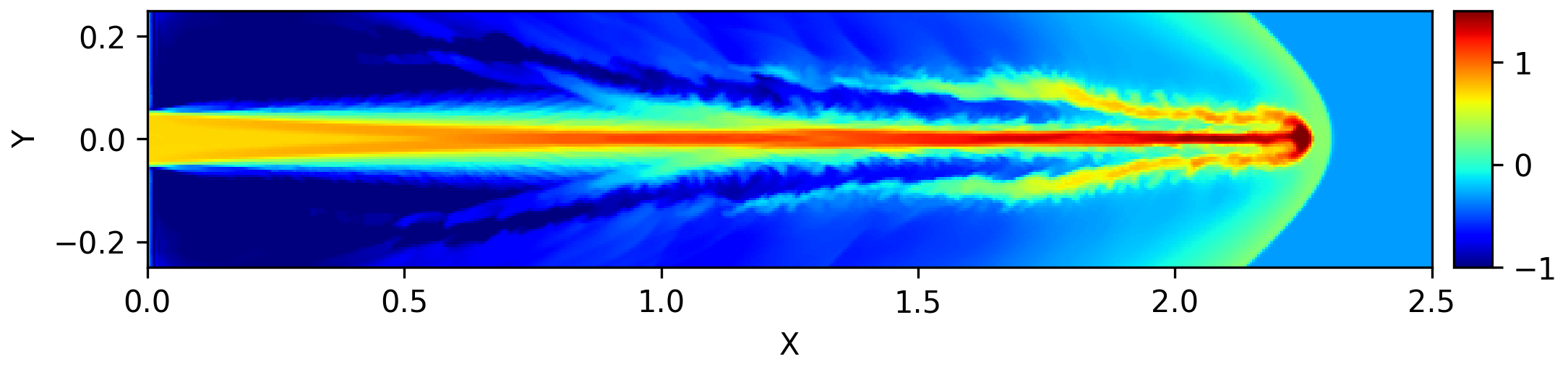}
    \includegraphics[width=0.32\textwidth]{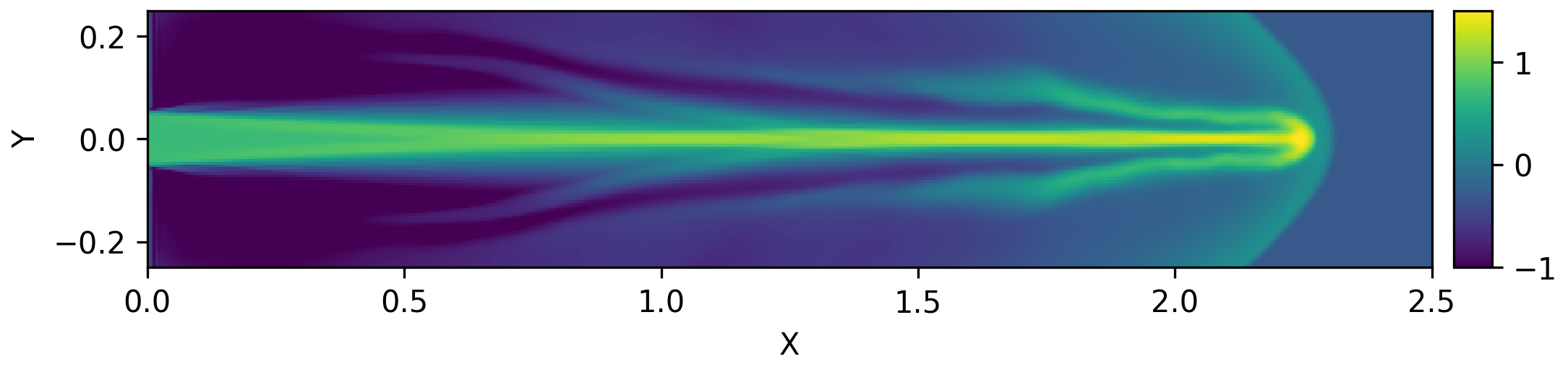}
    \includegraphics[width=0.32\textwidth]{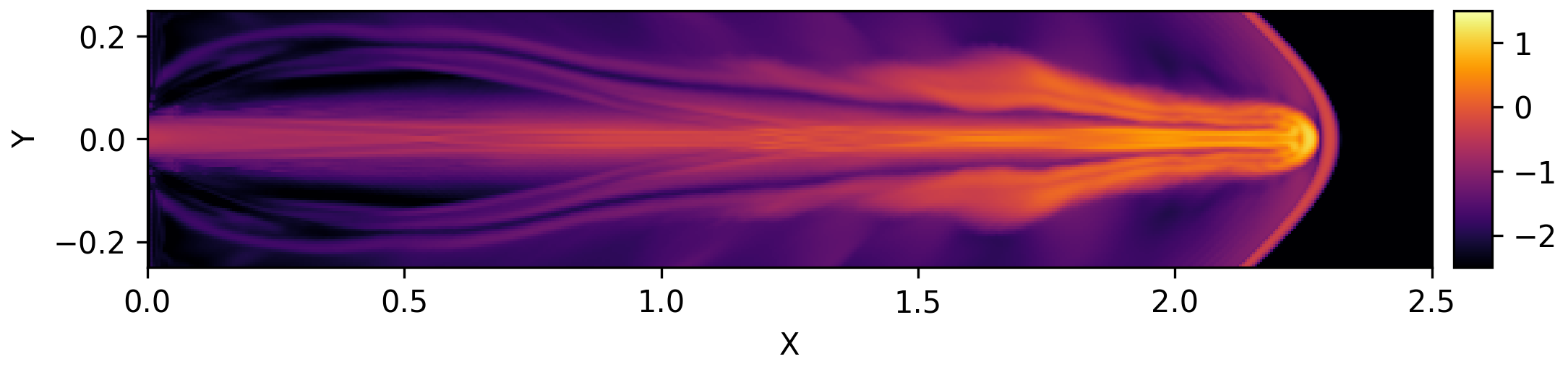} \\
    \includegraphics[width=0.32\textwidth]{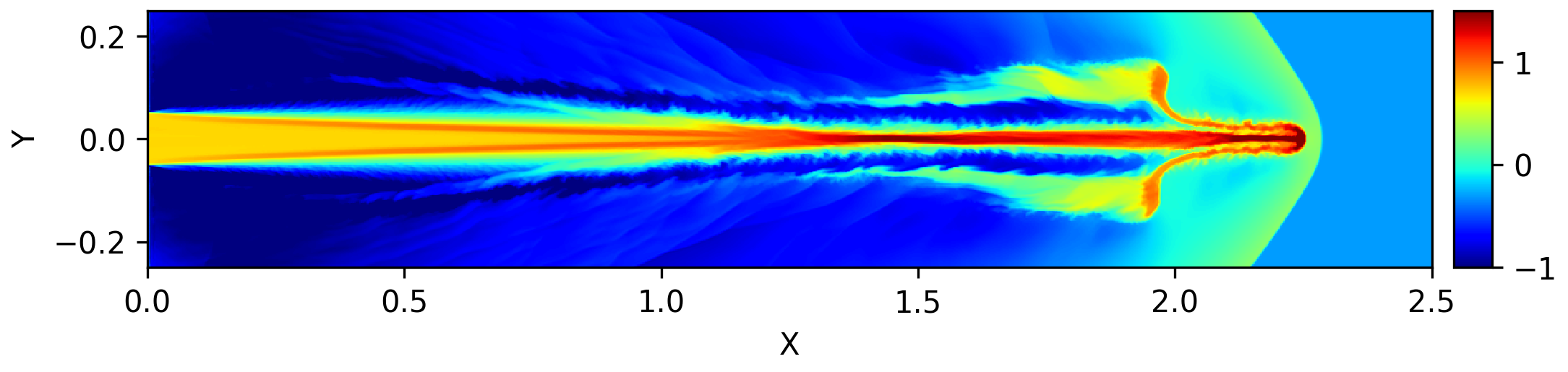}
    \includegraphics[width=0.32\textwidth]{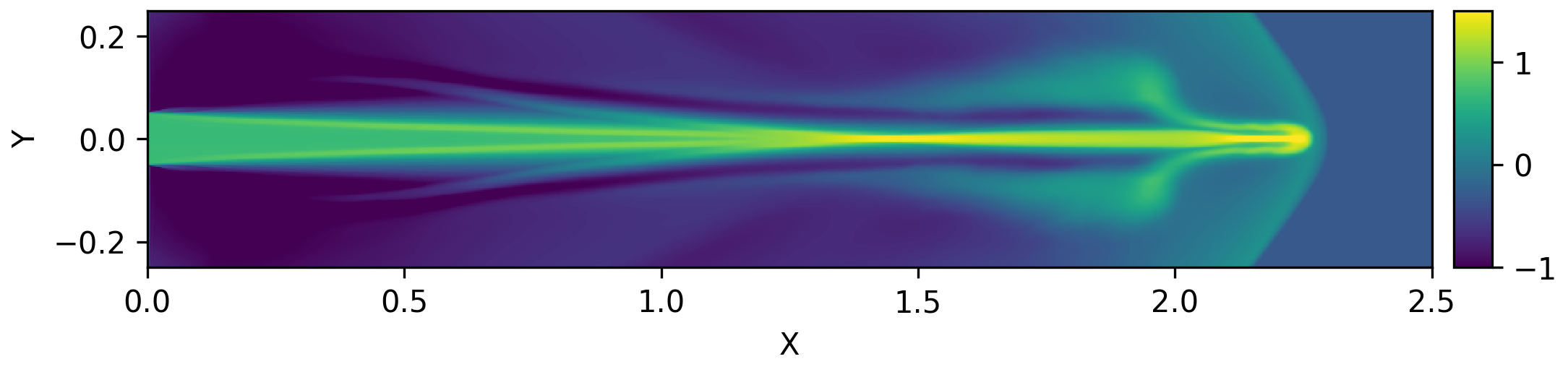}
    \includegraphics[width=0.32\textwidth]{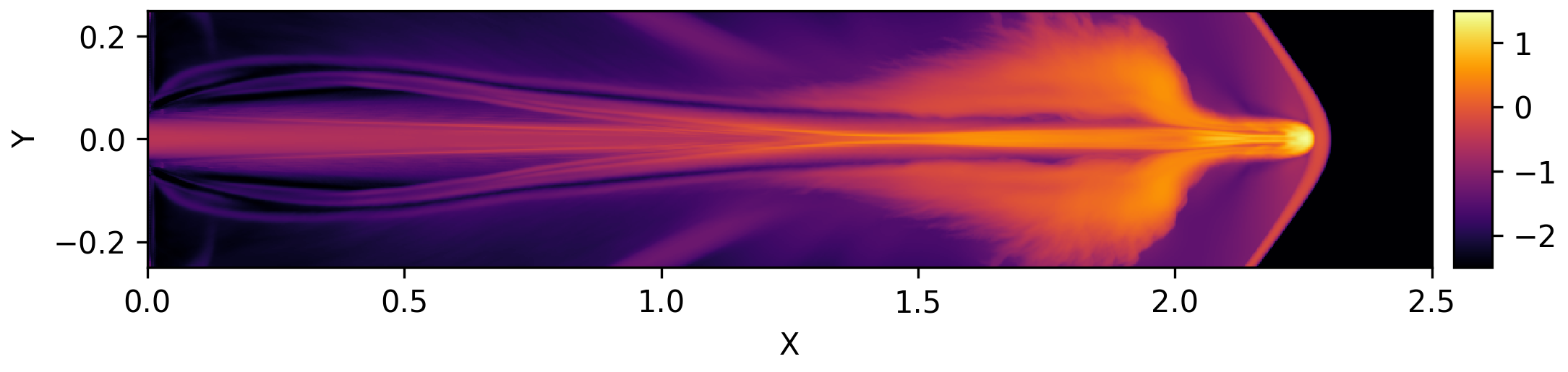}\\
    \includegraphics[width=0.32\textwidth]{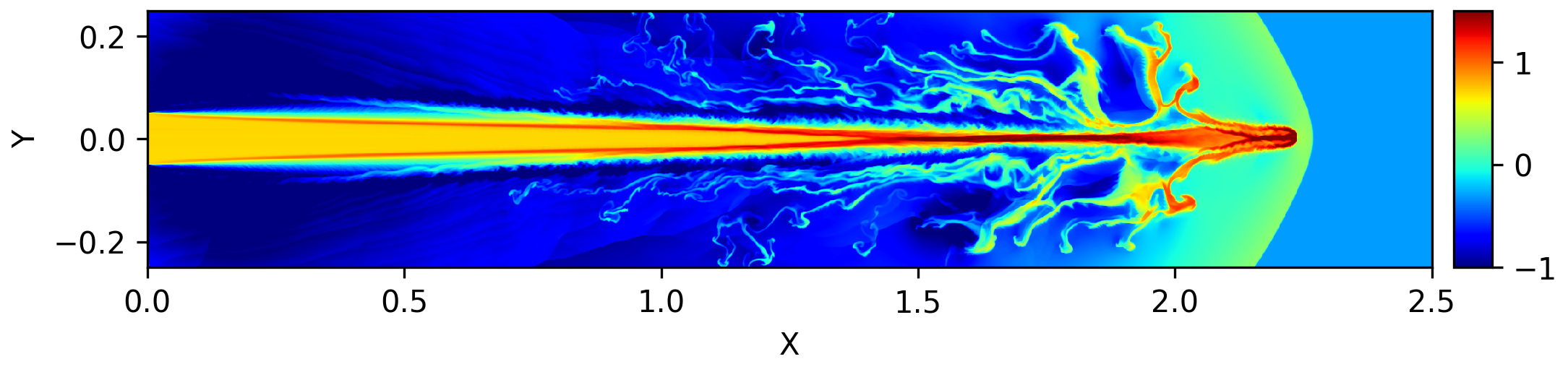}
    \includegraphics[width=0.32\textwidth]{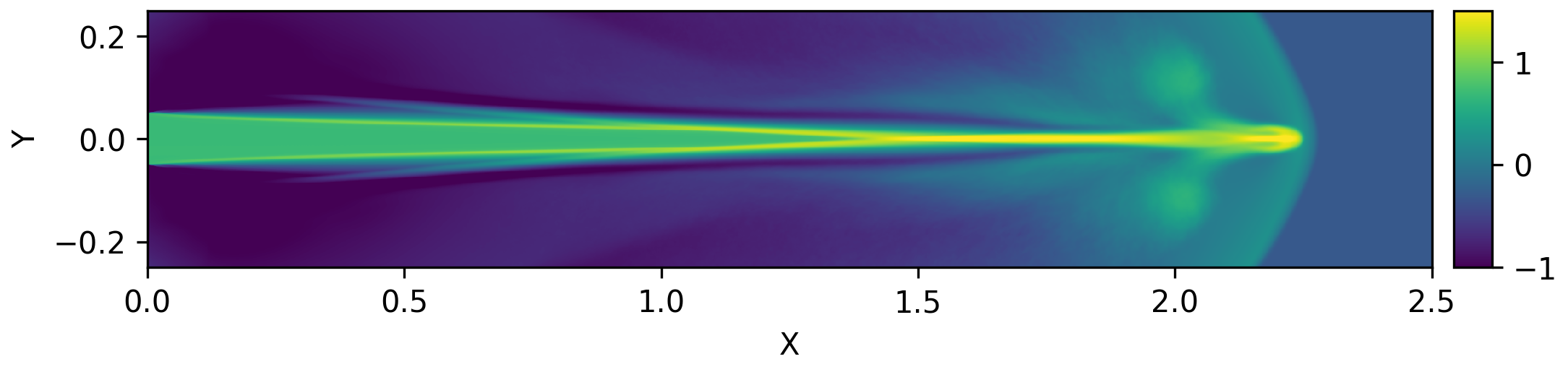}
    \includegraphics[width=0.32\textwidth]{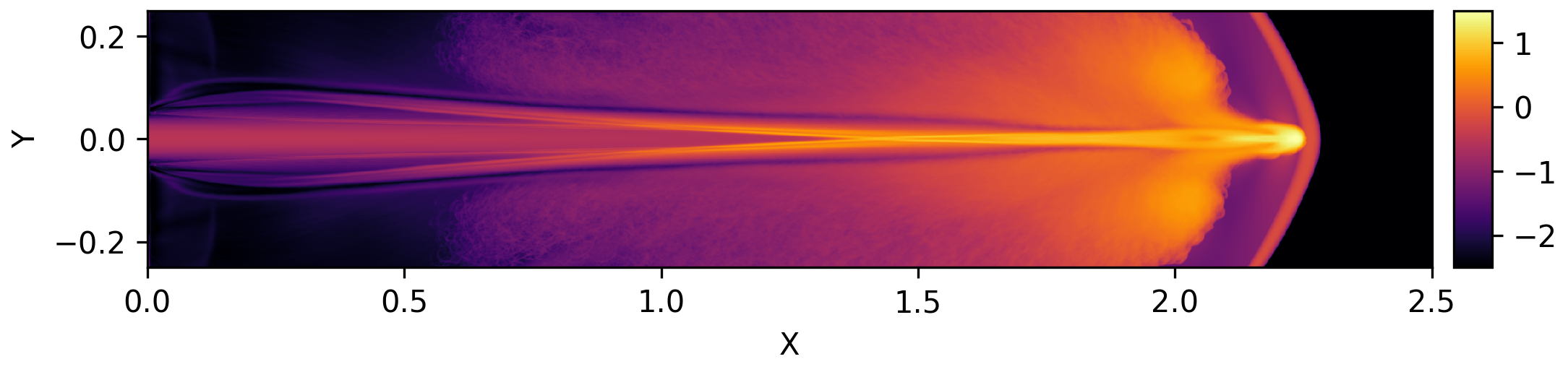}\\
    \includegraphics[width=0.32\textwidth]{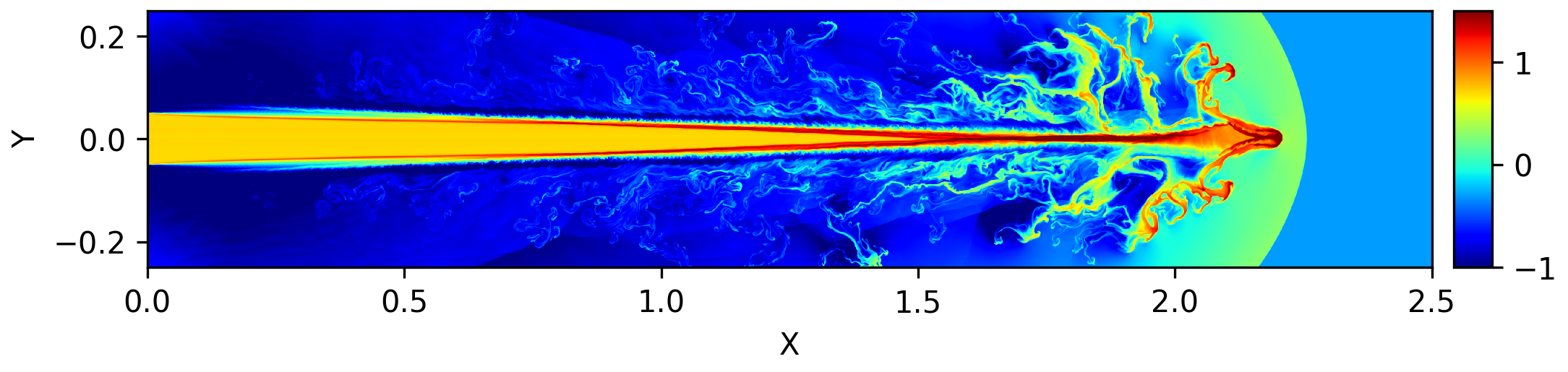}
    \includegraphics[width=0.32\textwidth]{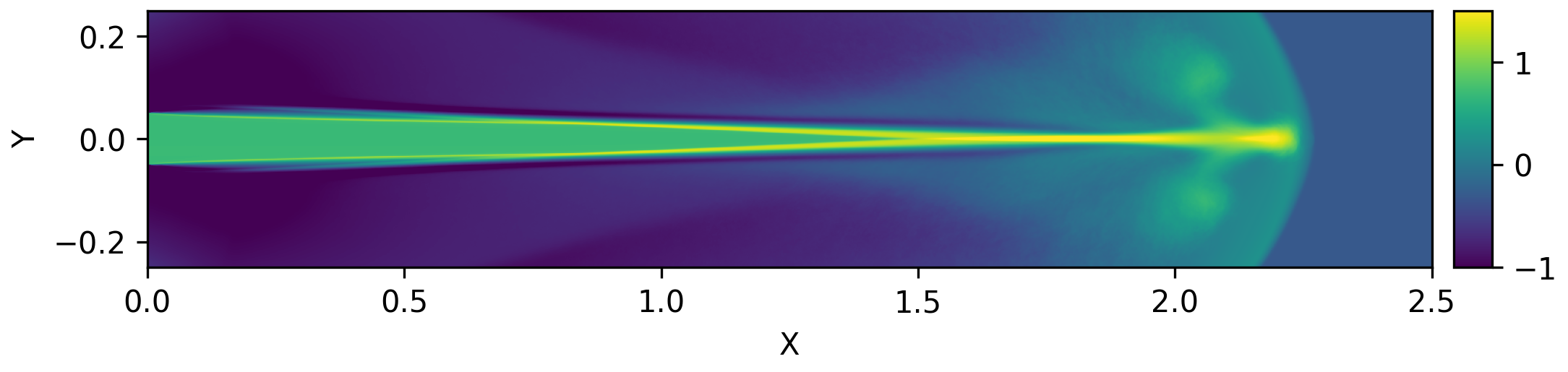}
    \includegraphics[width=0.32\textwidth]{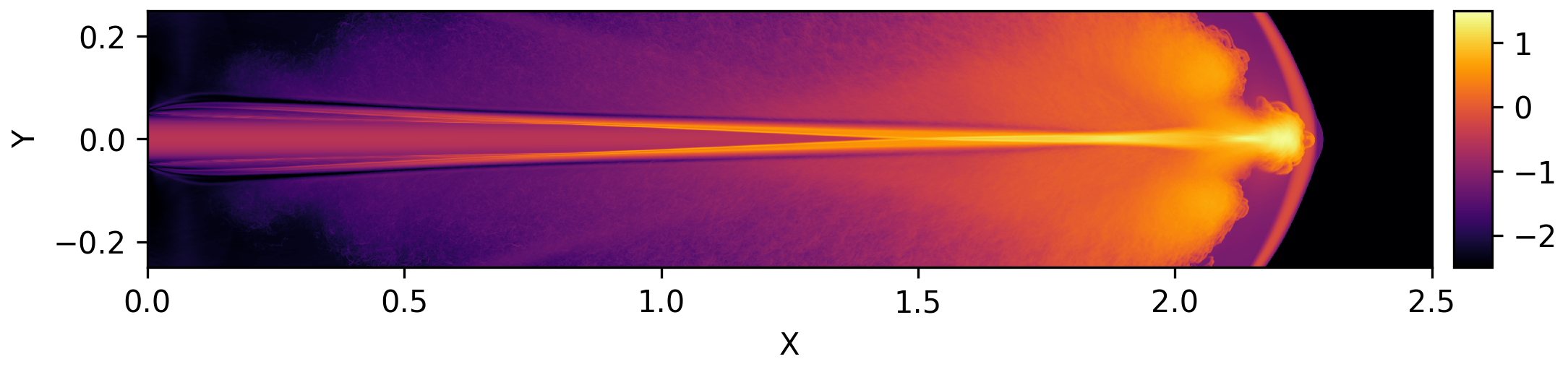}\\
    \caption{
        Single sample and statistical moments (mean and standard deviation for \(M=1000\)) of the density at $t=0.0035\,\text{s}$ (left to right).
        The layout displays \(x\)-resolutions \(N=500, 1000, 2000, 4000\) (top to bottom) on the clipped domain \(\overline{\mathcal{D}} = [0,2.5] \times [-0.25, 0.25]\).
        Colorbars are in logarithmic scale.
    }
    \label{fig:jet_convergence_t35}
\end{figure}

Figure~\ref{fig:overview_final} visualizes two single samples, the mean, and the standard deviation of the computed density evolution at the last timestep $t=0.0035$ for the highest resolution $4000\times 3200$ and $M=1000$ samples.
Figures~\ref{fig:jet_convergence_t1},~\ref{fig:jet_convergence_t2},~\ref{fig:jet_convergence_t3}, and~\ref{fig:jet_convergence_t35} visualize a single representative sample, the mean, and the standard deviation of the computed density evolution for several timesteps, $M =1000$ samples, and all probed resolutions, respectively.
Notably, for the later timesteps, the samples differ visually across increasing grid resolutions, whereas the mean and the standard deviation converge.
Furthermore, as shown in Figures~\ref{fig:jet_convergence_t1},~\ref{fig:jet_convergence_t2},~\ref{fig:jet_convergence_t3}, and~\ref{fig:jet_convergence_t35}, the non-zero variance demonstrates that the statistical solution does not collapse to a Dirac measure but instead captures a range of physically admissible flow states.
This confirms that the numerical approximation is effectively resolving a distribution-valued solution rather than a single deterministic realization; Sections~\ref{subsec:centerline} and~\ref{subsec:pdfheatmap} sharpen this two-dimensional, moment-level picture to a spatially resolved, distribution-level analysis along the jet axis, and Section~\ref{subsec:discussion} quantifies it through the convergence metrics of Section~\ref{sec:weak-strong}.

\subsection{Spatial shock structure along the jet axis}\label{subsec:centerline}

\begin{figure}[ht!]
    \pgfplotsset{
        centerline style/.style={
            width=0.32\textwidth,
            height=0.24\textwidth,
            xmin=-0.1, xmax=2.6,
            ymin=-10, ymax=110,
            ylabel={$\rho(x,0)$},
            xlabel={$x$},
            ytick={-50,0,50,100,150},
            xtick={-0.5,0,0.5,1.0,1.5,2.0,2.5,3.0},
            label style={font=\footnotesize},
            tick label style={font=\footnotesize},
            grid=both,
            yminorgrids=true,
            xminorgrids=true,
            grid=both,
            minor x tick num=4,
            minor x grid style={gray!40},
            major x grid style={black!60},
            minor y tick num=4,
            minor y grid style={gray!40},
            major y grid style={black!60},
            enlarge x limits=false,
        }
    }
    \centering
    \subfloat[$t=0.001$, $N_x=1000$]{%
    \begin{tikzpicture}
    \begin{axis}[centerline style,
                legend to name=centerlineLegend,
                legend style={
                    legend columns=3,
                    font=\footnotesize,
                    draw=darkgray!60!black,
                    fill=white,
                    legend cell align=left,
                    /tikz/every even column/.append style={column sep=0.2cm}
                    }]
    \addplot[name path=up95, draw=none, forget plot, each nth point=4, filter discard warning=false] table[x index=0, y index=7] {./figures/centerline_res400_ck2_t0p0010.dat};
    \addplot[name path=lo05, draw=none, forget plot, each nth point=4, filter discard warning=false] table[x index=0, y index=6] {./figures/centerline_res400_ck2_t0p0010.dat};
    \addplot[blue!25] fill between[of=up95 and lo05];
    \addlegendentry{5\,--\,95\,\% quantile band}
    \addplot[name path=up84, draw=none, forget plot, each nth point=4, filter discard warning=false] table[x index=0, y index=5] {./figures/centerline_res400_ck2_t0p0010.dat};
    \addplot[name path=lo16, draw=none, forget plot, each nth point=4, filter discard warning=false] table[x index=0, y index=4] {./figures/centerline_res400_ck2_t0p0010.dat};
    \addplot[blue!45] fill between[of=up84 and lo16];
    \addlegendentry{16\,--\,84\,\% quantile band}
    \addplot[orange, line width=0.3pt, each nth point=4, filter discard warning=false] table[x index=0, y index=1] {./figures/centerline_res400_ck2_t0p0010.dat};
    \addlegendentry{single sample}
    \addplot[blue, line width=1.1pt, each nth point=4, filter discard warning=false] table[x index=0, y index=2] {./figures/centerline_res400_ck2_t0p0010.dat};
    \addlegendentry{mean $\mathbb{E}[\rho]$}
    \addplot[red, densely dotted, line width=0.9pt, each nth point=4, filter discard warning=false] table[x index=0, y index=3] {./figures/centerline_res400_ck2_t0p0010.dat};
    \addlegendentry{std $\sqrt{\mathrm{Var}[\rho]}$}
    \end{axis}
    \end{tikzpicture}
    }%
    \subfloat[$t=0.001$, $N_x=2000$]{%
    \begin{tikzpicture}
    \begin{axis}[centerline style]
    \addplot[name path=up95, draw=none, forget plot, each nth point=4, filter discard warning=false] table[x index=0, y index=7] {./figures/centerline_res800_ck2_t0p0010.dat};
    \addplot[name path=lo05, draw=none, forget plot, each nth point=4, filter discard warning=false] table[x index=0, y index=6] {./figures/centerline_res800_ck2_t0p0010.dat};
    \addplot[blue!25, forget plot] fill between[of=up95 and lo05];
    \addplot[name path=up84, draw=none, forget plot, each nth point=4, filter discard warning=false] table[x index=0, y index=5] {./figures/centerline_res800_ck2_t0p0010.dat};
    \addplot[name path=lo16, draw=none, forget plot, each nth point=4, filter discard warning=false] table[x index=0, y index=4] {./figures/centerline_res800_ck2_t0p0010.dat};
    \addplot[blue!45, forget plot] fill between[of=up84 and lo16];
    \addplot[orange, line width=0.3pt, forget plot, each nth point=4, filter discard warning=false] table[x index=0, y index=1] {./figures/centerline_res800_ck2_t0p0010.dat};
    \addplot[blue, line width=1.1pt, forget plot, each nth point=4, filter discard warning=false] table[x index=0, y index=2] {./figures/centerline_res800_ck2_t0p0010.dat};
    \addplot[red, densely dotted, line width=0.9pt, forget plot, each nth point=4, filter discard warning=false] table[x index=0, y index=3] {./figures/centerline_res800_ck2_t0p0010.dat};
    \end{axis}
    \end{tikzpicture}
    }%
    \subfloat[$t=0.001$, $N_x=4000$]{%
    \begin{tikzpicture}
    \begin{axis}[centerline style]
    \addplot[name path=up95, draw=none, forget plot, each nth point=4, filter discard warning=false] table[x index=0, y index=7] {./figures/centerline_res1600_ck2_t0p0010.dat};
    \addplot[name path=lo05, draw=none, forget plot, each nth point=4, filter discard warning=false] table[x index=0, y index=6] {./figures/centerline_res1600_ck2_t0p0010.dat};
    \addplot[blue!25, forget plot] fill between[of=up95 and lo05];
    \addplot[name path=up84, draw=none, forget plot, each nth point=4, filter discard warning=false] table[x index=0, y index=5] {./figures/centerline_res1600_ck2_t0p0010.dat};
    \addplot[name path=lo16, draw=none, forget plot, each nth point=4, filter discard warning=false] table[x index=0, y index=4] {./figures/centerline_res1600_ck2_t0p0010.dat};
    \addplot[blue!45, forget plot] fill between[of=up84 and lo16];
    \addplot[orange, line width=0.3pt, forget plot, each nth point=4, filter discard warning=false] table[x index=0, y index=1] {./figures/centerline_res1600_ck2_t0p0010.dat};
    \addplot[blue, line width=1.1pt, forget plot, each nth point=4, filter discard warning=false] table[x index=0, y index=2] {./figures/centerline_res1600_ck2_t0p0010.dat};
    \addplot[red, densely dotted, line width=0.9pt, forget plot, each nth point=4, filter discard warning=false] table[x index=0, y index=3] {./figures/centerline_res1600_ck2_t0p0010.dat};
    \end{axis}
    \end{tikzpicture}
    }
    \\[-0.3em]
    \subfloat[$t=0.002$, $N_x=1000$]{%
    \begin{tikzpicture}
    \begin{axis}[centerline style]
    \addplot[name path=up95, draw=none, forget plot, each nth point=4, filter discard warning=false] table[x index=0, y index=7] {./figures/centerline_res400_ck4_t0p0020.dat};
    \addplot[name path=lo05, draw=none, forget plot, each nth point=4, filter discard warning=false] table[x index=0, y index=6] {./figures/centerline_res400_ck4_t0p0020.dat};
    \addplot[blue!25, forget plot] fill between[of=up95 and lo05];
    \addplot[name path=up84, draw=none, forget plot, each nth point=4, filter discard warning=false] table[x index=0, y index=5] {./figures/centerline_res400_ck4_t0p0020.dat};
    \addplot[name path=lo16, draw=none, forget plot, each nth point=4, filter discard warning=false] table[x index=0, y index=4] {./figures/centerline_res400_ck4_t0p0020.dat};
    \addplot[blue!45, forget plot] fill between[of=up84 and lo16];
    \addplot[orange, line width=0.3pt, forget plot, each nth point=4, filter discard warning=false] table[x index=0, y index=1] {./figures/centerline_res400_ck4_t0p0020.dat};
    \addplot[blue, line width=1.1pt, forget plot, each nth point=4, filter discard warning=false] table[x index=0, y index=2] {./figures/centerline_res400_ck4_t0p0020.dat};
    \addplot[red, densely dotted, line width=0.9pt, forget plot, each nth point=4, filter discard warning=false] table[x index=0, y index=3] {./figures/centerline_res400_ck4_t0p0020.dat};
    \end{axis}
    \end{tikzpicture}
    }%
    \subfloat[$t=0.002$, $N_x=2000$]{%
    \begin{tikzpicture}
    \begin{axis}[centerline style]
    \addplot[name path=up95, draw=none, forget plot, each nth point=4, filter discard warning=false] table[x index=0, y index=7] {./figures/centerline_res800_ck4_t0p0020.dat};
    \addplot[name path=lo05, draw=none, forget plot, each nth point=4, filter discard warning=false] table[x index=0, y index=6] {./figures/centerline_res800_ck4_t0p0020.dat};
    \addplot[blue!25, forget plot] fill between[of=up95 and lo05];
    \addplot[name path=up84, draw=none, forget plot, each nth point=4, filter discard warning=false] table[x index=0, y index=5] {./figures/centerline_res800_ck4_t0p0020.dat};
    \addplot[name path=lo16, draw=none, forget plot, each nth point=4, filter discard warning=false] table[x index=0, y index=4] {./figures/centerline_res800_ck4_t0p0020.dat};
    \addplot[blue!45, forget plot] fill between[of=up84 and lo16];
    \addplot[orange, line width=0.3pt, forget plot, each nth point=4, filter discard warning=false] table[x index=0, y index=1] {./figures/centerline_res800_ck4_t0p0020.dat};
    \addplot[blue, line width=1.1pt, forget plot, each nth point=4, filter discard warning=false] table[x index=0, y index=2] {./figures/centerline_res800_ck4_t0p0020.dat};
    \addplot[red, densely dotted, line width=0.9pt, forget plot, each nth point=4, filter discard warning=false] table[x index=0, y index=3] {./figures/centerline_res800_ck4_t0p0020.dat};
    \end{axis}
    \end{tikzpicture}
    }%
    \subfloat[$t=0.002$, $N_x=4000$]{%
    \begin{tikzpicture}
    \begin{axis}[centerline style]
    \addplot[name path=up95, draw=none, forget plot, each nth point=4, filter discard warning=false] table[x index=0, y index=7] {./figures/centerline_res1600_ck4_t0p0020.dat};
    \addplot[name path=lo05, draw=none, forget plot, each nth point=4, filter discard warning=false] table[x index=0, y index=6] {./figures/centerline_res1600_ck4_t0p0020.dat};
    \addplot[blue!25, forget plot] fill between[of=up95 and lo05];
    \addplot[name path=up84, draw=none, forget plot, each nth point=4, filter discard warning=false] table[x index=0, y index=5] {./figures/centerline_res1600_ck4_t0p0020.dat};
    \addplot[name path=lo16, draw=none, forget plot, each nth point=4, filter discard warning=false] table[x index=0, y index=4] {./figures/centerline_res1600_ck4_t0p0020.dat};
    \addplot[blue!45, forget plot] fill between[of=up84 and lo16];
    \addplot[orange, line width=0.3pt, forget plot, each nth point=4, filter discard warning=false] table[x index=0, y index=1] {./figures/centerline_res1600_ck4_t0p0020.dat};
    \addplot[blue, line width=1.1pt, forget plot, each nth point=4, filter discard warning=false] table[x index=0, y index=2] {./figures/centerline_res1600_ck4_t0p0020.dat};
    \addplot[red, densely dotted, line width=0.9pt, forget plot, each nth point=4, filter discard warning=false] table[x index=0, y index=3] {./figures/centerline_res1600_ck4_t0p0020.dat};
    \end{axis}
    \end{tikzpicture}
    }
    \\[-0.3em]
    \subfloat[$t=0.003$, $N_x=1000$]{%
    \begin{tikzpicture}
    \begin{axis}[centerline style]
    \addplot[name path=up95, draw=none, forget plot, each nth point=4, filter discard warning=false] table[x index=0, y index=7] {./figures/centerline_res400_ck6_t0p0030.dat};
    \addplot[name path=lo05, draw=none, forget plot, each nth point=4, filter discard warning=false] table[x index=0, y index=6] {./figures/centerline_res400_ck6_t0p0030.dat};
    \addplot[blue!25, forget plot] fill between[of=up95 and lo05];
    \addplot[name path=up84, draw=none, forget plot, each nth point=4, filter discard warning=false] table[x index=0, y index=5] {./figures/centerline_res400_ck6_t0p0030.dat};
    \addplot[name path=lo16, draw=none, forget plot, each nth point=4, filter discard warning=false] table[x index=0, y index=4] {./figures/centerline_res400_ck6_t0p0030.dat};
    \addplot[blue!45, forget plot] fill between[of=up84 and lo16];
    \addplot[orange, line width=0.3pt, forget plot, each nth point=4, filter discard warning=false] table[x index=0, y index=1] {./figures/centerline_res400_ck6_t0p0030.dat};
    \addplot[blue, line width=1.1pt, forget plot, each nth point=4, filter discard warning=false] table[x index=0, y index=2] {./figures/centerline_res400_ck6_t0p0030.dat};
    \addplot[red, densely dotted, line width=0.9pt, forget plot, each nth point=4, filter discard warning=false] table[x index=0, y index=3] {./figures/centerline_res400_ck6_t0p0030.dat};
    \end{axis}
    \end{tikzpicture}
    }%
    \subfloat[$t=0.003$, $N_x=2000$]{%
    \begin{tikzpicture}
    \begin{axis}[centerline style]
    \addplot[name path=up95, draw=none, forget plot, each nth point=4, filter discard warning=false] table[x index=0, y index=7] {./figures/centerline_res800_ck6_t0p0030.dat};
    \addplot[name path=lo05, draw=none, forget plot, each nth point=4, filter discard warning=false] table[x index=0, y index=6] {./figures/centerline_res800_ck6_t0p0030.dat};
    \addplot[blue!25, forget plot] fill between[of=up95 and lo05];
    \addplot[name path=up84, draw=none, forget plot, each nth point=4, filter discard warning=false] table[x index=0, y index=5] {./figures/centerline_res800_ck6_t0p0030.dat};
    \addplot[name path=lo16, draw=none, forget plot, each nth point=4, filter discard warning=false] table[x index=0, y index=4] {./figures/centerline_res800_ck6_t0p0030.dat};
    \addplot[blue!45, forget plot] fill between[of=up84 and lo16];
    \addplot[orange, line width=0.3pt, forget plot, each nth point=4, filter discard warning=false] table[x index=0, y index=1] {./figures/centerline_res800_ck6_t0p0030.dat};
    \addplot[blue, line width=1.1pt, forget plot, each nth point=4, filter discard warning=false] table[x index=0, y index=2] {./figures/centerline_res800_ck6_t0p0030.dat};
    \addplot[red, densely dotted, line width=0.9pt, forget plot, each nth point=4, filter discard warning=false] table[x index=0, y index=3] {./figures/centerline_res800_ck6_t0p0030.dat};
    \end{axis}
    \end{tikzpicture}
    }%
    \subfloat[$t=0.003$, $N_x=4000$]{%
    \begin{tikzpicture}
    \begin{axis}[centerline style]
    \addplot[name path=up95, draw=none, forget plot, each nth point=4, filter discard warning=false] table[x index=0, y index=7] {./figures/centerline_res1600_ck6_t0p0030.dat};
    \addplot[name path=lo05, draw=none, forget plot, each nth point=4, filter discard warning=false] table[x index=0, y index=6] {./figures/centerline_res1600_ck6_t0p0030.dat};
    \addplot[blue!25, forget plot] fill between[of=up95 and lo05];
    \addplot[name path=up84, draw=none, forget plot, each nth point=4, filter discard warning=false] table[x index=0, y index=5] {./figures/centerline_res1600_ck6_t0p0030.dat};
    \addplot[name path=lo16, draw=none, forget plot, each nth point=4, filter discard warning=false] table[x index=0, y index=4] {./figures/centerline_res1600_ck6_t0p0030.dat};
    \addplot[blue!45, forget plot] fill between[of=up84 and lo16];
    \addplot[orange, line width=0.3pt, forget plot, each nth point=4, filter discard warning=false] table[x index=0, y index=1] {./figures/centerline_res1600_ck6_t0p0030.dat};
    \addplot[blue, line width=1.1pt, forget plot, each nth point=4, filter discard warning=false] table[x index=0, y index=2] {./figures/centerline_res1600_ck6_t0p0030.dat};
    \addplot[red, densely dotted, line width=0.9pt, forget plot, each nth point=4, filter discard warning=false] table[x index=0, y index=3] {./figures/centerline_res1600_ck6_t0p0030.dat};
    \end{axis}
    \end{tikzpicture}
    }
    \\[-0.3em]
    \subfloat[$t=0.0035$, $N_x=1000$]{%
    \begin{tikzpicture}
    \begin{axis}[centerline style]
    \addplot[name path=up95, draw=none, forget plot, each nth point=4, filter discard warning=false] table[x index=0, y index=7] {./figures/centerline_res400_ck7_t0p0035.dat};
    \addplot[name path=lo05, draw=none, forget plot, each nth point=4, filter discard warning=false] table[x index=0, y index=6] {./figures/centerline_res400_ck7_t0p0035.dat};
    \addplot[blue!25, forget plot] fill between[of=up95 and lo05];
    \addplot[name path=up84, draw=none, forget plot, each nth point=4, filter discard warning=false] table[x index=0, y index=5] {./figures/centerline_res400_ck7_t0p0035.dat};
    \addplot[name path=lo16, draw=none, forget plot, each nth point=4, filter discard warning=false] table[x index=0, y index=4] {./figures/centerline_res400_ck7_t0p0035.dat};
    \addplot[blue!45, forget plot] fill between[of=up84 and lo16];
    \addplot[orange, line width=0.3pt, forget plot, each nth point=4, filter discard warning=false] table[x index=0, y index=1] {./figures/centerline_res400_ck7_t0p0035.dat};
    \addplot[blue, line width=1.1pt, forget plot, each nth point=4, filter discard warning=false] table[x index=0, y index=2] {./figures/centerline_res400_ck7_t0p0035.dat};
    \addplot[red, densely dotted, line width=0.9pt, forget plot, each nth point=4, filter discard warning=false] table[x index=0, y index=3] {./figures/centerline_res400_ck7_t0p0035.dat};
    \end{axis}
    \end{tikzpicture}
    }%
    \subfloat[$t=0.0035$, $N_x=2000$]{%
    \begin{tikzpicture}
    \begin{axis}[centerline style]
    \addplot[name path=up95, draw=none, forget plot, each nth point=4, filter discard warning=false] table[x index=0, y index=7] {./figures/centerline_res800_ck7_t0p0035.dat};
    \addplot[name path=lo05, draw=none, forget plot, each nth point=4, filter discard warning=false] table[x index=0, y index=6] {./figures/centerline_res800_ck7_t0p0035.dat};
    \addplot[blue!25, forget plot] fill between[of=up95 and lo05];
    \addplot[name path=up84, draw=none, forget plot, each nth point=4, filter discard warning=false] table[x index=0, y index=5] {./figures/centerline_res800_ck7_t0p0035.dat};
    \addplot[name path=lo16, draw=none, forget plot, each nth point=4, filter discard warning=false] table[x index=0, y index=4] {./figures/centerline_res800_ck7_t0p0035.dat};
    \addplot[blue!45, forget plot] fill between[of=up84 and lo16];
    \addplot[orange, line width=0.3pt, forget plot, each nth point=4, filter discard warning=false] table[x index=0, y index=1] {./figures/centerline_res800_ck7_t0p0035.dat};
    \addplot[blue, line width=1.1pt, forget plot, each nth point=4, filter discard warning=false] table[x index=0, y index=2] {./figures/centerline_res800_ck7_t0p0035.dat};
    \addplot[red, densely dotted, line width=0.9pt, forget plot, each nth point=4, filter discard warning=false] table[x index=0, y index=3] {./figures/centerline_res800_ck7_t0p0035.dat};
    \end{axis}
    \end{tikzpicture}
    }%
    \subfloat[$t=0.0035$, $N_x=4000$]{%
    \begin{tikzpicture}
    \begin{axis}[centerline style]
    \addplot[name path=up95, draw=none, forget plot, each nth point=4, filter discard warning=false] table[x index=0, y index=7] {./figures/centerline_res1600_ck7_t0p0035.dat};
    \addplot[name path=lo05, draw=none, forget plot, each nth point=4, filter discard warning=false] table[x index=0, y index=6] {./figures/centerline_res1600_ck7_t0p0035.dat};
    \addplot[blue!25, forget plot] fill between[of=up95 and lo05];
    \addplot[name path=up84, draw=none, forget plot, each nth point=4, filter discard warning=false] table[x index=0, y index=5] {./figures/centerline_res1600_ck7_t0p0035.dat};
    \addplot[name path=lo16, draw=none, forget plot, each nth point=4, filter discard warning=false] table[x index=0, y index=4] {./figures/centerline_res1600_ck7_t0p0035.dat};
    \addplot[blue!45, forget plot] fill between[of=up84 and lo16];
    \addplot[orange, line width=0.3pt, forget plot, each nth point=4, filter discard warning=false] table[x index=0, y index=1] {./figures/centerline_res1600_ck7_t0p0035.dat};
    \addplot[blue, line width=1.1pt, forget plot, each nth point=4, filter discard warning=false] table[x index=0, y index=2] {./figures/centerline_res1600_ck7_t0p0035.dat};
    \addplot[red, densely dotted, line width=0.9pt, forget plot, each nth point=4, filter discard warning=false] table[x index=0, y index=3] {./figures/centerline_res1600_ck7_t0p0035.dat};
    \end{axis}
    \end{tikzpicture}
    }
    \\[-0.3em]
    \vspace{.5em}
    \ref*{centerlineLegend}
    \caption{Centerline density $\rho$ (at $y=0$) along the flow direction $x$
for the Mach~2000 jet, arranged as a matrix of resolutions (columns
$N_x=1000,2000,4000$) and physical times (rows
$t=0.001,0.002,0.003,0.0035$). Each panel shows a single representative
sample realization (orange), the ensemble mean $\mathbb{E}[\rho]$ (blue),
the ensemble standard deviation $\sqrt{\mathrm{Var}[\rho]}$ (red dotted),
and two nested quantile bands of the $M=1000$-sample ensemble, computed
from the empirical quantiles \eqref{eq:empcdf}: the inner,
darker band between the empirical $16\%$ and $84\%$ quantiles (containing
$68\%$ of the realizations pointwise) and the outer, lighter band between
the $5\%$ and $95\%$ quantiles (containing $90\%$). All panels share
identical axes so that the growth of the turbulent tip and the sharpening
of fine-scale structure with resolution and time are directly comparable.
For the purpose of illustration, only every fourth data point is shown.}
    \label{fig:centerline_matrix}
\end{figure}
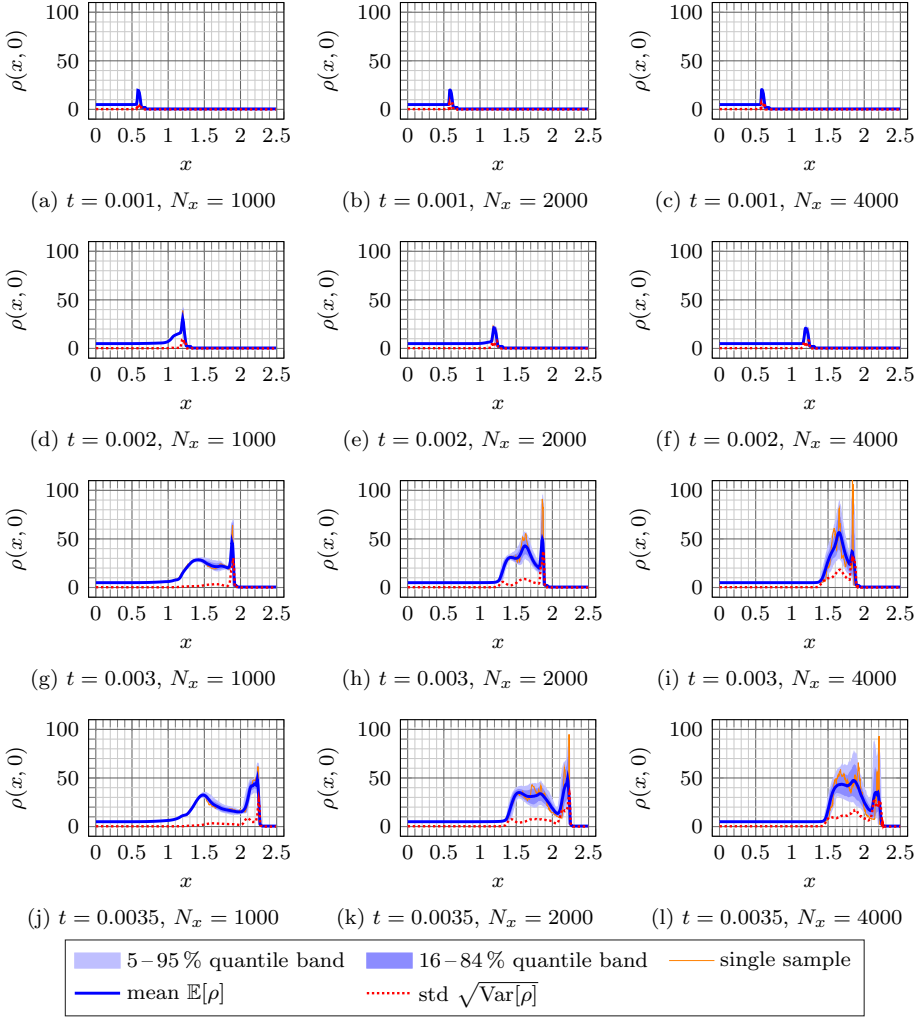

Figure~\ref{fig:centerline_matrix} complements the global convergence rates of Section~\ref{subsec:discussion} with a spatially resolved view of the statistical solution, showing the density along the jet centerline ($y=0$) as a function of the streamwise coordinate~$x$. Reading the matrix down its rows exposes the temporal development of the flow, while reading across its columns exposes the effect of mesh refinement at fixed time; the shared axes make both comparisons immediate. In every panel, we plot a single representative realization (orange), the ensemble mean $\mathbb{E}_M[\rho]$ (blue), the standard deviation $\sqrt{\mathrm{Var}_M[\rho]}$ (red dotted), and two nested quantile bands, the inner between the empirical $16\%$ and $84\%$ quantiles (dark shading) and the outer between the $5\%$ and $95\%$ quantiles (light shading), all computed from the full $M=1000$-sample ensemble of \eqref{eq:empirical}.

For definiteness, the bands are level sets of the empirical cumulative distribution function of the centerline density: evaluating the ensemble at the point $(x,0)$,
\begin{equation}\label{eq:empcdf}
    \hat F_{t,x}(\xi) \;=\; \frac{1}{M}\sum_{m=1}^{M} \mathbb{I}_{\{\rho^{\triangle x}(\omega_m,(x,0),t)\,\le\,\xi\}},
    \qquad
    p_{\alpha}(x) \;=\; \inf\bigl\{\xi\in\mathbb{R} : \hat F_{t,x}(\xi)\ge \alpha\bigr\},
\end{equation}
so that $p_{\alpha}$ is the empirical $\alpha$-quantile of the one-point law at $(x,0)$; these are precisely the quantile functions whose $L^1$ distance constitutes the one-point Wasserstein metric via \eqref{eq:w1cdf}--\eqref{eq:w1disc}. The bands are pointwise prediction bands for a single realization, containing $68\%$ and $90\%$ of the samples, respectively; for a Gaussian law the inner band would essentially coincide with mean\,$\pm$\,std. Unlike a mean\,$\pm$\,std envelope, however, the quantile bands are constrained to the physical range $\rho>0$ and presume neither symmetry nor unimodality, so their shape carries distributional information beyond the second moment. As a parameter-free diagnostic of that shape we use the quantile width ratio
\begin{equation}\label{eq:widthratio}
    \mathcal{Q}(x) \;=\; \frac{p_{95}(x)-p_{5}(x)}{p_{84}(x)-p_{16}(x)},
    \qquad
    \mathcal{Q}_{\mathcal{N}} \;=\; \frac{\Phi^{-1}(0.95)-\Phi^{-1}(0.05)}{\Phi^{-1}(0.84)-\Phi^{-1}(0.16)} \;\approx\; 1.65,
\end{equation}
where $\Phi$ denotes the standard normal cumulative distribution function; $\mathcal{Q}\equiv\mathcal{Q}_{\mathcal{N}}$ for every Gaussian law, independently of its mean and variance, so departures of $\mathcal{Q}$ from $\mathcal{Q}_{\mathcal{N}}$ measure anomalous tail weight without reference to any fitted model.

Three features organize the picture. First, upstream of the jet head, the flow is statistically quiescent: for $x\lesssim1.4$ at the latest time, the mean density sits at the injected core value $\bar{\rho}_{\mathrm{jet}}\approx5$ with a standard deviation of only $\approx0.3$, and the quantile bands collapse onto thin, symmetric sleeves around the mean. This residual spread is not noise but the injected inlet randomness itself: on the axis, the perturbation \eqref{eq:karhunenExp} carries a standard deviation of $A\bar{\rho}_{\mathrm{jet}}\bigl(\tfrac{1}{3}\sum_{k=1}^{K}k^{-4}\bigr)^{1/2}\approx0.30$, in agreement with the measured value, so the laminar core advects the inlet randomness downstream without amplifying it.

Second, the statistical spread above this advected inlet level is generated by the Kelvin--Helmholtz dynamics of the beam shear layers together with the mixing at the jet head, rather than at the head alone. The standard deviation rises sharply where the shear-layer instability attains nonlinear amplitude (near $x\approx1.43$ at $t=0.0035$), grows through the sheared cocoon region trailing the head, and peaks in the bow-shock and mixing region at the tip, where it reaches $\approx27$ against a local mean of order $20$--$45$; near the very head, the relative fluctuation $\sqrt{\mathrm{Var}_M[\rho]}/\mathbb{E}_M[\rho]$ approaches $0.7$, signaling that the density there is genuinely random from realization to realization rather than a sharp deterministic front. The quantile bands make the character of this randomness visible and quantifiable. They are asymmetric about the mean. At the finest resolution and latest time, the median of $(p_{95}-\mathbb{E}_M[\rho])/(\mathbb{E}_M[\rho]-p_{5})$ along the head is $\approx1.4$, which is the signature of a positively skewed one-point law generated by intermittent high-density shocklets. The nested bands further separate two regimes of randomness. In the interior of the mixing region, the width ratio \eqref{eq:widthratio} stays at its Gaussian value (median $\mathcal{Q}\approx1.67$ along the head), so the one-point law there, although skewed, carries no anomalous tail weight. At the leading fringe of the head ($x\approx2.27$ at $t=0.0035$), by contrast, the ratio grows to order $10^{3}$: the inner band collapses onto the undisturbed ambient state, i.e., the central $68\%$ of realizations are still unperturbed, while the $95\%$ quantile has already detached, meaning the jet has reached these locations in more than $5\%$ but fewer than $16\%$ of the realizations. The one-point law there is an intermittent two-state mixture governed by the random position of the front. This is the classical external intermittency of turbulent/non-turbulent interfaces \cite{corrsin1955, dasilva2014}, arising here across the ensemble at fixed time rather than along a time series at a fixed probe, and it is the same bimodal structure resolved directly by the density histograms of Figure~\ref{fig:pdf_heatmap} (cf.~\eqref{eq:mixture} below). The single-sample trace makes the skewness concrete: the plotted realization spikes to $\rho\approx94$ at $x\approx2.20$, where the ensemble mean is only $\approx27$ and the $95\%$ quantile $\approx81$; such excursions beyond the outer band occur, by construction, at only a few percent of locations for any fixed realization, yet somewhere along the head in essentially every realization.

Third, the matrix demonstrates directly why the statistical solution converges while the pathwise field does not. Moving down each column, the turbulent tip advances downstream and the quantile bands broaden and move with it, tracing the growth of the statistically uncertain region as the jet penetrates the ambient medium. Moving across each row, from $N_x=1000$ to $4000$, the ensemble mean, the standard-deviation profile, and both quantile bands are visually stable, i.e., their large-scale shape, amplitude, and the location of the tip are change marginally, even as the single-sample trace becomes progressively rougher and resolves ever finer shocklets and shear-layer filaments. In other words, mesh refinement adds structure to individual realizations without shifting the low-order statistics: the coarse-grained mean, variance, and quantile fields are grid-converged in the sense measured by the reference protocol \eqref{eq:ref}, while the pathwise field is not, exactly as expected for a non-Dirac statistical solution. The near-collapse of the mean and the quantile bands across the three columns is the visual counterpart of the stable reference rates of $\mathcal{W}_{1,1}$ and the mean reported in Section~\ref{subsec:discussion}, and the simultaneous roughening of the orange trace is the counterpart of the strong error's~\eqref{eq:strongError} failure to converge.

\subsection{One-point distribution of the centerline density}\label{subsec:pdfheatmap}

\begin{figure}[ht!]
    \pgfplotsset{
        pdfaxis style/.style={
                scale only axis,
                width=0.75\textwidth, 
                height=0.09\textwidth,
       	        label style={font=\footnotesize},
	            tick label style={font=\footnotesize},
                xlabel={$x$}, 
                ylabel={$\rho(x,0)$},
                xmin=0, 
                xmax=2.5,
                ymin=0, 
                ymax=150,
                enlargelimits=false, 
                axis on top,
                colormap/cool,
                point meta min=-3.410,
                point meta max=-0.410
        }
    }
    \centering
    \subfloat[\(t=0.0015\)]{%
        \begin{tikzpicture}
            \begin{axis}[pdfaxis style]
                \addplot graphics[xmin=0, xmax=2.5,
                    ymin=0, ymax=150] {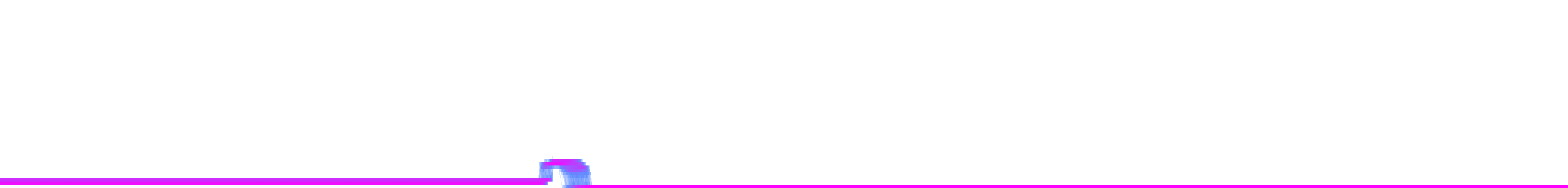};
            \end{axis}
        \end{tikzpicture}}
    \\
    \subfloat[\(t=0.0020\)]{%
        \begin{tikzpicture}
            \begin{axis}[pdfaxis style]
                \addplot graphics[xmin=0, xmax=2.5,
                    ymin=0, ymax=150] {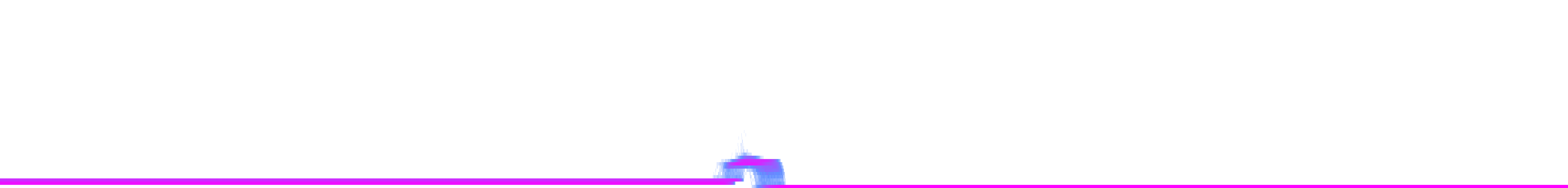};
            \end{axis}
        \end{tikzpicture}}
    \\
    \subfloat[\(t=0.0025\)]{%
        \begin{tikzpicture}
            \begin{axis}[pdfaxis style]
                \addplot graphics[xmin=0, xmax=2.5,
                    ymin=0, ymax=150] {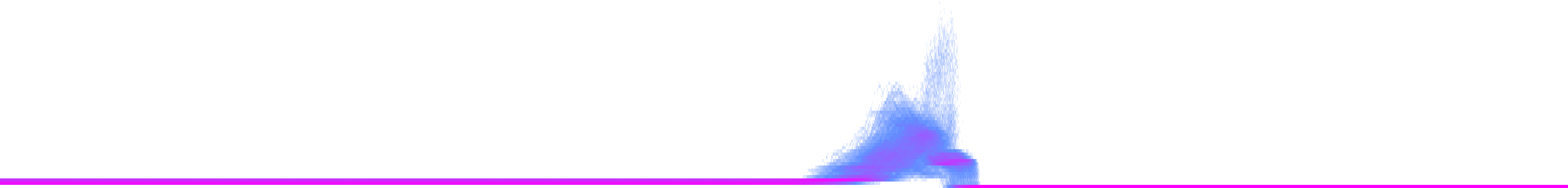};
            \end{axis}
        \end{tikzpicture}}
    \\
    \subfloat[\(t=0.0030\)]{%
        \begin{tikzpicture}
            \begin{axis}[pdfaxis style]
                \addplot graphics[xmin=0, xmax=2.5,
                    ymin=0, ymax=150] {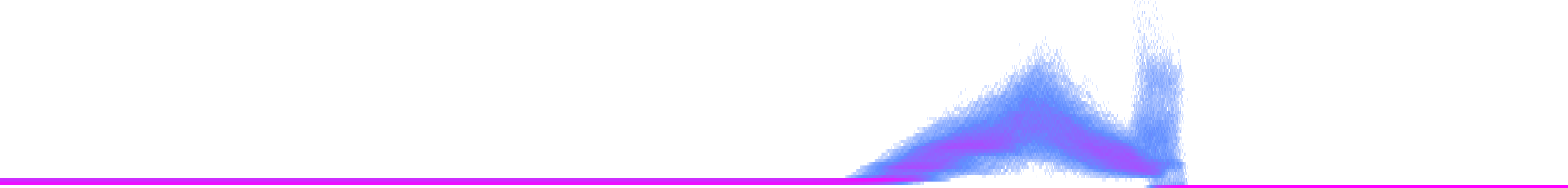};
            \end{axis}
        \end{tikzpicture}}
    \\
    \subfloat[\(t=0.0035\)]{%
        \begin{tikzpicture}
            \begin{axis}[pdfaxis style]
                \addplot graphics[xmin=0, xmax=2.5,
                    ymin=0, ymax=150] {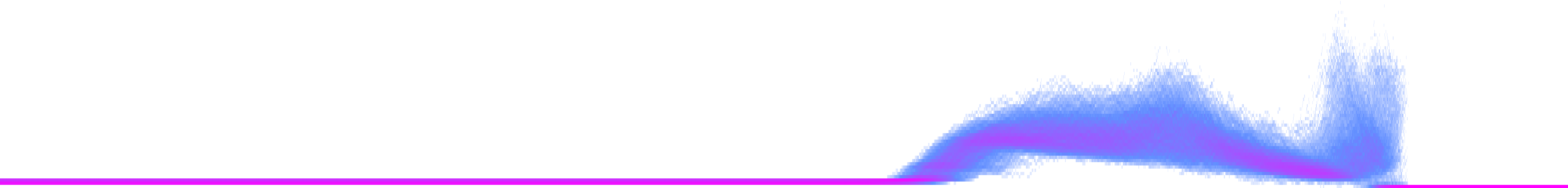};
            \end{axis}
        \end{tikzpicture}}
    \\ \vspace{.2em}
    \centerline{%
        \begin{tikzpicture}
            \pgfplotscolorbardrawstandalone[
                colormap/cool,
                point meta min=-3.410,
                point meta max=-0.410,
                colorbar horizontal,
                colorbar style={
                    ylabel={$\log_{10} f(x,\rho(x,0))$},
                    height={0.2cm},
                    width={6cm},
                    ylabel style={rotate=-90, font=\footnotesize},
                    xtick={-3.41,-3,-2,-1,-0.41},
                    xticklabel style={font=\footnotesize}
                }
            ]
        \end{tikzpicture}}
    \caption{Per-slice normalized probability density \eqref{eq:hist} of
    the centerline density $\rho(x,0)$ across the $M=1000$ samples at
    several timesteps for the \(N_{x} = 4000\) grid. The colormap spans
    the global minimum and maximum of \(\log_{10} f\) over all positive
    values and all displayed times; bins carrying no samples are shown
    in white.}
    \label{fig:pdf_heatmap}
\end{figure}

Figure~\ref{fig:pdf_heatmap} sharpens the moment- and quantile-level description of
Figure~\ref{fig:centerline_matrix} to the level of the full one-point
distribution. For a fixed time $t$ and streamwise location $x$, the
ensemble of solutions induces a probability measure on the local density
values, namely the law of the random variable $\rho(x,0,t)$,
\begin{equation}\label{eq:onepointlaw}
    \nu_{t,x} \;=\; \operatorname{law}\bigl(\rho(x,0,t)\bigr)
    \;\in\; \mathcal{P}(\mathbb{R}_{+}),
\end{equation}
which is precisely the first one-point correlation marginal of the
statistical solution \eqref{eq:mut_exact} restricted to the density
component and evaluated on the jet axis. We estimate $\nu_{t,x}$
from the $M=1000$ samples by the one-point marginal of the empirical
measure \eqref{eq:empirical},
$\smash{\hat\nu_{t,x} = M^{-1}\sum_{m=1}^{M}\delta_{\rho^{\triangle x}(\omega_m,(x,0),t)}}$,
and display its histogram density on a fixed partition $\{B_j\}$ of the
density axis,
\begin{equation}\label{eq:hist}
    f(x,\rho) \;=\; \frac{1}{M\,|B_j|}\,
      \#\bigl\{\, m : \rho^{\triangle x}(\omega_m,(x,0),t)\in B_j \,\bigr\}
    \qquad \text{for } \rho\in B_j,
\end{equation}
normalized slice by slice so that $\int f(x,\rho)\,\mathrm{d}\rho = 1$
for every $x$. Each vertical slice of a panel is therefore a probability
density in $\rho$, and the panels as a whole resolve the object whose
first two moments, $\mathbb{E}[\rho](x) = \int \rho
\,\nu_{t,x}(\mathrm{d}\rho)$ and $\mathrm{Var}[\rho](x) = \int (\rho -
\mathbb{E}[\rho](x))^2 \,\nu_{t,x}(\mathrm{d}\rho)$, and whose quantile
functions \eqref{eq:empcdf} were plotted in
Figure~\ref{fig:centerline_matrix}. Because $f$ varies over several
orders of magnitude, the color encodes $\log_{10} f$; the scale is
common to all panels and spans the global minimum and maximum of
$\log_{10} f$ over all positive values and all displayed times, so that
color is directly comparable across timesteps. Bins carrying no samples
are shown in white.

The heatmaps make visible what the moments necessarily compress. Away
from the jet head, each marginal concentrates in essentially a single
bin around the core value $\bar\rho_{\mathrm{jet}}\approx 5$: the thin
magenta filament along the bottom of every panel is a numerically
Dirac one-point law, in agreement with the vanishing standard deviation
observed there and with the interpretation of the statistical solution
as deterministic in the quiescent region. At the turbulent head, by
contrast, the marginals are emphatically non-atomic. Reading the panels
in time, the support of $\nu_{t,x}$ fans out from a compact spot at
$t=0.0015$ into a broad, structured cloud that by $t=0.0035$ extends
over more than two decades of density, up to $\rho \approx 140$, while
the bulk of the probability, traced by the magenta ridge, remains near
$\rho \approx 40$. The marginals there are strongly skewed and, in the
bow-shock region, bimodal: a random shock position $X_s$ turns the
pointwise density into a two-state mixture,
\begin{equation}\label{eq:mixture}
    \nu_{t,x} \;\approx\; p(x)\,\nu^{\mathrm{pre}}
      + \bigl(1-p(x)\bigr)\,\nu^{\mathrm{post}},
    \qquad p(x)=\mathbb{P}(X_s > x),
\end{equation}
broadened by turbulent fluctuations of the pre- and post-shock states
themselves; the divergence of the quantile width ratio
\eqref{eq:widthratio} at the leading fringe observed in
Section~\ref{subsec:centerline} is the quantile-level fingerprint of
exactly this mixture. The structure \eqref{eq:mixture} also explains
quantitatively why the moment-level picture of
Figure~\ref{fig:centerline_matrix} looks the way it does: the ensemble
mean interpolates between the mixture components and thus runs through
a region of low probability density, single realizations overshoot it
by a factor of $\approx3.5$ (the spike to $\rho\approx 94$ against a
local mean of $\approx 27$), and any summary of the form
mean\,$\pm$\,std, which implicitly presumes a unimodal, roughly
symmetric law, would misrepresent a distribution that is neither. This is the
motivation for the nested quantile bands employed in
Figure~\ref{fig:centerline_matrix}.

The figure also identifies the object on which the convergence theory
of Section~\ref{sec:weak-strong} operates. The one-point Wasserstein
error \eqref{eq:w1cont}--\eqref{eq:w1disc} is, by the identity
\eqref{eq:w1cdf}, the spatially averaged $L^1$ distance between the
cumulative distribution functions of exactly the marginals displayed
here. Convergence in this metric requires the histograms of
Figure~\ref{fig:pdf_heatmap} , i.e., their support, modal ridge, and tail
weight, to stabilize under mesh refinement, but places no constraint on
where within that distribution an individual realization falls. This is
the measure-theoretic content of the dichotomy observed globally: the
strong, pathwise error \eqref{eq:strongError} stalls because refined
samples decorrelate and merely resample the (broad) law at the
turbulent head, while the statistical solution converges because the
law itself is stable. The near-Dirac marginals upstream and the broad
marginals at the head thus delineate, within a single figure, the
regions where the computation is effectively deterministic and those
where only statistical statements are meaningful.

\subsection{Divergence of the strong error versus statistical convergence}\label{subsec:discussion}

Throughout this discussion, we read the \textit{reference} protocol \eqref{eq:ref} as the estimate of convergence toward the (finest-grid surrogate of the) statistical solution, and the \textit{Cauchy} protocol \eqref{eq:cauchy} as the estimate of internal, self-referential convergence between successive resolutions; the systematic offset between the two is expected and is interpreted below.

\begin{figure}[ht!]
	\pgfplotsset{
        common rate style/.style={
            width=0.28\textwidth,
	        height=0.33\textwidth,
	        xlabel={$t \,\text{[s]}$},
	        ylabel={$\mathfrak{r}_{\mathcal{W}_{1,1}}(t)$},
	        xmin=-0.0001,
	        xmax=0.0041,
	        ymin=-1.2,
	        ymax=1.7,
	        label style={font=\footnotesize},
            xlabel style={at={(0.5,-0.4)}},
            tick label style={font=\footnotesize},
	        yminorgrids=true,
	        xminorgrids=true,
	        grid=both,
	        minor x tick num=4,
	        minor x grid style={gray!40},
	        major x grid style={black!60},
	        minor y tick num=4,
	        minor y grid style={gray!40},
	        major y grid style={black!60}
        }
    }
	\centering
	\subfloat[\(\mathcal{W}_{1,1}\)-rate]{
		\begin{tikzpicture}
		    \begin{axis}[
		    	common rate style,
		        ylabel={$\mathfrak{r}_{\mathcal{W}_{1,1}}(t)$},
		        legend to name=ratesLegend,
		        legend style={
                    legend columns=2,
                    /tikz/every even column/.append style={column sep=0.25cm},
                    draw=darkgray!60!black,
                    fill=white,
                    legend cell align=left,
                    font = \footnotesize}
		    ]
		    \addplot[color=teal, solid, mark=square, line width=1.5pt] coordinates {
		        (0.0005,0.374) (0.0010,0.339) (0.0015,0.394) (0.0020,0.313) (0.0025,0.337) (0.0030,0.336) (0.0035,0.342)
		    };
		    \addlegendentry{\textit{Cauchy} protocol}
		    \addplot[color=orange, solid, mark=o, line width=1.5pt] coordinates {
		        (0.0005,0.774) (0.0010,0.739) (0.0015,0.758) (0.0020,0.686) (0.0025,0.637) (0.0030,0.655) (0.0035,0.674)
		    };
		    \addlegendentry{\textit{reference} protocol}
		    \end{axis}
		\end{tikzpicture}
		\label{fig:w11-over-time}
		}%
		\subfloat[\(\mathcal{W}_{1,2}\)-rate]{
		\begin{tikzpicture}
		    \begin{axis}[common rate style, ylabel={$\mathfrak{r}_{\mathcal{W}_{1,2}}(t)$}]
		    \addplot[forget plot, color=teal, solid, mark=square, line width=1.5pt] coordinates {
		        (0.0005,-1.091) (0.0010,0.211) (0.0015,0.163) (0.0020,0.299) (0.0025,0.417) (0.0030,1.052) (0.0035,1.200)
		    };
		    \addplot[forget plot, color=orange, solid, mark=o, line width=1.5pt] coordinates {
		        (0.0005,0.331) (0.0010,0.637) (0.0015,0.624) (0.0020,0.888) (0.0025,0.817) (0.0030,0.872) (0.0035,1.029)
		    };
		    \end{axis}
		\end{tikzpicture}
		\label{fig:w12-over-time}
		}%
		\subfloat[$\mathcal{E}_{\mathrm{mean}}$ rate]{
		\begin{tikzpicture}
		    \begin{axis}[common rate style, ylabel={$\mathfrak{r}_{\mathcal{E}_{\mathrm{mean}}}(t)$}]
		    \addplot[forget plot, color=teal, solid, mark=square, line width=1.5pt] coordinates {
		        (0.0005,0.377) (0.0010,0.374) (0.0015,0.420) (0.0020,0.337) (0.0025,0.353) (0.0030,0.345) (0.0035,0.351)
		    };
		    \addplot[forget plot, color=orange, solid, mark=o, line width=1.5pt] coordinates {
		        (0.0005,0.772) (0.0010,0.750) (0.0015,0.764) (0.0020,0.692) (0.0025,0.634) (0.0030,0.642) (0.0035,0.659)
		    };
		    \end{axis}
		\end{tikzpicture}
		\label{fig:mean-over-time}
		}%
		\\
		\subfloat[$\mathcal{E}_{\mathrm{std}}$-rate]{
		\begin{tikzpicture}
		    \begin{axis}[common rate style, ylabel={$\mathfrak{r}_{\mathcal{E}_{\mathrm{std}}}(t)$}]
		    \addplot[forget plot, color=teal, solid, mark=square, line width=1.5pt] coordinates {
		        (0.0005,0.408) (0.0010,-0.109) (0.0015,0.021) (0.0020,0.228) (0.0025,0.221) (0.0030,0.175) (0.0035,0.053)
		    };
		    \addplot[forget plot, color=orange, solid, mark=o, line width=1.5pt] coordinates {
		        (0.0005,0.779) (0.0010,0.461) (0.0015,0.639) (0.0020,0.775) (0.0025,0.772) (0.0030,0.744) (0.0035,0.609)
		    };
		    \end{axis}
		\end{tikzpicture}
		\label{fig:std-over-time}
		}%
		\subfloat[$\mathcal{E}_{\mathrm{strong}}$-rate]{
				\begin{tikzpicture}
				    \begin{axis}[common rate style, ylabel={$\mathfrak{r}_{\mathcal{E}_{\mathrm{strong}}}(t)$}]
				    \addplot[forget plot, color=teal, solid, mark=square, line width=1.5pt] coordinates {
				        (0.0005,0.366) (0.0010,0.253) (0.0015,0.218) (0.0020,0.110) (0.0025,0.095) (0.0030,0.036) (0.0035,0.026)
				    };
				    \addplot[forget plot, color=orange, solid, mark=o, line width=1.5pt] coordinates {
				        (0.0005,0.760) (0.0010,0.653) (0.0015,0.586) (0.0020,0.476) (0.0025,0.393) (0.0030,0.344) (0.0035,0.344)
				    };
				    \end{axis}
				\end{tikzpicture}
			\label{fig:strong-over-time}	
		}
		\\
		\vspace{.5em}
		\ref*{ratesLegend}
		\caption{Fitted convergence rate $\mathfrak{r}_{\mathfrak{O}}(t)$ of \eqref{eq:ratefit} as a function of time for each observable, one panel per observable, with both comparison operators (\textit{Cauchy}~\eqref{eq:cauchy}, teal squares; \textit{reference}~\eqref{eq:ref}, orange circles) shown at fixed sample size $M=1000$. Panels: (a) one-point Wasserstein $\mathcal{W}_{1,1}$, (b) two-point Wasserstein $\mathcal{W}_{1,2}$, (c) mean, (d) standard deviation, (e) strong pathwise error. All panels share the axes and grid of panel~(a).}
	\label{fig:protocol_rates}
\end{figure}
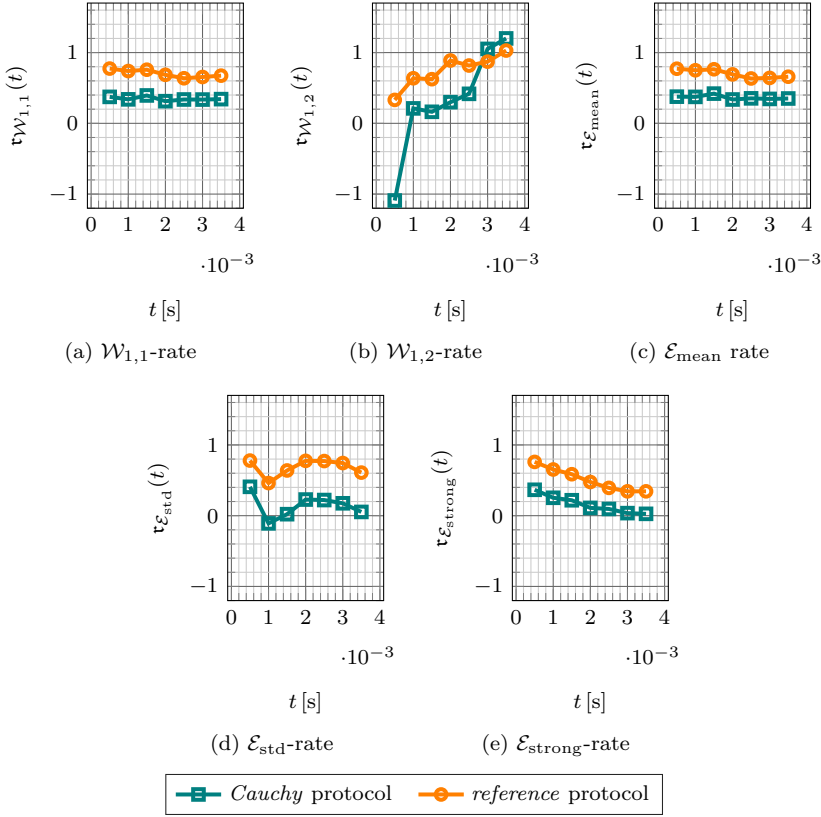

Figure~\ref{fig:protocol_rates} presents the fitted rate $\mathfrak{r}_{\mathfrak{O}}(t)$ of \eqref{eq:ratefit} as a function of physical time for each observable, the five panels sharing common axes to make the qualitative contrasts directly comparable. The panels are ordered as the one-point Wasserstein distance $\mathcal{W}_{1,1}$~(Figure~\ref{fig:w11-over-time}), the two-point Wasserstein distance $\mathcal{W}_{1,2}$~(Figure~\ref{fig:w12-over-time}), the mean~(Figure~\ref{fig:mean-over-time}), the standard deviation~(Figure~\ref{fig:std-over-time}), and the strong pathwise error~(Figure~\ref{fig:strong-over-time}); teal squares denote the \textit{Cauchy} protocol and orange circles the \textit{reference} protocol throughout. Read together, the figure exhibits the qualitative signature that organizes the remainder of this section: the statistical observables in Figures~\ref{fig:w11-over-time}--\ref{fig:std-over-time} settle onto stable, positive \textit{reference} rates as the flow becomes turbulent, whereas the strong error in Figure~\ref{fig:strong-over-time} decays monotonically toward zero under the \textit{Cauchy} comparison. The two families of curves move in opposite directions, statistical rates holding up while the pathwise rate collapses; the tables below quantify this behavior, and Sections~\ref{subsec:centerline} and~\ref{subsec:pdfheatmap} give its spatially resolved counterpart.

\paragraph{Early acoustic transient.} At very early times, the flow is dominated by the initial Riemann problem at the inlet and the launch of the primary bow shock into the quiescent ambient medium. In this phase the perturbation-induced ensemble spread is still small, and several observables are correspondingly ill-conditioned. This is most visible in the two-point Wasserstein distance $\mathcal{W}_{1,2}$ of \eqref{eq:w2}, whose \textit{Cauchy} rate at $t=0.0005$ is a spurious $-1.091$: with a nearly deterministic, still-laminar core the two-point correlations carry almost no statistical signal, and the fitted slope is meaningless rather than physical. The first- and second-moment observables \eqref{eq:mean}--\eqref{eq:std} and the one-point distance $\mathcal{W}_{1,1}$ of \eqref{eq:w1disc} are already well-behaved at this stage, with \textit{reference} rates clustered near $0.77$--$0.78$ (Table~\ref{tab:err-t0p0005}), reflecting the smoothness of the pre-turbulent ensemble. We therefore regard $t\lesssim0.0015$ as a transient window in which the two-point statistics, in particular, should not be over-interpreted.

\paragraph{Onset of shear-layer turbulence and the pathwise breakdown.} As the jet penetrates the ambient medium, the shear layers between the supersonic core and the surrounding gas become Kelvin--Helmholtz unstable, and by $t\gtrsim0.0015$ the head of the jet rolls up into fully developed, chaotic turbulence. The signature of this transition in our metrics is the monotone collapse of the \emph{strong} sample-wise error \eqref{eq:strongError}. Under the \textit{Cauchy} comparison, the strong rate $\mathfrak{r}_{\mathcal{E}_{\mathrm{strong}}}$ decays from $0.366$ at $t=0.0005$ through $0.110$ at $t=0.0020$ to essentially zero at late times ($0.036$ at $t=0.0030$ and $0.026$ at $t=0.0035$; Table~\ref{tab:err-t0p0035}). A convergence rate indistinguishable from zero means that refining the mesh no longer brings an individual realization closer to its finer-grid counterpart: the pointwise trajectory of the flow does not converge. This is precisely the pathwise non-uniqueness anticipated for the multi-dimensional compressible Euler equations in the turbulent regime, in which arbitrarily fine vortical structures continue to appear under refinement and no strong limit exists. We emphasize that the \textit{Cauchy} strong rate is the appropriate diagnostic for this statement, as it is free of any reference-grid bias; the \textit{reference} strong rate remains positive ($0.344$ at $t=0.0035$) only because it measures the approach of each coarse field to the fixed finest grid, which is a coarse-graining statement and not evidence of a continuum pathwise limit.

\paragraph{Convergence of the statistical solution.} The statistical observables behave in the opposite way and converge throughout the turbulent phase. The one-point Wasserstein distance $\mathcal{W}_{1,1}$ maintains a stable \textit{reference} rate close to $0.66$--$0.69$ across the entire late window, reaching $0.674$ at $t=0.0035$, and the ensemble mean tracks it almost exactly ($0.659$ at the final time). The near-coincidence of the $\mathcal{W}_{1,1}$ and mean rates is expected since both are low-order functionals of the same one-point marginal of the statistical solution (cf.~\eqref{eq:w1cdf} and \eqref{eq:mean}), and their agreement is a useful internal consistency check. The corresponding \textit{Cauchy} rates for these two observables sit near $0.34$--$0.35$; the roughly factor-of-two offset between the \textit{reference} and \textit{Cauchy} rates is the systematic lever-arm effect of comparing against the distant finest grid versus comparing adjacent resolutions, and it is stable precisely because $\mathcal{W}_{1,1}$ and the mean are clean, well-resolved quantities. This juxtaposition, a pathwise error that refuses to converge next to a probability measure that converges at a definite positive rate, is the central finding of this work and the signature of a non-Dirac statistical solution: the flow does not settle on a unique deterministic state, yet its law is well-defined and mesh-convergent.

\paragraph{Second moment and the two-point distance.} The standard deviation and the two-point Wasserstein distance are more delicate, and their behavior is informative rather than anomalous. The \textit{reference} standard-deviation rate is strong and stable in the turbulent window ($0.775$, $0.772$, $0.744$ at $t=0.0020,0.0025,0.0030$, easing to $0.609$ at $t=0.0035$), confirming that the second moment of the one-point law converges. Its \textit{Cauchy} counterpart, however, is erratic and small (ranging over $-0.109$ to $0.228$), and at early times even changes sign. This is the expected consequence of the second moment being the most sampling- and resolution-sensitive of the one-point functionals: the consecutive-grid difference of two standard-deviation fields is a small quantity formed from the tails of the marginal, where both the Monte Carlo sampling error at $M=1000$ and the residual grid-scale variance are largest, so the \textit{Cauchy} standard-deviation comparison has the poorest signal-to-noise of all observables and should be read from the \textit{reference} protocol only. The two-point distance $\mathcal{W}_{1,2}$ shows the opposite time trend to the strong error: having been meaningless during the transient, it becomes the fastest-converging observable at late times, with \textit{reference} rates rising to $0.888$ at $t=0.0020$ and $1.029$ at $t=0.0035$, and \textit{Cauchy} rates climbing even higher, to $1.052$ and $1.200$. That the \textit{Cauchy} two-point rate overtakes the \textit{reference} rate late in the simulation, which is the reverse of the ordering seen for $\mathcal{W}_{1,1}$ and the mean, indicates that the dominant two-point correlations are large-scale features (the coherent roll-up and bow-shock structure) on which consecutive grids agree increasingly well as they refine, whereas the finest grid still carries additional fine-scale correlation content that penalizes the reference comparison. We note that the mid-time $\mathcal{W}_{1,2}$ fits are based on a single ensemble draw and exhibit the largest scatter of any observable, so only the late-time two-point rates should be assigned quantitative weight.

\paragraph{Synthesis.} The five observables are consistent with the theory of statistical solutions for the compressible Euler equations. The vanishing \textit{Cauchy} strong rate certifies that individual realizations do not converge; the stable positive \textit{reference} rates certify that the underlying probability measure does. The spatially resolved analysis of Sections~\ref{subsec:centerline} and~\ref{subsec:pdfheatmap} identifies the mechanism: the one-point marginals are near-Dirac upstream and broad, skewed, intermittent mixtures \eqref{eq:mixture} in the sheared turbulent region and at the head, so refined samples decorrelate and merely resample a law that is itself stable under refinement. The persistence of a non-degenerate, mesh-convergent law in the face of pathwise divergence is exactly the behavior predicted when strong solutions break down and are replaced by measure-valued statistical solutions~\cite{fjordholm,fjordholm2017,lye2020}. We do not claim convergence to a specific literature rate: the measured one-point rates of $\approx0.66$--$0.69$ (reference) and $\approx0.34$ (Cauchy) are properties of this Mach~$2000$ configuration and the VLBM discretization, and they should be read as evidence for the \emph{qualitative} dichotomy of pathwise divergence with statistical convergence, rather than as an estimate of a universal exponent.

\begin{table}[ht!]
\centering
\caption{Per-resolution $L^1(\mathcal{D})$ error values and fitted convergence rates $\mathfrak{r}$ of \eqref{eq:ratefit} at $t=0.0005$, with observables in columns and protocols in row blocks. All fields are conservatively block-averaged to the common coarsest grid of $500\times400$ cells before evaluation, and $M=1000$ samples are used at every resolution. Within each protocol block the first three rows are the per-pair errors (grid pair given as coarse/fine streamwise cell count $N_x$) and the bold $\mathfrak{r}$ row is the least-squares slope of $\log(\text{error})$ against $\log(\triangle x)$.}
\label{tab:err-t0p0005}
\footnotesize
\setlength{\tabcolsep}{2pt}
\begin{tabular}{llr rrrrr}
\toprule
$t$ & Protocol & Grid pair & $\mathcal{W}_{1,1}$ & $\mathcal{W}_{1,2}$ & $\mathcal{E}_{\mathrm{mean}}$ & $\mathcal{E}_{\mathrm{std}}$ & $\mathcal{E}_{\mathrm{strong}}$ \\
\midrule
\multirow{8}{*}{$0.0005$} & \multirow{4}{*}{\textit{Cauchy}} & 500/1000 & $5.813 \times 10^{-3}$ & $9.502 \times 10^{-4}$ & $5.623 \times 10^{-3}$ & $1.108 \times 10^{-3}$ & $1.098 \times 10^{-2}$ \\
 &  & 1000/2000 & $4.157 \times 10^{-3}$ & $3.416 \times 10^{-3}$ & $4.011 \times 10^{-3}$ & $6.170 \times 10^{-4}$ & $7.829 \times 10^{-3}$ \\
 &  & 2000/4000 & $3.461 \times 10^{-3}$ & $4.309 \times 10^{-3}$ & $3.333 \times 10^{-3}$ & $6.295 \times 10^{-4}$ & $6.607 \times 10^{-3}$ \\
 & & $\mathfrak{r}$ (rate) & \textbf{0.374} & \textbf{-1.091} & \textbf{0.377} & \textbf{0.408} & \textbf{0.366} \\
\cmidrule(l){2-8}
 & \multirow{4}{*}{\textit{reference}} & 500/4000 & $1.012 \times 10^{-2}$ & $6.817 \times 10^{-3}$ & $9.723 \times 10^{-3}$ & $1.853 \times 10^{-3}$ & $1.895 \times 10^{-2}$ \\
 &  & 1000/4000 & $6.520 \times 10^{-3}$ & $7.146 \times 10^{-3}$ & $6.252 \times 10^{-3}$ & $1.054 \times 10^{-3}$ & $1.225 \times 10^{-2}$ \\
 &  & 2000/4000 & $3.461 \times 10^{-3}$ & $4.309 \times 10^{-3}$ & $3.333 \times 10^{-3}$ & $6.295 \times 10^{-4}$ & $6.607 \times 10^{-3}$ \\
 & & $\mathfrak{r}$ (rate) & \textbf{0.774} & \textbf{0.331} & \textbf{0.772} & \textbf{0.779} & \textbf{0.760} \\
\bottomrule
\end{tabular}
\end{table}

\begin{table}[ht!]
\centering
\caption{Per-resolution $L^1(\mathcal{D})$ error values and fitted convergence rates $\mathfrak{r}$ of \eqref{eq:ratefit} at $t=0.0010$, with observables in columns and protocols in row blocks. All fields are conservatively block-averaged to the common coarsest grid of $500\times400$ cells before evaluation, and $M=1000$ samples are used at every resolution. Within each protocol block the first three rows are the per-pair errors (grid pair given as coarse/fine streamwise cell count $N_x$) and the bold $\mathfrak{r}$ row is the least-squares slope of $\log(\text{error})$ against $\log(\triangle x)$.}
\label{tab:err-t0p001}
\footnotesize
\setlength{\tabcolsep}{2pt}
\begin{tabular}{llr rrrrr}
\toprule
$t$ & Protocol & Grid pair & $\mathcal{W}_{1,1}$ & $\mathcal{W}_{1,2}$ & $\mathcal{E}_{\mathrm{mean}}$ & $\mathcal{E}_{\mathrm{std}}$ & $\mathcal{E}_{\mathrm{strong}}$ \\
\midrule
\multirow{8}{*}{$0.0010$} & \multirow{4}{*}{\textit{Cauchy}} & 500/1000 & $1.647 \times 10^{-2}$ & $2.588 \times 10^{-2}$ & $1.593 \times 10^{-2}$ & $3.049 \times 10^{-3}$ & $2.932 \times 10^{-2}$ \\
 &  & 1000/2000 & $1.480 \times 10^{-2}$ & $3.498 \times 10^{-2}$ & $1.422 \times 10^{-2}$ & $2.839 \times 10^{-3}$ & $2.689 \times 10^{-2}$ \\
 &  & 2000/4000 & $1.030 \times 10^{-2}$ & $1.931 \times 10^{-2}$ & $9.477 \times 10^{-3}$ & $3.547 \times 10^{-3}$ & $2.065 \times 10^{-2}$ \\
 & & $\mathfrak{r}$ (rate) & \textbf{0.339} & \textbf{0.211} & \textbf{0.374} & \textbf{-0.109} & \textbf{0.253} \\
\cmidrule(l){2-8}
 & \multirow{4}{*}{\textit{reference}} & 500/4000 & $2.868 \times 10^{-2}$ & $4.674 \times 10^{-2}$ & $2.680 \times 10^{-2}$ & $6.719 \times 10^{-3}$ & $5.104 \times 10^{-2}$ \\
 &  & 1000/4000 & $2.098 \times 10^{-2}$ & $3.205 \times 10^{-2}$ & $1.948 \times 10^{-2}$ & $5.435 \times 10^{-3}$ & $3.829 \times 10^{-2}$ \\
 &  & 2000/4000 & $1.030 \times 10^{-2}$ & $1.931 \times 10^{-2}$ & $9.477 \times 10^{-3}$ & $3.547 \times 10^{-3}$ & $2.065 \times 10^{-2}$ \\
 & & $\mathfrak{r}$ (rate) & \textbf{0.739} & \textbf{0.637} & \textbf{0.750} & \textbf{0.461} & \textbf{0.653} \\
\bottomrule
\end{tabular}
\end{table}

\begin{table}[ht!]
\centering
\caption{Per-resolution $L^1(\mathcal{D})$ error values and fitted convergence rates $\mathfrak{r}$ of \eqref{eq:ratefit} at $t=0.0015$, with observables in columns and protocols in row blocks. All fields are conservatively block-averaged to the common coarsest grid of $500\times400$ cells before evaluation, and $M=1000$ samples are used at every resolution. Within each protocol block the first three rows are the per-pair errors (grid pair given as coarse/fine streamwise cell count $N_x$) and the bold $\mathfrak{r}$ row is the least-squares slope of $\log(\text{error})$ against $\log(\triangle x)$.}
\label{tab:err-t0p0015}
\footnotesize
\setlength{\tabcolsep}{2pt}
\begin{tabular}{llr rrrrr}
\toprule
$t$ & Protocol & Grid pair & $\mathcal{W}_{1,1}$ & $\mathcal{W}_{1,2}$ & $\mathcal{E}_{\mathrm{mean}}$ & $\mathcal{E}_{\mathrm{std}}$ & $\mathcal{E}_{\mathrm{strong}}$ \\
\midrule
\multirow{8}{*}{$0.0015$} & \multirow{4}{*}{\textit{Cauchy}} & 500/1000 & $3.695 \times 10^{-2}$ & $4.183 \times 10^{-2}$ & $3.569 \times 10^{-2}$ & $7.524 \times 10^{-3}$ & $6.179 \times 10^{-2}$ \\
 &  & 1000/2000 & $3.022 \times 10^{-2}$ & $5.441 \times 10^{-2}$ & $2.806 \times 10^{-2}$ & $9.442 \times 10^{-3}$ & $5.381 \times 10^{-2}$ \\
 &  & 2000/4000 & $2.141 \times 10^{-2}$ & $3.338 \times 10^{-2}$ & $1.993 \times 10^{-2}$ & $7.309 \times 10^{-3}$ & $4.569 \times 10^{-2}$ \\
 & & $\mathfrak{r}$ (rate) & \textbf{0.394} & \textbf{0.163} & \textbf{0.420} & \textbf{0.021} & \textbf{0.218} \\
\cmidrule(l){2-8}
 & \multirow{4}{*}{\textit{reference}} & 500/4000 & $6.128 \times 10^{-2}$ & $7.934 \times 10^{-2}$ & $5.748 \times 10^{-2}$ & $1.772 \times 10^{-2}$ & $1.029 \times 10^{-1}$ \\
 &  & 1000/4000 & $4.113 \times 10^{-2}$ & $8.340 \times 10^{-2}$ & $3.789 \times 10^{-2}$ & $1.379 \times 10^{-2}$ & $7.232 \times 10^{-2}$ \\
 &  & 2000/4000 & $2.141 \times 10^{-2}$ & $3.338 \times 10^{-2}$ & $1.993 \times 10^{-2}$ & $7.309 \times 10^{-3}$ & $4.569 \times 10^{-2}$ \\
 & & $\mathfrak{r}$ (rate) & \textbf{0.758} & \textbf{0.624} & \textbf{0.764} & \textbf{0.639} & \textbf{0.586} \\
\bottomrule
\end{tabular}
\end{table}

\begin{table}[ht!]
\centering
\caption{Per-resolution $L^1(\mathcal{D})$ error values and fitted convergence rates $\mathfrak{r}$ of \eqref{eq:ratefit} at $t=0.0020$, with observables in columns and protocols in row blocks. All fields are conservatively block-averaged to the common coarsest grid of $500\times400$ cells before evaluation, and $M=1000$ samples are used at every resolution. Within each protocol block the first three rows are the per-pair errors (grid pair given as coarse/fine streamwise cell count $N_x$) and the bold $\mathfrak{r}$ row is the least-squares slope of $\log(\text{error})$ against $\log(\triangle x)$.}
\label{tab:err-t0p002}
\footnotesize
\setlength{\tabcolsep}{2pt}
\begin{tabular}{llr rrrrr}
\toprule
$t$ & Protocol & Grid pair & $\mathcal{W}_{1,1}$ & $\mathcal{W}_{1,2}$ & $\mathcal{E}_{\mathrm{mean}}$ & $\mathcal{E}_{\mathrm{std}}$ & $\mathcal{E}_{\mathrm{strong}}$ \\
\midrule
\multirow{8}{*}{$0.0020$} & \multirow{4}{*}{\textit{Cauchy}} & 500/1000 & $6.059 \times 10^{-2}$ & $5.249 \times 10^{-2}$ & $5.845 \times 10^{-2}$ & $1.661 \times 10^{-2}$ & $9.847 \times 10^{-2}$ \\
 &  & 1000/2000 & $4.976 \times 10^{-2}$ & $8.743 \times 10^{-2}$ & $4.518 \times 10^{-2}$ & $1.892 \times 10^{-2}$ & $8.690 \times 10^{-2}$ \\
 &  & 2000/4000 & $3.927 \times 10^{-2}$ & $3.466 \times 10^{-2}$ & $3.663 \times 10^{-2}$ & $1.210 \times 10^{-2}$ & $8.459 \times 10^{-2}$ \\
 & & $\mathfrak{r}$ (rate) & \textbf{0.313} & \textbf{0.299} & \textbf{0.337} & \textbf{0.228} & \textbf{0.110} \\
\cmidrule(l){2-8}
 & \multirow{4}{*}{\textit{reference}} & 500/4000 & $1.016 \times 10^{-1}$ & $1.187 \times 10^{-1}$ & $9.566 \times 10^{-2}$ & $3.546 \times 10^{-2}$ & $1.636 \times 10^{-1}$ \\
 &  & 1000/4000 & $7.085 \times 10^{-2}$ & $1.007 \times 10^{-1}$ & $6.519 \times 10^{-2}$ & $2.646 \times 10^{-2}$ & $1.199 \times 10^{-1}$ \\
 &  & 2000/4000 & $3.927 \times 10^{-2}$ & $3.466 \times 10^{-2}$ & $3.663 \times 10^{-2}$ & $1.210 \times 10^{-2}$ & $8.459 \times 10^{-2}$ \\
 & & $\mathfrak{r}$ (rate) & \textbf{0.686} & \textbf{0.888} & \textbf{0.692} & \textbf{0.775} & \textbf{0.476} \\
\bottomrule
\end{tabular}
\end{table}

\begin{table}[ht!]
\centering
\caption{Per-resolution $L^1(\mathcal{D})$ error values and fitted convergence rates $\mathfrak{r}$ of \eqref{eq:ratefit} at $t=0.0025$, with observables in columns and protocols in row blocks. All fields are conservatively block-averaged to the common coarsest grid of $500\times400$ cells before evaluation, and $M=1000$ samples are used at every resolution. Within each protocol block the first three rows are the per-pair errors (grid pair given as coarse/fine streamwise cell count $N_x$) and the bold $\mathfrak{r}$ row is the least-squares slope of $\log(\text{error})$ against $\log(\triangle x)$.}
\label{tab:err-t0p0025}
\footnotesize
\setlength{\tabcolsep}{2pt}
\begin{tabular}{llr rrrrr}
\toprule
$t$ & Protocol & Grid pair & $\mathcal{W}_{1,1}$ & $\mathcal{W}_{1,2}$ & $\mathcal{E}_{\mathrm{mean}}$ & $\mathcal{E}_{\mathrm{std}}$ & $\mathcal{E}_{\mathrm{strong}}$ \\
\midrule
\multirow{8}{*}{$0.0025$} & \multirow{4}{*}{\textit{Cauchy}} & 500/1000 & $9.703 \times 10^{-2}$ & $6.152 \times 10^{-2}$ & $9.311 \times 10^{-2}$ & $2.884 \times 10^{-2}$ & $1.484 \times 10^{-1}$ \\
 &  & 1000/2000 & $7.553 \times 10^{-2}$ & $1.023 \times 10^{-1}$ & $6.886 \times 10^{-2}$ & $3.431 \times 10^{-2}$ & $1.286 \times 10^{-1}$ \\
 &  & 2000/4000 & $6.083 \times 10^{-2}$ & $3.453 \times 10^{-2}$ & $5.710 \times 10^{-2}$ & $2.123 \times 10^{-2}$ & $1.301 \times 10^{-1}$ \\
 & & $\mathfrak{r}$ (rate) & \textbf{0.337} & \textbf{0.417} & \textbf{0.353} & \textbf{0.221} & \textbf{0.095} \\
\cmidrule(l){2-8}
 & \multirow{4}{*}{\textit{reference}} & 500/4000 & $1.471 \times 10^{-1}$ & $1.071 \times 10^{-1}$ & $1.376 \times 10^{-1}$ & $6.188 \times 10^{-2}$ & $2.242 \times 10^{-1}$ \\
 &  & 1000/4000 & $1.091 \times 10^{-1}$ & $1.109 \times 10^{-1}$ & $1.004 \times 10^{-1}$ & $4.764 \times 10^{-2}$ & $1.777 \times 10^{-1}$ \\
 &  & 2000/4000 & $6.083 \times 10^{-2}$ & $3.453 \times 10^{-2}$ & $5.710 \times 10^{-2}$ & $2.123 \times 10^{-2}$ & $1.301 \times 10^{-1}$ \\
 & & $\mathfrak{r}$ (rate) & \textbf{0.637} & \textbf{0.817} & \textbf{0.634} & \textbf{0.772} & \textbf{0.393} \\
\bottomrule
\end{tabular}
\end{table}

\begin{table}[ht!]
\centering
\caption{Per-resolution $L^1(\mathcal{D})$ error values and fitted convergence rates $\mathfrak{r}$ of \eqref{eq:ratefit} at $t=0.0030$, with observables in columns and protocols in row blocks. All fields are conservatively block-averaged to the common coarsest grid of $500\times400$ cells before evaluation, and $M=1000$ samples are used at every resolution. Within each protocol block the first three rows are the per-pair errors (grid pair given as coarse/fine streamwise cell count $N_x$) and the bold $\mathfrak{r}$ row is the least-squares slope of $\log(\text{error})$ against $\log(\triangle x)$.}
\label{tab:err-t0p003}
\footnotesize
\setlength{\tabcolsep}{2pt}
\begin{tabular}{llr rrrrr}
\toprule
$t$ & Protocol & Grid pair & $\mathcal{W}_{1,1}$ & $\mathcal{W}_{1,2}$ & $\mathcal{E}_{\mathrm{mean}}$ & $\mathcal{E}_{\mathrm{std}}$ & $\mathcal{E}_{\mathrm{strong}}$ \\
\midrule
\multirow{8}{*}{$0.0030$} & \multirow{4}{*}{\textit{Cauchy}} & 500/1000 & $1.233 \times 10^{-1}$ & $1.578 \times 10^{-1}$ & $1.182 \times 10^{-1}$ & $3.966 \times 10^{-2}$ & $1.802 \times 10^{-1}$ \\
 &  & 1000/2000 & $1.097 \times 10^{-1}$ & $9.063 \times 10^{-2}$ & $1.024 \times 10^{-1}$ & $5.307 \times 10^{-2}$ & $1.791 \times 10^{-1}$ \\
 &  & 2000/4000 & $7.736 \times 10^{-2}$ & $3.671 \times 10^{-2}$ & $7.326 \times 10^{-2}$ & $3.110 \times 10^{-2}$ & $1.715 \times 10^{-1}$ \\
 & & $\mathfrak{r}$ (rate) & \textbf{0.336} & \textbf{1.052} & \textbf{0.345} & \textbf{0.175} & \textbf{0.036} \\
\cmidrule(l){2-8}
 & \multirow{4}{*}{\textit{reference}} & 500/4000 & $1.918 \times 10^{-1}$ & $1.230 \times 10^{-1}$ & $1.784 \times 10^{-1}$ & $8.726 \times 10^{-2}$ & $2.764 \times 10^{-1}$ \\
 &  & 1000/4000 & $1.525 \times 10^{-1}$ & $1.152 \times 10^{-1}$ & $1.432 \times 10^{-1}$ & $6.919 \times 10^{-2}$ & $2.360 \times 10^{-1}$ \\
 &  & 2000/4000 & $7.736 \times 10^{-2}$ & $3.671 \times 10^{-2}$ & $7.326 \times 10^{-2}$ & $3.110 \times 10^{-2}$ & $1.715 \times 10^{-1}$ \\
 & & $\mathfrak{r}$ (rate) & \textbf{0.655} & \textbf{0.872} & \textbf{0.642} & \textbf{0.744} & \textbf{0.344} \\
\bottomrule
\end{tabular}
\end{table}

\begin{table}[ht!]
\centering
\caption{Per-resolution $L^1(\mathcal{D})$ error values and fitted convergence rates $\mathfrak{r}$ of \eqref{eq:ratefit} at $t=0.0035$, with observables in columns and protocols in row blocks. All fields are conservatively block-averaged to the common coarsest grid of $500\times400$ cells before evaluation, and $M=1000$ samples are used at every resolution. Within each protocol block the first three rows are the per-pair errors (grid pair given as coarse/fine streamwise cell count $N_x$) and the bold $\mathfrak{r}$ row is the least-squares slope of $\log(\text{error})$ against $\log(\triangle x)$.}
\label{tab:err-t0p0035}
\footnotesize
\setlength{\tabcolsep}{2pt}
\begin{tabular}{llr rrrrr}
\toprule
$t$ & Protocol & Grid pair & $\mathcal{W}_{1,1}$ & $\mathcal{W}_{1,2}$ & $\mathcal{E}_{\mathrm{mean}}$ & $\mathcal{E}_{\mathrm{std}}$ & $\mathcal{E}_{\mathrm{strong}}$ \\
\midrule
\multirow{8}{*}{$0.0035$} & \multirow{4}{*}{\textit{Cauchy}} & 500/1000 & $1.497 \times 10^{-1}$ & $2.403 \times 10^{-1}$ & $1.423 \times 10^{-1}$ & $5.203 \times 10^{-2}$ & $2.120 \times 10^{-1}$ \\
 &  & 1000/2000 & $1.410 \times 10^{-1}$ & $1.837 \times 10^{-1}$ & $1.311 \times 10^{-1}$ & $7.067 \times 10^{-2}$ & $2.232 \times 10^{-1}$ \\
 &  & 2000/4000 & $9.319 \times 10^{-2}$ & $4.553 \times 10^{-2}$ & $8.745 \times 10^{-2}$ & $4.835 \times 10^{-2}$ & $2.044 \times 10^{-1}$ \\
 & & $\mathfrak{r}$ (rate) & \textbf{0.342} & \textbf{1.200} & \textbf{0.351} & \textbf{0.053} & \textbf{0.026} \\
\cmidrule(l){2-8}
 & \multirow{4}{*}{\textit{reference}} & 500/4000 & $2.373 \times 10^{-1}$ & $1.896 \times 10^{-1}$ & $2.181 \times 10^{-1}$ & $1.124 \times 10^{-1}$ & $3.294 \times 10^{-1}$ \\
 &  & 1000/4000 & $1.940 \times 10^{-1}$ & $2.143 \times 10^{-1}$ & $1.805 \times 10^{-1}$ & $9.802 \times 10^{-2}$ & $2.860 \times 10^{-1}$ \\
 &  & 2000/4000 & $9.319 \times 10^{-2}$ & $4.553 \times 10^{-2}$ & $8.745 \times 10^{-2}$ & $4.835 \times 10^{-2}$ & $2.044 \times 10^{-1}$ \\
 & & $\mathfrak{r}$ (rate) & \textbf{0.674} & \textbf{1.029} & \textbf{0.659} & \textbf{0.609} & \textbf{0.344} \\
\bottomrule
\end{tabular}
\end{table}

The convergence rates of all five observables under both comparison operators are collected in Figure~\ref{fig:protocol_rates}, and the underlying per-resolution error values, together with the fitted rates~$\mathfrak{r}$, are tabulated at each time step in Tables~\ref{tab:err-t0p0005}--\ref{tab:err-t0p0035}. The statistical estimators employed here follow those of~\cite{simonis2026}.
\newsavebox{\mybox}

\subsection{Temporal trends of the convergence rates}\label{subsec:ratetrends}

While Tables~\ref{tab:err-t0p0005}--\ref{tab:err-t0p0035} report the convergence rate $\mathfrak{r}_{\mathfrak{O}}(t)$ at each checkpoint separately, the time series of these rates itself carries information: whether an observable's convergence is stationary, improving, or deteriorating as the flow develops from the laminar transient into fully developed turbulence. To quantify this, Table~\ref{tab:rate-trends} reports a linear least-squares fit $\mathfrak{r}_{\mathfrak{O}}(t)\approx a\,t+b$ through the seven per-checkpoint rates of Figure~\ref{fig:protocol_rates}, for each observable and each comparison protocol. Alongside the raw slope $a$ and intercept $b$, the table lists the total drift $\triangle\mathfrak{r}=a\cdot0.003$ accumulated across the analysis window $t\in[0.0005,0.0035]$, which expresses the same trend in directly readable rate units, and the coefficient of determination $R^2$, which separates genuine temporal trends from scatter about a flat mean. The final column extrapolates the fit to the last checkpoint.

The fits divide the ten series into three clearly distinguishable classes. The first class comprises the strong pathwise error under both protocols, which exhibits the steepest and cleanest decline in the entire dataset ($\triangle\mathfrak{r}=-0.338$ with $R^2=0.940$ for the \textit{Cauchy} comparison, $\triangle\mathfrak{r}=-0.441$ with $R^2=0.954$ for the \textit{reference} comparison). The near-perfect linearity of this decay is itself a physical statement: the loss of pathwise convergence is not an abrupt event but a steady erosion that tracks the gradual development of the shear-layer turbulence, with the fitted \textit{Cauchy} rate reaching $-0.011\approx0$ at $t_{\mathrm{end}}$. The extrapolated fit thus confirms that individual realizations have ceased to converge by the end of the run.

The second class comprises the two-point Wasserstein distance, the only observable whose rate \emph{rises}, and it does so strongly and reliably under both protocols ($\triangle\mathfrak{r}=+1.888$, $R^2=0.828$ \textit{Cauchy}; $\triangle\mathfrak{r}=+0.591$, $R^2=0.848$ \textit{reference}). The opposite sign relative to the strong error is the trend-level expression of the crossover discussed in Section~\ref{subsec:discussion}: as coherent large-scale correlations (roll-up, bow shock) become established, consecutive grids agree on the two-point structure increasingly well, and the \textit{Cauchy} rate not only rises but overtakes the \textit{reference} rate, its fitted end value of $1.265$ exceeding the reference value of $1.038$. The large positive $\triangle\mathfrak{r}$ of the \textit{Cauchy} series should nonetheless be read with the early-time caveat in mind: the fit window includes the ill-conditioned transient value $-1.091$ at $t=0.0005$, which steepens the fitted slope beyond what the developed-turbulence regime alone would give.

The third class comprises the one-point statistics. The \textit{reference} rates of $\mathcal{W}_{1,1}$ and the mean drift down only mildly ($\triangle\mathfrak{r}=-0.126$ and $-0.147$, both with good fits, $R^2\approx0.72$--$0.79$), settling toward fitted end values of $0.640$ and $0.628$; their \textit{Cauchy} counterparts are statistically flat ($\triangle\mathfrak{r}=-0.034$ and $-0.043$ at $R^2=0.205$ and $0.308$), i.e.\ stationary near $0.33$--$0.35$ throughout. The mild reference-side decline is consistent with the flow becoming statistically harder as ever more of the domain is occupied by the turbulent head, without any indication of breakdown. Finally, the standard deviation is the deliberate outlier of the table: under both protocols its $R^2$ is below $0.02$, meaning the fitted line explains essentially none of the variation, and the tabulated drifts ($-0.064$ and $+0.041$) are noise rather than trend. This is the trend-level counterpart of the second moment's erratic per-checkpoint behavior discussed in Section~\ref{subsec:discussion}: the std rate fluctuates about its mean level (roughly $0.14$ \textit{Cauchy}, $0.68$ \textit{reference}) without any resolvable temporal drift, and we report its fit only for completeness.

In summary, the trend analysis compresses the qualitative picture of Figure~\ref{fig:protocol_rates} into a small set of reliable numbers: a steadily collapsing pathwise rate, a steadily improving two-point statistical rate, stationary-to-mildly-declining one-point rates, and a second moment whose temporal behavior is dominated by noise. The opposing signs of the two reliable extremes, i.e., $\triangle\mathfrak{r}<0$ for the strong error versus $\triangle\mathfrak{r}>0$ for the two-point distance, restate the central dichotomy of this work at the level of temporal trends: as the turbulence develops, deterministic convergence degrades precisely while statistical convergence consolidates.

\begin{table}[ht!]
\centering
\caption{Linear least-squares fit $\mathfrak{r}_{\mathfrak{O}}(t)\approx a\,t+b$ of the time-dependent convergence rate \eqref{eq:ratefit} for each observable under each comparison protocol, over $t\in[0.0005,0.0035]$. The column $\triangle\mathfrak{r}$ is the total change of the fitted rate across the run ($a$ times the window length $0.003$), i.e.\ the drift from the first to the last time; $R^2$ is the goodness of fit; $a$ and $b$ are the slope (rate per unit time) and intercept; and $\mathfrak{r}(t_{\mathrm{end}})$ is the fitted rate at $t_{\mathrm{end}}=0.0035$. Low $R^2$ (notably $\mathcal{E}_{\mathrm{std}}$ under both protocols and the Cauchy $\mathcal{W}_{1,1}$) indicates a rate that is effectively flat or noise-dominated rather than genuinely trending, so the corresponding $\triangle\mathfrak{r}$ should not be over-interpreted.}
\label{tab:rate-trends}
\footnotesize
\setlength{\tabcolsep}{4pt}
\begin{tabular}{ll rr rrr}
\toprule
Observable & Protocol & $\triangle\mathfrak{r}$ & $R^2$ & $a$ (slope) & $b$ (intercept) & $\mathfrak{r}(t_{\mathrm{end}})$ \\
\midrule
$\mathcal{W}_{1,1}$ & \textit{Cauchy} & $-0.034$ & $0.205$ & $-1.136\times10^{1}$ & $0.371$ & $0.331$ \\
 & \textit{reference} & $-0.126$ & $0.722$ & $-4.207\times10^{1}$ & $0.787$ & $0.640$ \\
\addlinespace
$\mathcal{W}_{1,2}$ & \textit{Cauchy} & $+1.888$ & $0.828$ & $6.292\times10^{2}$ & $-0.937$ & $1.265$ \\
 & \textit{reference} & $+0.591$ & $0.848$ & $1.969\times10^{2}$ & $0.349$ & $1.038$ \\
\addlinespace
$\mathcal{E}_{\mathrm{mean}}$ & \textit{Cauchy} & $-0.043$ & $0.308$ & $-1.450\times10^{1}$ & $0.394$ & $0.344$ \\
 & \textit{reference} & $-0.147$ & $0.790$ & $-4.893\times10^{1}$ & $0.800$ & $0.628$ \\
\addlinespace
$\mathcal{E}_{\mathrm{std}}$ & \textit{Cauchy} & $-0.064$ & $0.018$ & $-2.121\times10^{1}$ & $0.185$ & $0.111$ \\
 & \textit{reference} & $+0.041$ & $0.015$ & $1.350\times10^{1}$ & $0.656$ & $0.703$ \\
\addlinespace
$\mathcal{E}_{\mathrm{strong}}$ & \textit{Cauchy} & $-0.338$ & $0.940$ & $-1.126\times10^{2}$ & $0.383$ & $-0.011$ \\
 & \textit{reference} & $-0.441$ & $0.954$ & $-1.471\times10^{2}$ & $0.802$ & $0.287$ \\
\bottomrule
\end{tabular}
\end{table}

\section{Conclusions}\label{sec:conclusions}
We have computed the statistical solution of a Mach~2000 astrophysical jet on a sequence of grids reaching $4000\times3200$ cells ($12.8$ million), with $M=1000$ Monte Carlo samples at every resolution.
By defining the numerical solution as the pushforward~\eqref{eq:pushforward} of an initial probability measure through the VLBM operator, we mitigated the traditional instability of deterministic schemes in near-vacuum regimes; every quantity reported is a functional of the empirical measure~\eqref{eq:empirical} rather than of a single trajectory.
The central empirical finding is a sharp dichotomy.
The sample-wise strong $L^1$ error~\eqref{eq:strongError} ceases to converge under refinement (Cauchy rate $\approx0.03$ at $t=0.0035$), certifying the loss of pathwise uniqueness, while the statistical observables converge, with reference rates of $\approx0.67$ in the one-point Wasserstein distance~\eqref{eq:w1disc}, $\approx0.66$ in the mean~\eqref{eq:mean}, $\approx0.61$ in the standard deviation~\eqref{eq:std}, and $\approx1.03$ in the two-point Wasserstein distance~\eqref{eq:w2}.
The spatially resolved analysis along the jet axis locates the origin of this dichotomy: the one-point marginals~\eqref{eq:onepointlaw} are numerically Dirac in the undisturbed core, skewed in the Kelvin--Helmholtz-sheared regions and at the turbulent head, and, at the leading fringe, intermittent two-state mixtures~\eqref{eq:mixture} governed by the random front position.
These laws remain stable under refinement even as the individual realizations sampling them decorrelate.
This provides empirical evidence for the statistical stability of chaotic high-Mach flows and suggests a practical route to uncertainty quantification in extreme compressible astrophysics: pointwise deterministic uniqueness is lost, but convergence to a statistical limit is preserved.

The weak--strong uniqueness principle is conditional on the existence of a regular strong solution, and no such solution exists here: the inlet data are discontinuous and shocks form immediately.
The analytical principle therefore makes no assertion about this flow, and we claim neither uniqueness of the limit measure nor convergence of the scheme in an analytical sense.
What the computations document is a numerical counterpart of the weak--strong dichotomy: the VLBM selects, through its numerical dissipation, a particular sequence of approximations whose pathwise limit does not exist but whose statistical limit is stable and non-degenerate across resolutions.
A different discretization could in principle select a different statistical solution; quantifying this scheme dependence is a natural next step. 
Beyond the Wasserstein observables used here, the ensembles also lend themselves to other statistical descriptions of turbulence. 
For instance, random matrix statistics of snapshot ensembles place both incompressible and compressible turbulence in the universality class of quantum chaotic systems~\cite{levi2023}, and it would be of interest to test whether this structure persists in the extreme compressible regime of the present jet.

In observational astrophysics, the visible structures of extreme jets (such as emission knots and bow shocks) are directly governed by the turbulent dissipation of kinetic energy at internal shock fronts \cite{rees1978, blandford2019}.
Our results demonstrate that while the exact deterministic location of these shocks is gradually lost upon grid refinement in single sample simulations, the statistical footprint of the turbulence, and thus the macroscopic structure of the observable emission, is recovered under the Wasserstein metric.

\section*{Acknowledgments}
The work of S.\ Simonis was supported by the PRIME programme of the German Academic Exchange Service (DAAD), with funding from the Federal Ministry of Research, Technology and Space.
The initiation of this work at UZH was supported by a KHYS ConYS grant at KIT in 2024. 
S.\ Simonis thanks M.J.\ Krause, M.\ Frank and S.\ Mishra for providing the excellent research infrastructure at KIT and ETH Zurich.
S.\ Simonis gratefully acknowledges the computing time provided on the high-performance computer HoreKa by the National High-Performance Computing Center at KIT (NHR@KIT).
This center is jointly supported with funding from the Federal Ministry of Research, Technology and Space and the Ministry of Science, Research and the Arts of Baden-Württemberg, as part of the National High-Performance Computing (NHR) joint funding program (\href{https://www.nhr-verein.de/en/our-partners}{https://www.nhr-verein.de/en/our-partners}).
HoreKa is partly funded by the German Research Foundation (DFG).

\section*{Author contribution statement}
S.\ Simonis:
    Conceptualization,
    Methodology,
    Software,
    Validation,
    Formal Analysis,
    Investigation,
    Resources,
    Data Curation,
    Writing - Original Draft,
    Writing - Review {\&} Editing,
    Visualization,
    Project administration,
    Funding Acquisition;
G.\ Wissocq:
    Methodology,
    Writing - Review {\&} Editing.
All authors read and agreed to the final version of the manuscript.

\section*{Code and data availability}
The source code is available as the open-source library acropora~\cite{acropora} at \href{https://github.com/stephan-911/acropora}{https://github.com/stephan-911/acropora} under the BSD-2-Clause-Patent license.
The Monte Carlo ensembles restricted to the $500\times400$ grid, the convergence and centerline results, and the random inflow coefficients of every sample are archived on Zenodo~\cite{acroporadata} under the CC BY 4.0 license. 
The full-resolution ensembles (about 1.1\,TB) can be regenerated with~\cite{acropora}.

\section*{Declarations}
The authors disclose that generative AI assistants (Google Gemini and Anthropic Claude) were used during the preparation of this work to assist with code development, data analysis, and text formatting.
All AI-assisted code and analyses were tested, reviewed, and verified by the authors, and all text was reviewed and edited by the authors.
The authors assume responsibility for all content.
The authors declare that they have no known competing financial interests or personal relationships that could have appeared to influence the work reported in this paper.

\end{document}